\documentclass[a4paper,11pt]{article}
\usepackage[a4paper,text={16cm,22.5cm},centering]{geometry}
\usepackage{amsmath,amssymb,mathtools}
\usepackage{mathrsfs,bbm}
\usepackage[utf8]{inputenc}
\usepackage[TS1,T1]{fontenc}
\usepackage{lmodern}
\usepackage{microtype}
\usepackage{tikz}
\usepackage{authblk}

\usetikzlibrary{arrows} 
\usetikzlibrary{patterns} 
\usetikzlibrary{decorations.markings} 
\usetikzlibrary{decorations.pathmorphing}

\definecolor{vert}{rgb}{0.1,0.7,0.15}

\tikzset{hatch distance/.store in=\hatchdistance, hatch distance=10pt, hatch thickness/.store in=\hatchthickness, hatch thickness=2pt}

\makeatletter 
\pgfdeclarepatternformonly[\hatchdistance,\hatchthickness]{flexible hatch} {\pgfqpoint{-1pt}{-1pt}} {\pgfqpoint{\hatchdistance}{\hatchdistance}} {\pgfpoint{\hatchdistance-1pt}{\hatchdistance-1pt}} {\pgfsetcolor{\tikz@pattern@color} 
\pgfsetlinewidth{\hatchthickness} \pgfpathmoveto{\pgfqpoint{-1pt}{-1pt}} \pgfpathlineto{\pgfqpoint{\hatchdistance}{\hatchdistance}} \pgfusepath{stroke}} 
\makeatother

\rmfamily
\DeclareFontShape{T1}{lmr}{b}{sc}{<->ssub*cmr/bx/sc}{}
\DeclareFontShape{T1}{lmr}{bx}{sc}{<->ssub*cmr/bx/sc}{}

\DeclareMathOperator*{\esup}{ess~sup}
\DeclareMathOperator*{\argmin}{arg~min}

\newtheorem{defi}{Definition} 
\newtheorem{prop}[defi]{Proposition} 
\newtheorem{lemm}[defi]{Lemma} 
\newtheorem{theo}[defi]{Theorem} 
\newtheorem{coro}[defi]{Corollary} 
\newcommand{\qed}{\hglue 0pt\hfill$\square$\par}

\title{Homogenization of friction in a 2D linearly elastic contact problem}
\author[$\dagger$]{Patrick \textsc{Ballard}}
\author[$\ddag$]{Flaviana \textsc{Iurlano}}
\affil[$\dagger$]{Sorbonne Université, CNRS, Université de Paris, Institut Jean le Rond d'Alembert}
\affil[$\ddag$]{Sorbonne Université, CNRS, Université de Paris, Laboratoire Jacques-Louis Lions}
\begin{document}

\maketitle

\begin{abstract}
	Contact problems with Coulomb friction in linear elasticity are notoriously difficult,
	and their mathematical analysis is still largely incomplete.  In this paper, a model problem with heterogeneous friction coefficient is considered in two-dimensional elasticity.  For this model problem, an existence and uniqueness result is proved, relying heavily on harmonic analysis.  A complete and rigorous homogenization analysis can be performed in the case of a highly oscillating friction coefficient, being the first result in that direction.  The Coulomb law is found to hold in the limit, and an explicit formula is provided to calculate the effective friction coefficient.  This effective friction coefficient is found to differ from the spatial average, showing an influence of the coupling between friction and elasticity on the homogenized limit.
\end{abstract}

\section{Introduction}

\emph{Rate-independent processes} are a special class of evolution problems that do not possess any internal time scale.  They are appropriate for studying quasistatic evolution in which viscous effects can be neglected.  Some examples are provided by perfect plasticity~\cite{DalMasoPlast}, brittle damage~\cite{FrancfortGarroni}, brittle fracture~\cite{FrancfortMarigo} or phase transitions, which all have received considerable attention in the past decades both in the solid mechanics and applied mathematics communities.  The above examples have in common to possess an underlying variational structure: when introducing a time-discretization, the problem to solve at each time-step is spontaneously a minimization problem.  Thus, the techniques of the calculus of variation (direct method, relaxation, $\Gamma$-convergence,\ldots) have considerably contributed to the understanding of the underlying mathematical structures, enlightening the main qualitative features of the solutions.

However, not all rate-independent processes need to spontaneously exhibit an underlying variational structure.  One example is provided by contact problems with Coulomb friction in linear elasticity whose mathematical analysis was undertaken as early as 1972 by Duvaut \& Lions \cite{DuvautLions}.  More recent examples are the so-called non-associative models of plasticity \cite{BabaFrancfortMora}.  As a consequence of the lack of a `spontaneous' underlying variational structure, the mathematical analysis of these problems turns out to be considerably more involved.

Denoting by \(\mathbf{u}=u_{\rm n}\mathbf{n}+\mathbf{u}_{\rm t}\) and \(\mathbf{t}=t_{\rm n}\mathbf{n}+\mathbf{t}_{\rm t}\) the boundary displacement and the surface traction, and their splitting into normal and tangential parts (here \(\mathbf{n}\) denotes the outward unit normal), contact problems with Coulomb friction amount to consider a boundary condition in linear elasticity in which the Signorini contact condition:
\[
u_{\rm n}\leq \bar{g},\qquad t_{\rm n}\leq 0, \qquad (u_{\rm n}-\bar{g})t_{\rm n}=0,
\]
(here \(\bar{g}\) is a given function representing the initial gap with the obstacle), is complemented with the Coulomb law of dry friction, whose formal pointwise formulation can be synthetically written as:
\begin{equation}
	\label{eq:Coulomb}
	\forall \hat{\mathbf{v}},\qquad  \mathbf{t}_{\rm t}\cdot (\hat{\mathbf{v}} - \dot{\mathbf{u}}_{\rm t}) - \bar{f} t_{\rm n} \bigl(|\hat{\mathbf{v}}| - |\dot{\mathbf{u}}_{\rm t}|\bigr) \geq 0.
\end{equation}
Here, \(\bar{f}\geq 0\) is the given friction coefficient and the dot refers to the time-derivative.  The associated incremental problem (the problem to solve at each time-step in a time discretization) amounts to solve the same problem in which the velocity $\dot{\mathbf{u}}_{\rm t}$ in the Coulomb law is replaced by $\mathbf{u}_{\rm t}$.
Existence of a solution for the incremental problem was first proved in \cite{Jarusek} using a fixed point argument and extended in \cite{EckJarusek} based on the theory of pseudomonotone operators.  The existence of a solution for the evolution problem was first obtained in \cite{Andersson}.  All these results have in common to require that the friction coefficient is lower than a finite critical value and leave open the question of uniqueness for arbitrarily small friction coefficients.

The problem simplifies, though preserving mechanical meaningfulness, in the situation where a steady motion of the obstacle is prescribed, in which case the Coulomb law~\eqref{eq:Coulomb} of dry friction reduces to the linear condition:
\begin{equation}
	\label{eq:SteadyCoulomb}
	\mathbf{t}_{\rm t} = \bar{f}t_{\rm n}\boldsymbol{\tau},
\end{equation}
where \(\boldsymbol{\tau}\) is a given unit tangential vector (the normalized relative tangential velocity).  The corresponding elastostatic contact problem was first studied in \cite{bibiElas} in the case of a homogeneous friction coefficient.  In the case of a bounded body, it was proved using classical results on variational inequalities that there always exists a nonempty range \([0,\bar{f}_{\rm c}[\) so that the problem admits a unique solution for all \(\bar{f}\in [0,\bar{f}_{\rm c}[\).  An example of multiple solution was displayed, showing the possibility of a bifurcation with respect to the friction coefficient taken as a control parameter.  In the case where the body is a two-dimensional homogeneous isotropic half-space with homogeneous friction coefficient, the situation was proved simpler as a unique solution exists whatever the value of the friction coefficient is.

Hence, this background drove us naturally towards the geometry of the two-dimensional elastic half-space to investigate the homogenization of dry friction, that is, to study asymptotically highly oscillating heterogeneous friction coefficients.  Three natural questions arise in this context.
\begin{enumerate}
	\item Is the heterogeneous problem well-posed, that is, does it admit a unique solution?
	\item Is the steady friction law~\eqref{eq:SteadyCoulomb} stable by homogenization, that is, does it hold true in the limit with the heterogeneous friction coefficient \(\bar{f}(x)\) replaced with some constant effective friction coefficient \(\bar{f}_{\rm eff}\)?
	\item If so, what is the constant effective friction coefficient?
\end{enumerate}
In this paper, we bring a complete and rigorous answer to these three questions in the case of the simplest geometry, that is, that of a moving convex rigid indentor of finite width at the surface of the two-dimensional isotropic elastic half-space.  In particular, we prove that the frictional contact problem in the case of a highly oscillating periodic friction coefficient \(\bar{f}(x)\) has a unique solution, and that a constant effective friction coefficient \(\bar{f}_{\rm eff}\) appears in the limit, given by the formula:
\[
\bar{f}_{\rm eff} = {\textstyle\frac{2(1-\nu)}{1-2\nu}} \tan \biggl\langle \arctan \Bigl({\textstyle\frac{1-2\nu}{2(1-\nu)}}\,\bar{f}\Bigr)\biggr\rangle,
\]
where \(\nu\) is the Poisson ratio and \(\langle \cdot \rangle\) stands for the spatial average.  Hence, the effective friction coefficient depends on the elastic moduli and is different from the spatial average of the microscopic friction coefficient, except in the limit case of incompressibility (\(\nu \rightarrow 1/2\)).

The geometry of the problem to which this paper is restricted (a moving rigid indentor of finite width at the surface of the two-dimensional isotropic elastic half-space) truly encompasses all the couplings between unilateral contact, dry friction and elasticity.  Therefore, we expect that this result has a wide range of validity for steady sliding frictional contact problems in \emph{two-dimensional} isotropic homogeneous linear elasticity.  As an important difficulty in homogenization is often to identify the limit problem, we believe that it is the merit of the special geometry which is studied in this paper, to lead to the identification of the good candidate.

This whole investigation was at first sight deeply obstructed by a few obstacles, which have been overcome in the following way.
\begin{itemize}
	\item Although the problem with homogeneous friction is monotone and has therefore a unique solution for an indentor of arbitrary shape, the heterogeneous problem with highly oscillating periodic friction is never monotone.  However, we have been able to recover some monotonicity and solve it uniquely for an arbitrary \emph{convex} rigid indentor.
	\item The corresponding contact problem is non-variational.  Hence, one cannot expect to rely on $\Gamma$-convergence to study the homogenized limit.  The main problem is to handle the product of the normal component of the surface traction with the heterogeneous friction coefficient, for both of which only weak convergences are obtained.  This difficulty is overcome by exploiting the differential structure of the problem in an original, broad compensated compactness argument.
\end{itemize}

We first analyze a toy model (the flat indentor), for which an explicit solution can be obtained.  Drawing inspiration from this, we are able to uncover monotonicity in the general case of an arbitrary convex indentor, to solve uniquely the corresponding contact problem and prove the homogenization result.

The structure of the paper is the following.  Section~\ref{sec:statementResults} is devoted to the statement of our main results in the case of heterogeneous friction and a rigid flat indentor (Section~\ref{sec:statementFlat}), of homogeneous friction and arbitrary indentor (Section~\ref{sec:constantFriction}), and of heterogeneous friction and arbitrary indentor (existence results, Section~\ref{sec:statementHetero}) or convex indentor (uniqueness, regularity, and homogenization results, Section~\ref{sec:statementConvex}).  Proofs of the results of Sections \ref{sec:statementFlat}--\ref{sec:statementConvex} are presented respectively in Sections \ref{sec:proofFlat}--\ref{sec:proofConvexHetero}.  Finally, Appendices A--D contain some classical results concerning respectively the Hilbert transform, the spaces $H^{1/2}$ and $H^{-1/2}$, the Poisson integral, and pseudomonotone variational inequalities.

\section{Position of problem and statement of results}

\label{sec:statementResults}

\subsection{The formal problem}

We consider an isotropic homogeneous linearly elastic two-dimen\-sional half-space defined by \(y<0\) (\(y\) is the space variable along the direction perpendicular to the boundary).  The Poisson ratio is denoted by \(\nu\in \left]-1,1/2\right[\) and the force unit is chosen so that the Young modulus \(E=1\).  We denote by \(x\) the space variable along the boundary, by \(\mathbf{t}(x)\) the surface traction distribution on the boundary and the (outward) normal and tangential components will be addressed as \(t_\textrm{n}(x)\) and \(t_{\rm t}(x)\) (hence \(t_\textrm{n}=t_y\) and \(t_\textrm{t}=t_x\)).  The half-space is assumed to be free of body forces.  Use of Fourier transform with respect to \(x\) provides explicit knowledge of the fundamental solution \(\mathbf{\bar{u}}\) of the classical Neumann problem in linear elasticity (with no condition at infinity).  It is obtained up to an affine function involving four arbitrary constants.  Two of these constants (corresponding to an overall rotation and a component of stress at infinity) can be fixed by adding the following condition at infinity:
\[
\mathbf{\bar{u}}(x,y) = O\bigl(\log(x^2+y^2)\bigr), \quad \text{ as }\; x^2+y^2 \rightarrow \infty,
\]
where \(\mathbf{\bar{u}}\) denotes the displacement field.  In that case, the stress field goes to zero at infinity.  However, the displacement is generally infinite at infinity and there is no mean to fix the remaining two constants which correspond to an arbitrary translation.  Setting:
\[
\mathbf{u} := \frac{\mathbf{\bar{u}}}{2(1-\nu^2)},\qquad \text{and}\qquad \gamma := \frac{1-2\nu}{2(1-\nu)}\in \left]0,3/4\right[, 
\]
the boundary displacement \(\mathbf{u}(x,0)\) resulting from a prescribed surface traction distribution \(\mathbf{t}(x)\) with compact support, is given (see, for example,~\cite[theorem 1]{bibiJiri} for a detailed proof) by: 
\[
\begin{aligned}
	u_\textrm{n}' & = - \frac{1}{\pi}\,\text{pv}\frac{1}{x} * t_\textrm{n} - \gamma\,t_{\rm t},\\[1.2ex]
	u_{\rm t}' & = - \frac{1}{\pi}\,\text{pv}\frac{1}{x} * t_\textrm{t} + \gamma\,t_{\rm n}, 
\end{aligned}
\qquad\qquad \text{ in }\mathbb{R},
\]
where \(\text{pv}\,1/x\) denotes the distributional derivative of \(\log|x|\) and \(*\) is the convolution product with respect to \(x\).  The acronym `pv' stands for `(Cauchy) principal value' and the mapping:
\[
t \mapsto  - \frac{1}{\pi}\,\text{pv}\frac{1}{x} * t
\]
is known as the Hilbert transform\footnote{the multiplicative constant \(1/\pi\) is there to ensure that the Hilbert transform is an isometry of \(L^2(\mathbb{R})\), whose inverse mapping is minus itself.}.  In the case where \(t\) is a function, it is an example of what is sometimes called a `singular integral' or a `Cauchy principal value integral'. 

The above expression of \((u_\textrm{n},u_\textrm{t})\) in terms of \((t_\textrm{n},t_\textrm{t})\) is nothing but the Neumann-to-Dirichlet operator of the isotropic homogeneous elastic bidimensional half-space.  Each component of the boundary displacement is obtained up to an arbitrary additive constant, which is to  be interpreted as an overall rigid translation.

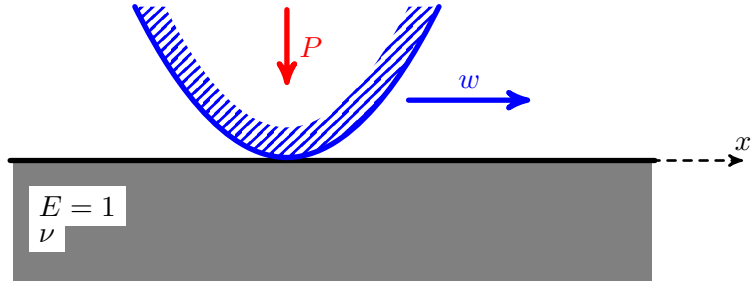
\begin{figure}[ht]
	\begin{center}
		\begin{tikzpicture}[line cap=round,scale=0.4]
			\draw[draw=none,fill=gray] (-12,-4) rectangle (9,0); 
			\draw[line width=2pt] (-12.1,0) -- (9.1,0); 
			\draw[->,>=stealth',line width=1pt,dashed] (-12.1,0) -- (12,0); 
			\path[pattern=flexible hatch,hatch distance=5pt,hatch thickness=1pt,pattern color=blue] (-8,5.1) parabola bend (-3,0.1) (2,5.1) -- (1,5.1) parabola bend (-3,1.1) (-7,5.1) -- cycle; 
			\draw[line width=2pt,blue] (-8,5.1) parabola bend (-3,0.1) (2,5.1); 
			\draw (-11.5,-1.5) node[fill=white,right] {\(E=1\)}; 
			\draw(-11.5,-2.5) node[fill=white,right] {\(\nu\)}; 
			\draw[->,>=stealth',line width=2pt,red] (-3,5) -- (-3,2.5) node[right,pos=0.5,red] {\(P\)}; 
			\draw[->,>=stealth',line width=2pt,blue] (1,2) -- (5,2) node[above,pos=0.5,blue] {\(w\)}; 
			\draw (12,0) node[above] {\(x\)}; 
		\end{tikzpicture}
		\caption{Geometry of the problem.} 
		\label{fig:GenHS2D}
	\end{center}
\end{figure}

We consider a rigid indentor which moves along the boundary of the half-space, at a constant velocity \(w>0\) (see Figure~\ref{fig:GenHS2D}).  The geometry of the indentor is described by a given function \(\bar{g}:\left]-1,1\right[\rightarrow\mathbb{R}\), so that the indentor fills the area defined by: \(y\geq\bar{g}(x)\) (\((x,y)\in\left]-1,1\right[\times\mathbb{R}\)).  Set:
\[
g := \frac{\bar{g}}{2(1-\nu^2)}. 
\]

The steady sliding frictional contact problem was introduced and studied in~\cite{bibiElas}.  In a frame moving with the indentor, it is formally that of finding \(\mathbf{t}(x),\mathbf{u}(x): \left]-1,1\right[\rightarrow\mathbb{R}^2\) such that: 
\begin{align*}
	\bullet \quad & - \frac{1}{\pi}\,\text{pv}\frac{1}{x} * t_\textrm{n} - \gamma\,t_{\rm t} = u_\textrm{n}',\qquad \text{in }\left]-1,1\right[,\\[1.2ex]
	\bullet \quad & - \frac{1}{\pi}\,\text{pv}\frac{1}{x} * t_\textrm{t} + \gamma\,t_\textrm{n} = u_{\rm t}',\qquad \text{in }\left]-1,1\right[,\\[1.2ex]
	\bullet \quad & u_\textrm{n} \leq g,\qquad t_\textrm{n}\leq 0,\qquad \bigl(u_\textrm{n} - g\bigr)\,t_\textrm{n}= 0 ,\qquad \text{in }\left]-1,1\right[,\\[1.2ex]
	\bullet \quad & t_{\rm t} = -\bar{f}t_\textrm{n},\qquad \text{in }\left]-1,1\right[,\\[1.2ex]
	\bullet \quad & \int_{-1}^{1} t_\textrm{n}(x)\,{\rm d}x=-P, 
\end{align*}
where \(P> 0\) is the given amplitude of the total normal force exerted from the moving indentor, \(\bar{f}\geq 0\) is a given friction coefficient.  The above convolution products are meant in terms of the extension by zero of \(t_\textrm{n}\), \(t_\textrm{t}\) to the whole real line.  Note that if \(g\) is changed into \(g+C\), where \(C\) denotes an arbitrary constant, then we get a solution for the new problem by just changing \(u_\textrm{n}\) into \(u_\textrm{n}+C\) in the solution.  This means that the penetration of the indentor into the half-space is undefined and this is due to the fact that the displacement field is infinite at infinity.  The problem can be parametrized by the total normal force \(P\) only, and not by the height of the moving obstacle, because it is undetermined.  This fact is intimately connected with the fact that the stress field in the half-space is not square integrable: the elastic energy of the solution is infinite and this is the reason why the problem has to be brought to the boundary by use of the fundamental solution of the Neumann problem for the half-space and singular integrals. 
Focusing on the normal components, this formal problem reduces to that of finding \(t_\textrm{n},u_\textrm{n}:\left]-1,1\right[\rightarrow\mathbb{R}\) such that: 
\[
\begin{aligned}
	\bullet \quad & - \frac{1}{\pi}\,\text{pv}\frac{1}{x} * t_\textrm{n} + \gamma \bar{f}\,t_\textrm{n} = u_\textrm{n}',\qquad \text{in }\left]-1,1\right[, \\[1.2ex]
	\bullet \quad & u_\textrm{n} \leq g,\qquad t_\textrm{n}\leq 0,\qquad \bigl(u_\textrm{n} - g\bigr)\,t_\textrm{n}= 0 ,\qquad \text{in }\left]-1,1\right[,\\[1.2ex]
	\bullet \quad & \int_{-1}^{1} t_\textrm{n}(x)\,{\rm d}x=-P. 
\end{aligned}
\]
Finally, we set \(f=\gamma\bar{f}\) and drop the index `n'.  Given two functions \(f,g:\left]-1,1\right[\rightarrow\mathbb{R}\) and \(P> 0\), the formal problem to be studied is now that of finding \(t,u:\left]-1,1\right[\rightarrow\mathbb{R}\) such that: 
\[
\begin{aligned}
	\bullet \quad & - \frac{1}{\pi}\,\text{pv}\frac{1}{x} * t + f\,t = u',\qquad \text{in }\left]-1,1\right[, \\[1.2ex]
	\bullet \quad & u \leq g,\qquad t\leq 0,\qquad \bigl(u - g\bigr)\,t= 0 ,\qquad \text{in }\left]-1,1\right[, \\[1.2ex]
	\bullet \quad & \int_{-1}^{1} t(x)\,{\rm d}x=-P, 
\end{aligned}
\]
where, again, the convolution product is meant in terms of the extension of \(t\) to the whole real line, by zero.

In the case where \(f\) is a constant, it was proved in \cite{bibiElas} that this formal problem can be put under the form of a monotone variational inequality and that the Lions-Stampacchia theorem yields a unique solution \((t,u)\in H^{-1/2}(-1,1)\times H^{1/2}(-1,1)\) for any given shape \(g\in H^{1/2}(-1,1)\) of the moving indentor and any positive total normal force \(P>0\).

In this paper, we are essentially concerned with the case of a heterogeneous friction coefficient, that is, the case where \(f:\left]-1,1\right[ \rightarrow \mathbb{R}\) is a given function of \(x\).  The situation turns out to be considerably more involved in this case, as the monotonicity property of the underlying linear operator is generally lost.  

As seen in the sequel, insight can be gained by studying first the particular case where the indentor is a rigid flat punch, that is, the case where \(g=0\).  The reason is that it can be proved that any solution of the unilateral problem achieves active contact everywhere.  The problem therefore reduces to a linear problem associated with the so-called Carleman (singular integral) equation.  In the case where the given function \(f\) is Lipschitz-continuous on \(\left]-1,1\right[\), some classical analysis of the Carleman equation is available (see for example \cite[Section 4-4]{tricomi}), based on Marcel Riesz's \(L^p\)-theory of the Hilbert transform in connection with the study of a class of holomorphic functions in the upper complex half-plane.  As the natural case of application of our analysis is the situation where \(f\) is a piecewise constant friction coefficient, we are driven to look for an extension of the classical theory of Carleman equation to the case where \(f\) is piecewise Lipschitz-continuous, that is, to the case where \(f\) is allowed to have jumps.  Interestingly enough, in that extended framework, uniqueness of solution is generally lost for the linear Carleman equation with generic appearance of a finite-dimensional kernel.  However, uniqueness is recovered for the unilateral problem, thanks to the inequality requirement.  In addition, the link made by Marcel Riesz's \(L^p\)-theory of the Hilbert transform between the Carleman equation and some class of holomorphic functions in the upper complex half-plane makes it easy to pass to the homogenized limit of highly oscillating friction coefficients, when the indentor is a rigid flat punch.

In the case of an indentor of arbitrary shape \(g\in H^{1/2}(-1,1)\) and heterogeneous friction coefficient \(f\in BV(-1,1)\), the contact problem is reduced to a variational inequality.  The underlying linear operator is monotone in the case where \(f\) is a non-decreasing function, but not in general.  In the general case, the variational inequality is proved to be pseudomonotone in the sense of Brézis \cite{BrezisIeq}, yielding the existence of a solution to the contact problem.  To obtain the uniqueness of solution and the homogenization analysis, the analysis is further restricted to the particular case of an indentor of arbitrary \emph{convex} shape with Lipschitz-continuous friction coefficient.  In this case, the analysis can be refined, yieding uniqueness of solution and homogenized limit, so that the analysis for the rigid flat punch is fully generalized.  We believe that our analysis may be extended, possibly with major technical adaptations, to the case of nonconvex indentors with piecewise Lipschitz-continuous friction coefficients.

\subsection{Statement of results in the case of the rigid flat punch}

\label{sec:statementFlat}

\begin{prop}
	\label{thm:ResPb}
	Let \(P>0\) and \(f:\left]-1,1\right[\rightarrow\mathbb{R}\) be a given piecewise Lipschitz-continuous function.  There exists a unique \((t,u)\in \cup_{p>1}L^p(-1,1)\times \cup_{p>1}W^{1,p}(-1,1)\) such that:
	\[
	\begin{aligned}
		\bullet \quad & - \frac{1}{\pi}\,\text{\rm pv}\frac{1}{x} * t + f\,t = u',\qquad \text{in }\mathscr{D}'(\left]-1,1\right[), \\[1.2ex]
		\bullet \quad & u \leq 0,\qquad t\leq 0,\qquad u\,t= 0 ,\qquad \text{a.e. in }\left]-1,1\right[, \\[1.2ex]
		\bullet \quad & \int_{-1}^{1} t(x)\,{\rm d}x=-P. 
	\end{aligned}
	\]
	The displacement \(u\) is actually always \(0\) on \(\left]-1,1\right[\) (active contact everywhere).  In the particular case where \(f\) is the piecewise constant function:
	\[
	f(x) = \sum_{i=1}^n f_i \, \chi_{\left]x_{i-1},x_i\right[}(x),
	\]
	where the \(f_i\) and \(-1=x_0 < x_1 < \cdots < x_n =1\) are given real constants and \(\chi_{\left]a,b\right[}\) denotes the characteristic function, then the unique solution \(t\) is explicitly given by:
	\begin{equation}
		\label{eq:explicitSol}
		t(x) = -\frac{P \prod_{i=1}^{n-1}|x-x_{i}|^{\alpha_{i}-\alpha_{i+1}}}{\pi \sqrt{1+f^2(x)}(1+x)^{\frac{1}{2} +\alpha_{1}}(1-x)^{\frac{1}{2}-\alpha_{n}}},   \qquad\text{where } \alpha_i:={\textstyle\frac{1}{\pi}} \arctan f_i.
	\end{equation}
\end{prop}

\begin{figure}[ht]
	\begin{center}
		\begin{tikzpicture}[hatch distance=6pt,hatch thickness=1pt,line cap=round,scale=0.4]
			\path[pattern=flexible hatch,pattern color=blue] (-12,5) rectangle (-11,0.25);
			\path[pattern=flexible hatch,pattern color=blue] (-12,1.25) rectangle (-10,0.25);
			\path[pattern=flexible hatch,pattern color=blue] (-8,1.25) rectangle (-6,0.25);
			\path[pattern=flexible hatch,pattern color=blue] (-4,1.25) rectangle (-2,0.25);
			\path[pattern=flexible hatch,pattern color=blue] (0,1.25) rectangle (2,0.25);
			\path[pattern=flexible hatch,pattern color=blue] (4,1.25) rectangle (6,0.25);
			\path[pattern=flexible hatch,pattern color=blue] (8,1.25) rectangle (10,0.25);
			\path[pattern=flexible hatch,pattern color=blue,
			hatch distance=9pt,hatch thickness=4pt] (12,4) rectangle (11,0.25);
			\path[pattern=flexible hatch,pattern color=blue,
			hatch distance=9pt,hatch thickness=4pt] (12,1.25) rectangle (10,0.25);
			\path[pattern=flexible hatch,pattern color=blue,
			hatch distance=9pt,hatch thickness=4pt] (8,1.25) rectangle (6,0.25);
			\path[pattern=flexible hatch,pattern color=blue,
			hatch distance=9pt,hatch thickness=4pt] (4,1.25) rectangle (2,0.25);
			\path[pattern=flexible hatch,pattern color=blue,
			hatch distance=9pt,hatch thickness=4pt] (0,1.25) rectangle (-2,0.25);
			\path[pattern=flexible hatch,pattern color=blue,
			hatch distance=9pt,hatch thickness=4pt] (-4,1.25) rectangle (-6,0.25);
			\path[pattern=flexible hatch,pattern color=blue,
			hatch distance=9pt,hatch thickness=4pt] (-8,1.25) rectangle (-10,0.25);
			\fill[gray,domain=-14.5:-12.0,samples=30] plot (\x-.2,{exp(0.5*ln(16*(-\x-11.9999)/3))}) -- (-14.7,0) -- (-14.7,3.6515) -- cycle;
			\fill[gray,domain=12.0:14.5,samples=30] plot (\x+.2,{exp(0.5*ln(16*(\x-11.9999)/3))}) -- (14.7,3.6515) -- (14.7,0) -- cycle;
			\fill[gray] (-14.7,0.1) rectangle (14.7,-2.5);
			\draw[line width=2pt,blue] (-12,9)--(-12,0.25)--(12,0.25)--(12,9);
			\draw[line width=2pt,black] (-12.2,0) -- (12.2,0);
			\draw[line width=2pt,black,domain=-15:-12,samples=30] plot
			(\x-.2,{exp(0.5*ln(16*(-\x-11.9999)/3))});
			\draw[line width=2pt,black,domain=12:15,samples=30] plot
			(\x+.2,{exp(0.5*ln(16*(\x-11.9999)/3))});
			\draw [line width=1pt,red] (-11.5,9)
			.. controls +(0.3,-1.5) and +(-0.3,-3) .. (-10.5,9);
			\draw [line width=1pt,red] (-9.5,9)
			.. controls +(0.5,-4) and +(0,6) .. (-8,0.3);
			\draw [line width=1pt,red] (-8,0.3)
			.. controls +(0,5) and +(-0.3,-5) .. (-6.5,9);
			\draw [line width=1pt,red] (-5.5,9)
			.. controls +(0.5,-5) and +(0,4) .. (-4,0.3);
			\draw [line width=1pt,red] (-4,0.3)
			.. controls +(0,3) and +(-0.3,-5) .. (-2.3,9);
			\draw [line width=1pt,red] (-1.8,9)
			.. controls +(0.5,-6) and +(0,3) .. (0,0.3);
			\draw [line width=1pt,red] (0,0.3)
			.. controls +(0,3) and +(-0.3,-6) .. (1.8,9);
			\draw [line width=1pt,red] (2.1,9)
			.. controls +(0.5,-6) and +(0,3) .. (4,0.3);
			\draw [line width=1pt,red] (4,0.3)
			.. controls +(0,2) and +(-0.3,-5) .. (5.8,9);
			\draw [line width=1pt,red] (6.3,9)
			.. controls +(0.5,-4) and +(0,3) .. (8,0.3);
			\draw [line width=1pt,red] (8,0.3)
			.. controls +(0,3) and +(-0.3,-5.5) .. (9.6,9);
			\draw [line width=1pt,red] (10.5,9)
			.. controls +(0.3,-5) and +(-0.3,-3) .. (11.5,9);
			\draw[line width=2pt,blue,->] (13,6) -- (16,6);
			\draw (14.5,6) node[above,blue] {$\mathbf{w}$};
			\draw[line width=2pt,red,->] (0,10) -- (0,8);
			\draw (0,10) node[above,red] {$P$};
		\end{tikzpicture}
		\caption{Normal contact force under a moving flat punch with piecewise constant friction coefficient.}
		\label{fig:flatPunch}
	\end{center}
\end{figure}
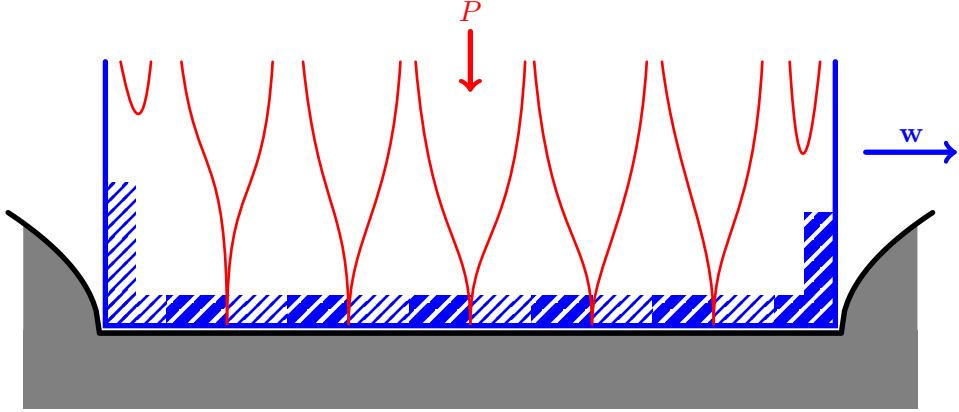

\bigskip
In the case where \(f\) is a constant, formula~\eqref{eq:explicitSol} reduces to the classical solution:
\[
t(x) = -\frac{P}{\pi\sqrt{1+f^2}(1+x)^{1/2+\alpha}(1-x)^{1/2-\alpha}}.
\]
In the case where \(f\) is piecewise constant, the explicit solution seems to be new.  The discontinuities of \(f\) introduce contributions of type \(|x-x_i|^{\alpha_{i}-\alpha_{i+1}}\) in \(t(x)\), so that the normal contact force \(t(x)\) goes to infinity at every increasing discontinuity (\(f(x_i+)>f(x_i-)\)) and goes to zero at every decreasing discontinuity (\(f(x_i+)<f(x_i-)\)) (see Figure~\ref{fig:flatPunch}).

In the case of a highly oscillating periodic piecewise constant friction coefficient, the contact force has therefore very large oscillations, as it goes to zero and to infinity on each period.  It is therefore very natural to investigate the question whether an effective constant friction coefficient can be calculated.  From the point of view of homogenization theory, this amounts to consider a \(2/n\)-periodic friction coefficient \(f_n(x)\) designed by repetition of a (suitably rescaled with respect to \(x\)) given period \(f_1\) and investigate whether the corresponding normal contact force \(t_n\) admits a limit in a suitable sense, as \(n\) goes to infinity.   A full answer is provided for the flat indentor by the following theorem.

\begin{theo}
	\label{thm:Homog}
	Let \(p:\left]-1,1\right[\rightarrow\mathbb{R}\) be a given piecewise Lipschitz-continuous function and extend it to \(\mathbb{R}\) by \(2\)-periodicity.  Let \(f_n:\left]-1,1\right[\rightarrow\mathbb{R}\) be the \(2/n\)-periodic piecewise Lipschitz-continuous function defined on \(\left]-1,1\right[\) by:
	\[
	f_n(x) := p(nx),
	\]
	and let \(t_n\in \cup_{p>1}L^p(-1,1)\) be the unique nonpositive solution of total mass \(-P\) of the equation:
	\[
	- \frac{1}{\pi}\,\text{\rm pv}\frac{1}{x} * t_n + f_n\,t_n = 0,\qquad \text{in }\mathscr{D}'(\left]-1,1\right[),
	\]
	provided by Proposition~\ref{thm:ResPb}.  Let also \(f_\textrm{\rm eff}\) be the constant:
	\[
	f_\textrm{\rm eff} := \tan \bigl\langle \arctan f_1 \bigr\rangle,
	\]
	where \(\langle \cdot\rangle\) denotes the average, that is, \(\langle h\rangle := (1/2)\int_{-1}^1 h\,{\rm d}x\), and let \(t_\textrm{\rm eff}\) be the unique solution of the problem associated with the constant \(f_\textrm{\rm eff}\):
	\[
	t_\textrm{\rm eff}(x) := -\frac{P}{\pi\sqrt{1+f_\textrm{\rm eff}^2}(1+x)^{1/2+\alpha_\textrm{\rm eff}}(1-x)^{1/2-\alpha_\textrm{\rm eff}}},\qquad \text{where } \alpha_\textrm{\rm eff} :={\textstyle\frac{1}{\pi}} \arctan f_\textrm{\rm eff}.
	\]
	Then, the sequences \((t_n)\) and \((f_nt_n)\) converge weakly in \(L^p(-1,1)\), respectively towards \(t_\textrm{\rm eff}\) and \(f_\textrm{\rm eff}t_\textrm{\rm eff}\), for all \(p\in \left]1,(1/2+\beta)^{-1}\right[\), where: 
	\[
	\beta := \frac{1}{\pi}\arctan \bigl\|f_1\bigr\|_{L^\infty(-1,1)} < \frac{1}{2}.
	\]
	In particular, the total tangential contact force \(-\int_{-1}^{1}f_nt_n\,{\rm d}x\) converges towards \(f_\textrm{\rm eff} P\).
\end{theo}

As seen in Section~\ref{sec:proofs}, the proofs of Proposition~\ref{thm:ResPb} and Theorem~\ref{thm:Homog} relies heavily on M. Riesz's \(L^p\)-theory of the Hilbert transform, and the link it induces between the Hilbert transform and a class of holomorphic functions in the upper complex half-plane.  A sketch of the main results of that theory are gathered in Appendix~A of this paper for easy reference.

\subsection{Statement of results in the case of homogeneous friction and indentor of arbitrary shape}

\label{sec:constantFriction}

In this section, the friction coefficient \(f\) is supposed to be an arbitrary real constant \(f\in \mathbb{R}\) (homogeneous friction).  In the case of an indentor with arbitrary shape \(g\in H^{1/2}(-1,1)\), the existence and uniqueness of solution \((t,u)\in H^{-1/2}(-1,1)\times H^{1/2}(-1,1)\) to the contact problem was proved in \cite{bibiJMPS}.  We outline here the structure of the proof for the sake of completeness, before stating new results on the \(L^p\) regularity of the unique solution, based on a Lewy-Stampacchia inequality.  In particular, it will be proved that if the shape of the indentor is Lipschitz-continuous \(g\in W^{1,\infty}(-1,1)\), the unique solution \((t,u)\in H^{-1/2}(-1,1)\times H^{1/2}(-1,1)\) to the contact problem has actually the regularity \((t,u)\in \cup_{p>1}L^p(-1,1)\times \cup_{p>1}W^{1,p}(-1,1)\).

As was seen in the study of the flat obstacle, the integrable function \(t_0:\left]-1,1\right[\rightarrow\mathbb{R}^+\) defined by:
\[
t_0(x) := \frac{1}{\pi\sqrt{1+f^2}{(1+x)}^{\frac{1}{2} +\alpha}{(1-x)}^{\frac{1}{2}-\alpha}} , \qquad \text{ with } \quad \alpha := \frac{1}{\pi} \arctan f \qquad\Bigl(\alpha\in\left] -1/2,1/2 \right[\,\Bigr), 
\]
satisfies:
\begin{equation}
	\label{eq:homogCarlTzero}
	- \frac{1}{\pi}\,\text{pv}\frac{1}{x} * t_0 + f\,t_0 = 0,\qquad \text{in }\left]-1,1\right[,\qquad \text{ and }\qquad \int_{-1}^1 t_0(x) \,{\rm d}x = 1,
\end{equation}
(a proof of the above facts is contained in Proposition~\ref{thm:defTzero} below).  It suggests the following shift on unknown in the contact problem:
\[
\tilde{t}(x) := t(x) + P \,t_0(x),
\]
so that \(\tilde{t}\) must now have zero integral over \(\left]-1,1\right[\).  Accordingly, we introduce the following closed subspace of \(H^{-1/2}(-1,1)\):
\[
H_{0}^* := \Bigl\{\hat{t} \in H^{-1/2}(-1,1)\;\bigm|\; \bigl\langle \hat{t}, 1\bigr\rangle = 0\Bigr\},
\]
where \(1\) denotes the constant function taking value \(1\) all over \(\left]-1,1\right[\).  The dual space of \(H_{0}^*\) is readily seen to be the quotient space \(H_{0}:=H^{1/2}(-1,1)/\mathbb{R}\).  An arbitrary \(\hat{t}\in H_0^*\) defines a distribution in \(H^{-1/2}(\mathbb{R})\) with support in \([-1,1]\) (a sketch of the basic definitions and facts about the spaces \(H^{1/2}\) and \(H^{-1/2}\) can be found in Appendix~B).  Its Fourier transform is a \(C^\infty\) function which vanishes at \(0\).  Using the classical expressions of the Fourier transform of \(\log|x|\) and \(\text{sgn}(x)\) (the sign function) recalled in Proposition~\ref{thm:FourierTransforms} in Appendix~C, this fact entails that the convolution products:
\[
t*\log|x| , \qquad t*\text{\rm sgn}(x), 
\]
are distributions in \(H^{1/2}(\mathbb{R})\) whose restrictions to the interval \(\left]-1,1\right[\) are therefore in \(H^{1/2}(-1,1)\).  In addition, the bilinear form defined by:
\[
t_{1},t_{2}\mapsto - \Bigl\langle t_{1}*\log|x|,t_{2}\Bigr\rangle, 
\]
is symmetric and is also positive definite on \(H_{0}^*\).  It therefore defines a scalar product on the space \(H_{0}^*\), and this scalar product induces a norm that is equivalent to that of \(H^{-1/2}(-1,1)\) (see~\cite[theorem 3]{bibiJiri} or lemma~\ref{thm:lemmHminusHalf} in Appendix~C).  The bilinear form:
\[
t_{1},t_{2}\mapsto \Bigl\langle t_{1}*\text{\rm sgn}(x),t_{2}\Bigr\rangle, 
\]
can easily be seen to be continuous and skew-symmetric on \(H_{0}^*\).  All in all, we have proved that the bilinear form:
\[
t_{1},t_{2}\mapsto - \frac{1}{\pi}\,\Bigl\langle t_{1}*\log|x|,t_{2}\Bigr\rangle + \frac{f}{2}\,\Bigl\langle t_{1}*\text{\rm sgn}(x),t_{2}\Bigr\rangle,
\]
is continuous and coercive on \(H_0^*\).  For an arbitrary \(\tilde{t}\in H_0^*\), the formula:
\[
- \frac{1}{\pi}\,\text{pv}\frac{1}{x} * \tilde{t} + f\,\tilde{t} = \mathring{u}',\qquad \text{in }\mathscr{D}'(\left]-1,1\right[),
\]
defines a unique equivalence class \(\mathring{u}=A\tilde{t}\) in \(H_{0}\), and the operator \(A:H_{0}^* \rightarrow H_0\) is continuous.  Reciprocally, for all \(\mathring{u}\in H_{0}\), the linear problem \(A\tilde{t}=\mathring{u}\) admits a unique solution \(\tilde{t}\in H_{0}^*\), thanks to the Lax-Milgram theorem.  Therefore, the operator \(A:H_{0}^* \rightarrow H_0\) is an isomorphism, and the bilinear form:
\[
\mathring{u}_1,\mathring{u}_2\mapsto \bigl\langle A^{-1}\mathring{u}_1,\mathring{u}_2\bigr\rangle,
\]
is continuous and coercive on \(H_{0}\), thanks to the open mapping theorem.  Incidentally, the restriction of \(A^{-1}\) to \(W^{1,\infty}(-1,1)/\mathbb{R}\) takes values in \(\cup_{p>1}L^p(-1,1)\) and is explicitly given by: 
\begin{equation}
	\label{eq:expliciteAminusOne}
	(A^{-1}\mathring{u})(x) = \frac{f\,\mathring{u}'(x)}{1+f^2}+t_0(x)\, \Biggl\{\text{pv}\frac{1}{x} * \biggl[{(1+x)}^{\frac{1}{2}+\alpha}{(1-x)}^{\frac{1}{2}-\alpha}\mathring{u}'(x)/\sqrt{1+f^2}\biggr]\Biggr\},
\end{equation}
where the integral of the above function over \(\left]-1,1\right[\) vanishes, thanks to formula~\eqref{eq:homogCarlTzero} (for a proof, see~\cite[Section 4-4]{tricomi}) or \cite[theorem 13]{bibiJiri}).

An equivalence class \(\mathring{u}\in H_{0}=H^{1/2}(-1,1)/\mathbb{R}\) will be said \emph{nonpositive} (notation \(\mathring{u}\leq 0\)) if one function at least in the equivalence class is nonpositive.  Hence, stating that \(\mathring{u}\in H_{0}=H^{1/2}(-1,1)/\mathbb{R}\) is nonpositive is equivalent to state that at least one function in the equivalence class is bounded by above, which is equivalent to state that all the functions in the equivalence class are bounded by above\footnote{we recall that \(H^{1/2}(-1,1)\) contains unbounded functions such as \(u(x)=\log|\log|x/2||\).}.

Picking an arbitrary \(g\in H^{1/2}(-1,1)\), it defines a unique equivalence class in \(H^{1/2}(-1,1)/\mathbb{R}\) (which is to be thought of as the function \(g\) up to an arbitrary additive constant) which will be denoted by \(\mathring{g}\in H_0\).  For \(\mathring{g}\in H_{0}\) and \(P\in \left[0,+\infty\right[\), the following subsets:
\begin{align*}
	K_{0} & := \Bigl\{ \mathring{u} \in H_{0}\;\bigm|\; \mathring{u}-\mathring{g}\leq 0\Bigr\}, 
	\\
	K_{0}^* & := \Bigl\{ t \in H_{0}^*\;\bigm|\; t- Pt_{0}\leq 0\Bigr\},
\end{align*}
are obviously non-empty, convex and closed in \(H_0\) and \(H_0^*\), respectively.  Finally, for \(P\geq 0\), the functional defined by:
\[
\varphi (\mathring{u}) := \left\{
\begin{array}{ll}
	\displaystyle \min_{\substack{u-g\in\mathring{u}-\mathring{g}\\
	u-g\leq 0}}\Bigl\langle -Pt_0\;,\; u-g \Bigr\rangle,\quad &\text{if }\mathring{u}\in K_{0},\\[5.0ex]
	+\infty , & \text{otherwise,} 
\end{array}
\right. 
\]
is proper, lower semicontinuous and convex on \(H_{0}\) (see~\cite[lemma 4]{bibiJMPS}).  It is the Legendre-Fenchel conjugate of:
\[
\varphi^*(t) := \Bigl\langle \mathring{g},t \Bigr\rangle + I_{K_0^*}(-t) = \left\{
\begin{array}{ll}
	\displaystyle \Bigl\langle \mathring{g},t \Bigr\rangle,\quad &\text{if }-t\in K_{0}^*,\\[2.0ex]
	+\infty , & \text{otherwise,} 
\end{array}
\right. 
\]

With this notation, the contact problem can be put under the form of any one of the three following equivalent problems.

\medskip
\noindent\textbf{Problem (i)}.  Find \(\mathring{u}\in H_{0}\) and \(\tilde{t}\in H_{0}^*\) such that:
\begin{itemize}
	\item \(\displaystyle - \frac{1}{\pi}\,\text{pv}\frac{1}{x} * \tilde{t} + f\,\tilde{t} = \mathring{u}',\qquad \text{in }\left]-1,1\right[\),
	\item \(\displaystyle \mathring{u}\in K_{0},\qquad \tilde{t}\in K_{0}^*, \qquad \min_{\substack{u-g\in\mathring{u}-\mathring{g}\\
	u-g\leq 0}}\Bigl\langle\tilde{t}-Pt_{0}\;,\;u-g\Bigr\rangle=0\).
\end{itemize}

\medskip
\noindent\textbf{Problem (ii)}.  Find \(\mathring{u}\in H_{0}\) and \(\tilde{t}\in H_{0}^*\) such that:
\begin{itemize}
	\item \(\displaystyle\forall\hat{t}\in H_{0}^*, \quad \bigl\langle A\tilde{t},\hat{t}-\tilde{t}\bigr\rangle + \varphi^*(-\hat{t}) - \varphi^*(-\tilde{t}) \geq 0 \).
	\item \(\displaystyle \mathring{u} = A\tilde{t}\),
\end{itemize}

\medskip
\noindent\textbf{Problem (iii)}.  Find \(\mathring{u}\in H_{0}\) and \(\tilde{t}\in H_{0}^*\) such that:
\begin{itemize}
	\item \(\displaystyle \forall \hat{u}\in H_{0},\qquad \bigl\langle A^{-1}\mathring{u},\hat{u}-\mathring{u}\bigr\rangle + \varphi(\hat{u}) - \varphi(\mathring{u}) \geq 0\).
	\item \(\displaystyle \tilde{t} = A^{-1}\mathring{u}\),
\end{itemize}

Problem~(iii) is the primal weak formulation under the form of a so-called variational inequality and problem~(ii) is the variational inequality corresponding to the dual weak formulation.  Unique solvability of problem~(ii) is a direct consequence of the Lions-Stampacchia theorem.  Hence, each of these three equivalent problems has the same unique solution \((\mathring{u},\tilde{t})\in K_{0}\times K_{0}^*\) (see~\cite{bibiJMPS} for a detailed proof).  Note that neither the operator \(A\), nor \(A^{-1}\) is symmetric (except in the frictionless case \(f=0\)), so that these variational inequalities are not associated with the minimization problem of an underlying energy.  

Starting from an arbitrary \(g\in H^{1/2}(-1,1)\), the above analysis consists in introducing first the corresponding equivalence class \(\mathring{g}\in H_0=H^{1/2}(-1,1)/\mathbb{R}\) and building then a unique solution \(\mathring{u}\in H_0=H^{1/2}(-1,1)/\mathbb{R}\).  Since one prefers to work with functions instead of equivalence classes, provided that \(P>0\), it is always possible to define the function of \(H^{1/2}(-1,1)\):
\[
u-g := \argmin_{\substack{u-g\in\mathring{u}-\mathring{g}\\
u-g\leq 0}}\Bigl\langle -Pt_0\;,\; u-g \Bigr\rangle.
\]
It corresponds to the unique function \(u-g\) in the equivalence class \(\mathring{u}-\mathring{g}\), whose supremum is \(0\).  Finally, given \(g\in H^{1/2}(-1,1)\), \(f\in\mathbb{R}\) and \(P>0\), we conclude that the following problem~$\mathscr{P}$ has a unique solution. 

\medskip
\noindent\textbf{Problem~$\pmb{\mathscr{P}}$.} Find \(u\in H^{1/2}(-1,1)\) and \(\tilde{t}\in H^{-1/2}(-1,1)\) such that:
\[
\begin{aligned}
	\bullet \quad & - \frac{1}{\pi}\,\text{pv}\frac{1}{x} * \tilde{t} + f\,\tilde{t} = u',\qquad \text{in }\left]-1,1\right[, \\[1.2ex]
	\bullet \quad & \tilde{t}-Pt_0\leq 0,\qquad u -g\leq 0, \qquad \bigl\langle \tilde{t}-Pt_0\;,\;u - g\bigr\rangle = 0 ,\\[1.2ex]
	\bullet \quad & \bigl\langle \tilde{t}\;,\;1\bigr\rangle=0. 
\end{aligned}
\]

\medskip
\noindent\textbf{Remark.} The solution of problem~$\pmb{\mathscr{P}}$ obviously satisfies \(\esup(u-g)=0\).

\medskip
Exact explicit solutions for the above steady sliding frictional contact problem are known in the particular cases:
\begin{itemize}
	\item \(\displaystyle g=0\) (moving rigid flat punch), the solution being given by \(u(x)=0\), \(\tilde{t}=0\), that is, \(t(x)=-P\,t_0(x)\),
	\item \(\displaystyle g=x^2/r\) (moving rigid parabola with radius of curvature \(r/(4(1-\nu^2))\) at apex), an extensive derivation of the explicit solution is to be found in~\cite{bibiElas}.  It shows in particular that the coincidence set or contact zone (the set of those \(x\in\left]-1,1\right[\) such that \(u(x)=g(x)\)) depends on the friction coefficient \(f\), in general.
\end{itemize}
These exact explicit solutions seem to have been first discovered by Galin in USSR just after World War II. 

Let us also mention that a catalogue of the universal singularities that can be displayed by the solution of the above steady sliding frictional contact problem is to be found in~\cite{bibiJMPS}.

In this paper, we first add the proof of some new facts about the unique solution of problem~$\pmb{\mathscr{P}}$.  All these facts are based on the following property for the operator \(A^{-1}\).  This property is sometimes called \emph{T-monotonicity} by some authors~\cite[p.~231]{Troianiello}.  It enables to adapt some techniques developed by Stampacchia for the obstacle problem to problem~$\pmb{\mathscr{P}}$.

\begin{theo}
	\label{thm:corodelamort}
	Let \(u\in H^{1/2}(-1,1)\).  Then the function \(u^+(x):=\max\{u(x),0\}\) is in \(H^{1/2}(-1,1)\) and we have:
	\[
	\Bigl\langle A^{-1}\mathring{u}\;,\; u^+ \Bigr\rangle \geq 0,
	\]
	where the equality is achieved if and only if \(u^+\) is constant.
\end{theo}

The detailed proof of Theorem~\ref{thm:corodelamort} is postponed to Section~\ref{sec:homogeneousFriction}

\begin{defi}
	A pair \((u,\tilde{t})\in H^{1/2}(-1,1)\times H^{-1/2}(-1,1)\) which satisfies all the requirements of problem~$\mathscr{P}$, except possibly for the requirement \(\langle \tilde{t}-Pt_0, u-g\rangle=0\) is called a \emph{subsolution} of problem~$\mathscr{P}$.  
\end{defi}

\begin{coro}
	\label{thm:subsolOrd}
	Let \((u,\tilde{t})\) be the solution of problem~$\mathscr{P}$ and \((v,\tilde{r})\) be an arbitrary subsolution.  Then:
	\[
	v\leq u.
	\]
\end{coro}

\noindent\textbf{Proof.}
Let \(z(x):=\max\{u(x),v(x)\}\), so that \(z\in H^{1/2}(-1,1)\) (because \(w\in H^{1/2}\Rightarrow w^+\in H^{1/2}\)).  As \(z-u ={(v-u)}^+\), \((u,\tilde{t})\) is the solution of problem~$\mathscr{P}$ and \(z\leq g\), we have:
\[
\bigl\langle \tilde{t}-Pt_0, {(v-u)}^+\bigr\rangle =\bigl\langle \tilde{t}-Pt_0, z-u \bigr\rangle \geq 0.
\]
Given an arbitrary subsolution \((v,\tilde{r})\), \(\tilde{r}-Pt_0\leq 0\) yields \(\langle \tilde{r}-Pt_0,{(v-u)}^+\rangle \leq 0\).  Gathering, we have:
\[
\bigl\langle A^{-1}(\mathring{v}-\mathring{u})\;,\;{(v-u)}^+\bigr\rangle = \big\langle \tilde{r}-\tilde{t}\;,\; {(v-u)}^+\bigr\rangle \leq 0,
\]
which implies that \({(v-u)}^+\) is a constant, thanks to Theorem~\ref{thm:corodelamort}.  If that constant were positive, then \(v-u\) would be a positive constant, which must be ruled out as we have both \(v\leq g\) and \(\esup (u-g)=0\).  Therefore, \({(v-u)}^+=0\) which is nothing but the claim. \qed

\begin{defi}
	Let \((u,\tilde{t})\) be the unique solution of problem~$\mathscr{P}$.  The \emph{contact zone} \(\mathscr{C}_P\) or \emph{coincidence set} is the closure in \(\left]-1,1\right[\) of the complement in \(\left]-1,1\right[\) of the support of \(u-g\).
\end{defi}

\begin{coro}
	\label{thm:contactZoneIncreasing}
	The contact zone \(\mathscr{C}_P\) is a nondecreasing function of \(P\):
	\[
	0< P_1 < P_2 \qquad \Rightarrow\qquad \mathscr{C}_{P_1}\subset \mathscr{C}_{P_2}. 
	\]
\end{coro}

\noindent\textbf{Proof.}
Let \((u_i,\tilde{t}_i)\) be the solution of problem~$\mathscr{P}$ with total force \(P_i\) (\(i=1,2\)).  Then, \((u_1,\tilde{t}_1)\) is a subsolution of problem~$\mathscr{P}$ with total force \(P_2\).  Corollary~\ref{thm:subsolOrd} yields \(u_1\leq u_2\leq g\), which entails the claim.\qed

\medskip
We now focus on the particular case where the obstacle is Lipschitz-continuous: \(g\in W^{1,\infty}(-1,1)\subset H^{1/2}(-1,1)\).  As the Hilbert transform is an isomorphism of \(L^p(\mathbb{R})\), for all \(p\in\left]1,\infty\right[\), and \(t_0\in L^p(-1,1)\), for all \(1\leq p<{(1/2+|\alpha|)}^{-1}\), formula~\eqref{eq:expliciteAminusOne} shows that, in this case, \(A^{-1}\mathring{g}\in L^p(-1,1)\), for all \(1\leq p<{(1/2+|\alpha|)}^{-1}\), where \(\alpha=(\arctan f)/\pi\).

\begin{theo}
	\label{thm:LewyStampacchia}
	Assuming that \(g\in W^{1,\infty}(-1,1)\subset H^{1/2}(-1,1)\), the solution \((u,\tilde{t})\) of problem~$\mathscr{P}$ satisfies the \emph{Lewy-Stampacchia} inequality:
	\[
	\min\Bigl\{Pt_0,A^{-1}\mathring{g}\Bigr\}\leq \tilde{t} \leq Pt_0.
	\]
\end{theo}

\noindent\textbf{Proof.}  We only have to prove the first inequality.  As \(Pt_0\) and \(A^{-1}\mathring{g}\) are both integrable functions, their pointwise minimum is well defined.  The closed convex subset of \(H_0^*\):
\begin{equation}
	\label{eq:convexSubset}
	\Bigl\{ t \in H_0^* \;\bigm|\; t \leq - \min\{Pt_0,A^{-1}\mathring{g}\} \Bigr\},
\end{equation}
is non-empty, as it contains \(-A^{-1}\mathring{g}\).  Therefore, as \(A:H_0^*\rightarrow H_0\) is continuous and coercive, the variational inequality in problem~(ii) above still has a unique solution when \(\mathring{g}\) is replaced by \(-\mathring{u}\) and \(Pt_0\) by \(-\min\{Pt_0,A^{-1}\mathring{g}\}\) (that is, \(K_0^*\) is replaced by \eqref{eq:convexSubset}), thanks to the Lions-Stampacchia theorem.  Hence, the contact problem (problem (i), (ii) or (iii) above) obtained by replacing \(\mathring{g}\) by \(-\mathring{u}\) and \(Pt_0\) by \(-\min\{Pt_0,A^{-1}\mathring{g}\}\) has a unique solution which is denoted by \((-\mathring{v},-A^{-1}\mathring{v})\).  It satisfies in  particular \(-\mathring{v}\leq-\mathring{u}\) and:
\begin{equation}
	\label{eq:conclforv}
	-A^{-1}\mathring{v} \leq  -\min\{Pt_0,A^{-1}\mathring{g}\}.
\end{equation}
We denote by \(u-v\) the unique function in the equivalence class \(\mathring{u} - \mathring{v}\) whose supremum is \(0\).  By construction of the nonpositive function \(u-v\), we have:
\begin{equation}
	\label{eq:pourLaSuite}
	\Bigl\langle A^{-1}\mathring{v} - \min\{Pt_0,A^{-1}\mathring{g}\}\;,\;v-u \Bigr\rangle = 0.
\end{equation}
Note that we also have:
\[
-g \leq -u,\qquad \text{and}\qquad -A^{-1}\mathring{g} \leq -\min\{Pt_0,A^{-1}\mathring{g}\} ,
\]
which shows that \((-g,-A^{-1}\mathring{g})\) is a subsolution of the contact problem solved by \(-v\).  Therefore, Corollary~\ref{thm:subsolOrd} yields:
\[
-g \leq -v ,\qquad \Rightarrow\qquad u\leq v \leq g.
\]
But, this entails that we can use \(v\) as a test function in the variational inequality that defines \(u\):
\[
\Bigl\langle A^{-1}\mathring{u} - Pt_0\;,\;v-u \Bigr\rangle = \Bigl\langle A^{-1}\mathring{u} - Pt_0\;,\;v-g\Bigr\rangle - \Bigl\langle A^{-1}\mathring{u} - Pt_0\;,\;u-g\Bigr\rangle \geq 0,
\]
which entails:
\[
\Bigl\langle A^{-1}\mathring{u} - Pt_0+{(Pt_0 - A^{-1}\mathring{g})}^+\;,\;v-u \Bigr\rangle \geq 0,
\]
that is:
\[
\Bigl\langle A^{-1}\mathring{u} - \min\{Pt_0,A^{-1}\mathring{g}\}\;,\;v-u \Bigr\rangle \geq 0.
\]
Taking the difference of this last inequality with identity~\eqref{eq:pourLaSuite}, we obtain:
\[
\Bigl\langle A^{-1}(\mathring{u}-\mathring{v}) \;,\;u-v \Bigr\rangle \leq 0,
\]
which shows that \(u-v\) is a constant function, as \(A^{-1}:H_0\rightarrow H_0^*\) is coercive.  As this constant function has \(0\) supremum, it vanishes identically.  So, we have proved that \(v=u\).  But then, the claimed inequality is given by inequality~\eqref{eq:conclforv}.\qed

\begin{coro}
	Let \(P>0\) and \(g\in W^{1,\infty}(-1,1)\).  Then, the solution \((u,\tilde{t})\) of problem~$\mathscr{P}$ fulfils actually the additional regularity: \(\tilde{t}\in L^p(-1,1)\), \(u\in W^{1,p}(-1,1)\), for all \(p\) such that \(1<p<{(1/2+|\alpha|)}^{-1}\), with \(\alpha=(\arctan f)/\pi\).
\end{coro}

\noindent\textbf{Proof.}
We have \(t_0\in L^p(-1,1)\), for all \(p\) such that \(1<p<{(1/2+|\alpha|)}^{-1}\).  As \(g\in W^{1,\infty}(-1,1)\), formula~\eqref{eq:expliciteAminusOne} shows that \(A^{-1}\mathring{g}\in L^p(-1,1)\) (we recall that the Hilbert transform maps continuously $L^q(\mathbb{R})$ onto itself, for all $q\in\left]1,\infty\right[$, by Theorem~\ref{thm:Riesz1} of Appendix~A).  By Theorem~\ref{thm:LewyStampacchia}, this is also true of \(\tilde{t}\), and by the first equation in problem~$\mathscr{P}$, this is true of $u'$.\qed

\begin{coro}
	\label{thm:convexIndentor}
	Let \(P>0\) and \(g\in W^{1,\infty}(-1,1)\) be \emph{convex}.  Then, the contact zone \(\mathscr{C}_P\) coincides with the support of \(t=\tilde{t}-Pt_0\) in \(\left]-1,1\right[\) and is connected, that is, it is an interval.  Denoting by \(a<b\in [-1,1]\) the bounds of that interval, we have:
	\begin{itemize}
		\item \(\displaystyle u < g,\qquad g' < u' < 0,\qquad \text{on }\left]-1,a\right[,\)
		\item \(\displaystyle u < g,\qquad 0 < u' < g',\qquad \text{on }\left]b,1\right[.\)
	\end{itemize}
	In particular, \(u\in W^{1,\infty}(-1,1)\) and \(\|u'\|_{L^\infty(-1,1)} = \|g'\|_{L^\infty(-1,1)}\).
\end{coro}

\noindent\textbf{Proof.}
If \(\mathscr{C}_P=\left]-1,1\right[\), then there is nothing to prove.  Otherwise, there exists \(x_0\) such that \(u(x_0)<g(x_0)\).  As \(u\) is continuous by the previous corollary, let \(\left]c,d\right[\) be the largest open interval containing \(x_0\) (\(x_0\in\left]c,d\right[\subset\left]-1,1\right[\)) in which \(u<g\).  We are going to prove that we cannot have \(-1<c<d<1\), which is sufficient to prove that the contact zone \(\mathscr{C}_P\) is connected.  So, let us suppose \(-1<c<d<1\).  By the continuity of \(u\) given by the preceding corollary, we have \(u(c)-g(c)=u(d)-g(d)=0\).  Furthermore, \(t=\tilde{t}-Pt_0\) must vanish on \(\left]c,d\right[\), so that, by using Theorem~\ref{thm:Riesz1} in Appendix~A:
\[
\forall x \in \left]c,d\right[,\qquad u'(x) = \frac{1}{\pi}\int_{-1}^{c}\frac{t(s)}{s-x}\,{\rm d}s + \frac{1}{\pi}\int_{d}^{1}\frac{t(s)}{s-x}\,{\rm d}s.
\]
As \(t\leq 0\), \(u'\) is decreasing on \(\left]c,d\right[\) and therefore \(u-g\) is concave on \(\left]c,d\right[\).  But \(u-g\) is also nonpositive and vanishes at \(c\) and \(d\).  Therefore, it must vanish identically on \(\left]c,d\right[\), yielding the expected contradiction.  Besides, note that the contradiction is reached under the sole hypothesis that there exist \(-1<c<d<1\) such that \(c,d\in\text{supp}\,t\) and \(t=0\) on \(\left]c,d\right[\).  Hence, both the contact zone and the support of \(t\) are intervals.  Let \(a<b \in [-1,1]\) be the bounds of the support of \(t\).  We have \(u=g\) on \(\left]a,b\right[\).  Also, \(u'\) is negative and decreasing on \(\left]-1,a\right[\), which entails \(u<g\) on \(\left]-1,a\right[\).  As \(u'(c-)\geq g'(c-)\), we must have \(g'<u'<0\) on \(\left]-1,a\right[\).  By the same argument, \(u<g\) and \(0<u'<g'\) on \(\left]b,1\right[\).  As \(u'=g'\) on \(\left]a,b\right[\), we get \(\|u'\|_{L^\infty(-1,1)} = \|g'\|_{L^\infty(-1,1)}\). \qed

\subsection{Existence of solutions in the case of heterogeneous friction and indentor of arbitrary shape}

\label{sec:statementHetero}

We now turn to the problem of generalizing the existence and uniqueness results for the case of homogeneous friction, to the case of heterogeneous friction.  The following two difficulties arise.
\begin{itemize}
	\item  As, in applications, heterogeneous friction coefficients are going to be piecewise constant, our framework should encompass discontinuous friction coefficients.  The restriction of a distribution \(t\in H^{-1/2}(-1,1)\) to \(\left]0,1\right[\) does not define a distribution on \(\left]-1,1\right[\), in general.  Hence, there are some \(t\in H^{-1/2}(-1,1)\) for which we cannot define the product \(ft\) and the functional framework must therefore be redesigned in this case.   Note, however, that as \(t\) must be nonpositive, it belongs to the Banach space of Radon measures \(\mathscr{M}([-1,1])\).  Any \(t\in H^{-1/2}\cap \mathscr{M}([-1,1])\) is a Radon measure without any atom, and a function \(f\in \mathit{BV}([-1,1])\) has a countable set of discontinuities.  Therefore, the product \(ft\) is a well-defined Radon measure in that case.  Hence, we are driven to look for \(t\) in \(H^{-1/2}\cap \mathscr{M}([-1,1])\), and one must therefore leave the framework of Hilbert spaces. 
	\item  In the case where  \(f\in \mathit{BV}([-1,1])\) and \(\tilde{t}\in H_0^*\cap \mathscr{M}([-1,1])\), it will be seen in the sequel that the identity:
	\[
	- \frac{1}{\pi}\,\text{pv}\frac{1}{x} * \tilde{t} + f\,\tilde{t} = \mathring{u}',\qquad \text{in }\left]-1,1\right[,
	\]
	defines uniquely \(\mathring{u}=A\tilde{t}\) in \( H_{0}+C^0([-1,1])\cap BV([-1,1])\subset H_{0}+C^0([-1,1])\).  However, it will be seen that the operator \(A\) is not strictly monotone (or even monotone) in general.  This could be anticipated from the analysis of the flat punch case, as Proposition~\ref{thm:defTi} shows that, in the case of a piecewise constant \(f\), the kernel of \(A\) can contain a finite-dimensional space.
\end{itemize}

However, we will prove that \(A:H_0^*\cap \mathscr{M}([-1,1])\rightarrow  H_{0}+C^0([-1,1])\) is \emph{pseudomonotone in the sense of Brézis}.  As a consequence, the contact problem always admits a solution under the only hypotheses that \(f\in \mathit{BV}([-1,1])\) and \(g\in H^{1/2}(-1,1)\).  In addition, in the particular case where the function \(f:[-1,1]\rightarrow\mathbb{R}\) is \emph{non-decreasing}, the operator \(A\) is strictly monotone and the solution is therefore unique.

To formulate the contact problem under the form of a variational inequality, we first extend the definition of \(t_0\) to less regular friction coefficient.  Interestingly enough, the proof of the following proposition (postponed to Section~\ref{sec:proofGeneralHetero}) relies on the same technique that yielded the homogenization of friction in the case of the flat indentor.

\begin{prop}
	\label{thm:defTzeroLinfty}
	Let \(f\in L^\infty(-1,1)\), extended by zero on \(\mathbb{R}\setminus\left]-1,1\right[\).  We denote by \(\tau:=-\frac{1}{\pi}\,\text{\rm pv}\frac{1}{x}*\arctan f\), the Hilbert transform of the function \(\arctan f\).  Let \(t_0\) be the function defined by:
	\[
	t_0(x):= \frac{e^{\tau(x)}}{\pi\sqrt{1-x^2}\sqrt{1+f^2(x)}},\qquad \text{for x}\in \left]-1,1\right[,
	\]
	and extended by zero on \(\mathbb{R}\setminus\left]-1,1\right[\).  Then, \(t_0\in \cup_{p>1}L^p(-1,1)\), it is obviously positive on \(\left]-1,1\right[\), has total mass \(\int_{-1}^1t_0\,{\rm d}x=1\) and solves the homogeneous Carleman equation:
	\[
	- \frac{1}{\pi}\,\text{\rm pv}\frac{1}{x} * t_0 + f\,t_0 = 0,\qquad \text{a.e. in }\left]-1,1\right[.
	\]
\end{prop}

Note that the above proposition provides incidentally a solution of the contact problem for the moving flat punch with arbitrary friction coefficient \(f\in L^\infty(-1,1)\).

\medskip
Let \(B_1\), \(B_2\) be Banach spaces that are both continuously embedded in some Hausdorff topological vector space.  Then, \(B_1+B_2\) and \(B_1\cap B_2\) are both Banach spaces for the natural norms:
\begin{align*}
	\|x\|_{B_1+B_2} & := \inf_{\substack{b_1\in B_1,b_2\in B_2,\\
	b_1+b_2=x}} \bigl\{\|b_1\|_{B_1}+\|b_2\|_{B_2}\bigr\},\\
	\|x\|_{B_1\cap B_2} & := \max\bigl\{ \|x\|_{B_1} , \|x\|_{B_2} \bigr\}.
\end{align*}
In addition, we have \({(B_1+B_2)}^*=B_1^*\cap B_2^*\).  Hence, the Banach space \(H_{0}^*\cap\mathscr{M}([-1,1])\) is the dual of the Banach space \(H_0+C^0([-1,1])\).  

\medskip
We are given \(P\geq 0\) and \(g\in H^{1/2}(-1,1)+\mathit{BV}([-1,1])\).  The set:
\begin{equation}
	\label{eq:defConv}
	K_{0}^* := \Bigl\{ t \in H_{0}^*\;\bigm|\; t- Pt_{0}\leq 0\Bigr\},
\end{equation}
is a nonempty (\(0\in K_{0}^*\)), convex, closed subset of \(H_{0}^*\cap\mathscr{M}([-1,1])\) (where \(H_0^*:=\{\hat{t}\in H^{-1/2}(-1,1)|\langle\hat{t},1\rangle=0\}\)).  We will associate with \(g\), a unique \(\mathring{g}\in H_0+\mathit{BV}([-1,1])\) as in Section~\ref{sec:constantFriction}, and write \(\mathring{u}\leq \mathring{g}\) to mean, as there, that one (and therefore all) function \(u-g\) in the equivalence class \(\mathring{u}-\mathring{g}\) is bounded by above.  The reason for allowing \(g\in H^{1/2}(-1,1)+\mathit{BV}([-1,1])\) is that this new setting encompasses indentors that may have steps.  This is interesting from the point of view of mechanics.  From the point of view of mathematics, there is no additional difficulty as any function in \(\mathit{BV}([-1,1])\) is the uniform limit of a sequence of step functions and is therefore integrable with respect to any Radon measure.  In addition, the linear functional \(t\mapsto\int g t\) is continuous on \(H_{0}^*\cap\mathscr{M}([-1,1])\), so that the notation \(t\mapsto\langle t,g\rangle\) can be used unambiguously.

Any \(t\in H^{-1/2}\cap \mathscr{M}([-1,1])\) is a Radon measure without any atom, and the same is true of \(ft\).  Hence, we can write unambiguously \(\int_0^x ft\), for any \(x\in [-1,1]\) and the function \(x\mapsto \int_0^x ft\) is a continuous function on \([-1,1]\) with bounded variation.  Also, \(\log|x|*t\) is in \(H^{1/2}(-1,1)\).  Therefore, given \(\tilde{t}\in H_0^*\cap \mathscr{M}([-1,1])\), the identity:
\[
- \frac{1}{\pi}\,\text{pv}\frac{1}{x} * \tilde{t} + f\,\tilde{t} = \mathring{u}',\qquad \text{in }\mathscr{D}'(\left]-1,1\right[),
\]
defines uniquely \(\mathring{u}=A\tilde{t}\) in \( H_{0}+\bigl(C^0([-1,1])\cap BV([-1,1])\bigr)\subset H_{0}+C^0([-1,1])\). 

With this notation, we can now easily generalize problems~(i) and (ii) of Section~\ref{sec:constantFriction} to the case where \(f\in\mathit{BV}([-1,1])\).

\medskip
\noindent\textbf{Problem (i$'$)}.  Find \(\mathring{u}\in H_{0}+C^0([-1,1])\cap BV([-1,1])\) and \(\tilde{t}\in H_{0}^*\cap\mathscr{M}([-1,1])\) such that:
\begin{itemize}
	\item \(\displaystyle - \frac{1}{\pi}\,\text{pv}\frac{1}{x} * \tilde{t} + f\,\tilde{t} = \mathring{u}',\qquad \text{in }\mathscr{D}'(\left]-1,1\right[)\),
	\item \(\displaystyle \mathring{u}\leq \mathring{g},\qquad \tilde{t}\leq Pt_0, \qquad \min_{\substack{u-g\in\mathring{u}-\mathring{g}\\
	u-g\leq 0}}\Bigl\langle\tilde{t}-Pt_{0}\;,\;u-g\Bigr\rangle=0\).
\end{itemize}

\smallskip
\noindent\textbf{Problem (ii$'$)}.  Find \(\mathring{u}\in H_{0}+C^0([-1,1])\cap BV([-1,1])\) and \(\tilde{t}\in K_{0}^*\) such that:
\begin{itemize}
	\item \(\displaystyle \forall\hat{t}\in K_{0}^*, \qquad \bigl\langle A\tilde{t},\hat{t}-\tilde{t}\bigr\rangle\geq \bigl\langle g, \hat{t}-\tilde{t}\bigr\rangle\).
	\item \(\displaystyle \mathring{u}=A\tilde{t}\),
\end{itemize}

\smallskip
\begin{prop}
	\label{thm:eqFav}
	Problems~{\rm (i$'$)} and {\rm (ii$'$)} are equivalent in the sense that any solution of problem~{\rm (i$'$)} is a solution of problem~{\rm (ii$'$)}, and reciprocally.
\end{prop}

\noindent\textbf{Proof.}
First, consider a solution \((\mathring{u},\tilde{t})\) of problem~(i$'$).  By the second condition of problem~(i$'$), there exists \(u\in \mathring{u}\) such that \(\langle \tilde{t}-Pt_{0}\;,\;u-g\rangle=0\) and \(u-g\leq 0\).  By the first condition, there exists \(C\in\mathbb{R}\) such that:
\[
-\frac{1}{\pi}\,\log|x|*\tilde{t}+\frac{1}{2}\,\text{sgn}(x)*(f\tilde{t}) = u+C, \qquad \text{in }\left]-1,1\right[,
\]
so that, for all \(\hat{t}\in K_0^*\):
\[
\Bigl\langle A\tilde{t},\hat{t}-\tilde{t}\Bigr\rangle-\Bigl\langle g, \hat{t}-\tilde{t}\Bigr\rangle = \Bigl\langle \hat{t}-Pt_{0}\;,\;u-g\Bigr\rangle - \Bigl\langle \tilde{t}-Pt_{0}\;,\;u-g\Bigr\rangle \geq 0,
\]
that is, \(\mathring{u}\), \(\tilde{t}\) solves problem~(ii$'$).

Reciprocally, consider a solution \((\mathring{u},\tilde{t})\) of problem~(ii$'$).  Pick an arbitrary \(u\in\mathring{u}\).  We have, for all \(\hat{t}\in K_0^*\):
\begin{equation}
	\label{eq:statContradict}
	\Bigl\langle \hat{t}-Pt_{0}\;,\;u-g\Bigr\rangle - \Bigl\langle \tilde{t}-Pt_{0}\;,\;u-g\Bigr\rangle = \Bigl\langle A\tilde{t},\hat{t}-\tilde{t}\Bigr\rangle-\Bigl\langle g, \hat{t}-\tilde{t}\Bigr\rangle \geq 0.
\end{equation}
Let us show that \(u-g\) is (essentially) bounded by above.  If it were not, then the sets:
\[
K_n:=\Bigl\{x\in [-1,1]\;\bigm|\; u(x)-g(x)\geq n\Bigr\},
\]
would all have a positive measure, and the sequence \((\hat{t}_n)\) in \(K_0^*\) defined by:
\[
\hat{t}_n := Pt_0-\frac{P}{|K_n|}\chi_{K_n},
\]
(where \(\chi_{K_n}\) is the characteristic function of \(K_n\)) would be such that \(\lim_{n \rightarrow +\infty}\langle \hat{t}_n-Pt_{0}\;,\;u-g\rangle = -\infty\) which would contradict~\eqref{eq:statContradict}.  Hence, \(u-g\) is (essentially) bounded by above, which shows that \(\mathring{u}\leq \mathring{g}\).  Therefore, by~\eqref{eq:statContradict}, for all \(\hat{t}\in K_0^*\):
\[
\min_{\substack{u-g\in\mathring{u}-\mathring{g}\\
u-g\leq 0}}\Bigl\langle\tilde{t}-Pt_{0}\;,\;u-g\Bigr\rangle=\Bigl\langle \tilde{t}-Pt_{0}\;,\;u-g-\sup_{\left]-1,1\right[}(u-g)\Bigr\rangle \leq \Bigl\langle \hat{t}-Pt_{0}\;,\;u-g-\sup_{\left]-1,1\right[}(u-g)\Bigr\rangle.
\]
Once again, one can easily construct a sequence \((\hat{t}_n)\) in \(K_0^*\) such that \(\lim_{n \rightarrow +\infty}\langle \hat{t}_n-Pt_{0}\;,\;u-g-\sup(u-g)\rangle = 0\), which shows that the value of the minimum is \(0\) and therefore that \((\mathring{u},\tilde{t})\) solves problem~(i$'$).
\qed

\begin{prop}
	\label{thm:boundedA}
	The linear operator \(A:H_{0}^*\cap\mathscr{M}([-1,1])\rightarrow  H_0+C^0([-1,1])\) is bounded.
\end{prop}

\noindent\textbf{Proof.} We have:
\[
\bigl\| At\bigr\|_{H_0+C^0} \leq \Bigl\| (-1/\pi) (\log|x|*t)\Bigr\|_{H^{1/2}} + \Bigl\| {\textstyle \int_0^x ft}\Bigr\|_{C^0} \leq M \Bigl( \bigl\| t \bigr\|_{H^{-1/2}} + \bigl\| t \bigr\|_{\mathscr{M}} \Bigr)\leq 2M  \bigl\| t \bigr\|_{H_0^*\cap\mathscr{M}},
\]
for some real constant \(M\) depending only on \(\|f\|_{L^\infty}\) and, in particular, independent of \(t\) (here the Fourier transform of the logarithm provided by Proposition~\ref{thm:FourierTransforms} in Appendix~C was used to obtain the second inequality).
\qed

\begin{theo}
	\label{thm:pseudomonot}
	In the case of an arbitrary function \(f\in \mathit{BV}([-1,1])\), the  bounded linear operator \(A:H_{0}^*\cap\mathscr{M}([-1,1])\rightarrow  H_0+C^0([-1,1])\) is pseudomonotone in the sense of Brézis (see Definition~\ref{thm:defPseudoMonot}), with the space \(H_{0}^*\cap\mathscr{M}([-1,1])\) endowed with weak-* topology and \(H_0+C^0([-1,1])\) endowed with the weak topology.  In addition, in the particular case where the function \(f\) is nondecreasing,  the bounded linear operator \(A:H_{0}^*\cap\mathscr{M}([-1,1])\rightarrow  H_0+C^0([-1,1])\) is strictly monotone (see Definition~\ref{thm:defMonot}).
\end{theo}

The proof of theorem~\ref{thm:pseudomonot} is postponed to Section~\ref{sec:proofGeneralHetero}.

\begin{coro}
	\label{thm:existUniq}
	In the case of an arbitrary function \(f\in \mathit{BV}([-1,1])\), problem~{\rm (ii$'$)} always has a solution.  In the particular case where \(f\) is nondecreasing, this solution is unique.
\end{coro}

\noindent\textbf{Proof.}
The set \(K_{0}^*\) is a closed convex subset of \(H_{0}^*\cap\mathscr{M}([-1,1])\) that contains \(0\).  In addition, the set \(\{-Pt_0\}+K_{0}^*\) contains only nonpositive measures of total mass \(-P\) and is therefore a subset of the closed ball of radius \(P\) in \(\mathscr{M}([-1,1])\).  Hence, a sequence in \(K_{0}^*\) whose \(H_{0}^*\cap\mathscr{M}([-1,1])\)-norm goes to infinity has also its \(H^{-1/2}(-1,1)\)-norm going to infinity.  Hence, we have:
\[
\lim_{\substack{\|t\|\rightarrow +\infty\\
t\in K_{0}^*}} \frac{\langle At,t\rangle+\langle g,t\rangle}{\|t\|} = \lim_{\substack{\|t\|\rightarrow +\infty\\
t\in K_{0}^*}} \frac{-\frac{1}{\pi} \Bigl\langle \log|x|*t,t\Bigr\rangle+\langle g,t\rangle}{\|t\|} = +\infty,
\]
where \(\|\cdot\|\) is the norm of \(H_{0}^*\cap\mathscr{M}([-1,1])\) (here, the second term in $\langle At,t\rangle$ has been removed since it is bounded).  The claim is now a direct consequence of Corollary~\ref{thm:coroBrez} in Appendix~D and Theorem~\ref{thm:pseudomonot}. \qed

\subsection{Regularity, uniqueness and homogenization in the case of heterogeneous friction and convex indentor}

\label{sec:statementConvex}

In this section, we are restricted to the case where the shape of the indentor \(g\in W^{1,\infty}(-1,1)\) is Lipschitz-continuous and \emph{convex}.  The heterogeneous friction coefficient \(f:\left]-1,1\right[\rightarrow\mathbb{R}\), will be supposed Lipschitz-continuous.  There is no doubt that the whole analysis could be generalized to the case of a \emph{piecewise} Lipschitz-continuous friction coefficient, as in the case of the moving flat punch, but we will restrict ourselves to the case of a Lipschitz-continous friction coefficient, for the sake of simplicity.  In this restricted case, we are going to be able to generalize all the results already proved in the case of the moving flat punch with piecewise Lipschitz-continuous friction coefficient: existence and uniqueness of a solution \((t,u)\in \cup_{p>1}L^p(-1,1)\times\cup_{p>1}W^{1,p}(-1,1)\) to the contact problem and homogenization analysis with the same value of the effective friction coefficient in the homogenized limit.  This generalization has an important added value as, in the case of the general convex indentor, contrary to the case of the flat punch, the contact zone is a true unknown of the problem and is different \emph{a priori} for each oscillating friction coefficient \(f_n\).  The fact that we obtain the same value for the effective friction coefficient in the homogenized limit shows that there is no influence of unilateral contact on homogenized friction.

In this section, \(f\) denotes a Lipschitz-continuous function on \(\left]-1,1\right[\), extended as usual by zero to the whole real line.  We also consider a function \(\bar{f}\), which is identically equal to \(f\) on \(\left]-1,1\right[\),  and chosen to be compactly supported in \(\mathbb{R}\) and Lipschitz-continuous on \(\mathbb{R}\) (which \(f\) is not, in general).  We denote by \(\tau:=-\frac{1}{\pi}\text{pv}\frac{1}{x}*\arctan f\) and \(\bar{\tau}:=-\frac{1}{\pi}\text{pv}\frac{1}{x}*\arctan \bar{f}\) the Hilbert transforms of \(\arctan f\) and \(\arctan\bar{f}\), respectively, and we keep the notation:
\[
t_0(x):= \frac{e^{\tau(x)}}{\pi\sqrt{1-x^2}\sqrt{1+f^2(x)}},\qquad \text{for }x\in \left]-1,1\right[,
\]
which, we recall, is a solution of the homogeneous Carleman equation in \(\cup_{p>1}L^p(-1,1)\) and has total mass \(1\).

We shall now prove that the original contact problem (problem~\(\mathscr{P}_{\!\rm o}\) below) which was seen in previous section to have a default of monotonicity, is equivalent to an auxiliary contact problem (problem~\(\mathscr{P}_{\!\rm a}\) below) associated with a monotone operator.  These two problems are defined as follows.

\bigskip
\noindent\textbf{Problem~\(\mathscr{P}_{\!\rm o}\).} Find \((\tilde{t},u)\in \cup_{p>1}L^p(-1,1)\times \cup_{p>1} W^{1,p}(-1,1)\) such that:
\[
\begin{aligned}
	\bullet \quad & - \frac{1}{\pi}\,\text{\rm pv}\frac{1}{x} * \tilde{t} + f\,\tilde{t} = u',\qquad \text{in }\mathscr{D}'(\left]-1,1\right[), \\[1.0ex]
	\bullet \quad & u-g \leq 0,\qquad \tilde{t}-Pt_0\leq 0,\qquad \bigl(u-g\bigr)\bigl(\tilde{t}-Pt_0\bigr)= 0 ,\qquad \text{a.e. in }\left]-1,1\right[, \\[1.0ex]
	\bullet \quad & \int_{-1}^{1} \tilde{t}\,{\rm d}x=0. 
\end{aligned}
\]

\bigskip
\noindent\textbf{Problem~\(\mathscr{P}_{\!\rm a}\).} Find \((\tilde{t},v)\in \cup_{p>1}L^p(-1,1)\times \cup_{p>1} W^{1,p}(-1,1)\) such that:
\[
\begin{aligned}
	\bullet \quad & \frac{e^{-\bar{\tau}}}{\scriptstyle\sqrt{1+f^2}}\biggl(- \frac{1}{\pi}\,\text{\rm pv}\frac{1}{x} * \tilde{t} + f\,\tilde{t}\biggr) = v',\qquad \text{in }\mathscr{D}'(\left]-1,1\right[), \\[1.0ex]
	\bullet \quad & v(x)-\int_{-1}^x \frac{g'e^{-\bar{\tau}}}{\scriptstyle\sqrt{1+f^2}}\,{\rm d}s \leq 0,\qquad \tilde{t}(x)-Pt_0(x)\leq 0,\\
	&\hspace*{3cm} \biggl(v(x)-\int_{-1}^x \frac{g'e^{-\bar{\tau}}}{\scriptstyle\sqrt{1+f^2}}\,{\rm d}s\biggr)\Bigl(\tilde{t}(x)-Pt_0(x)\Bigr)= 0 ,\qquad \text{a.e. in }\left]-1,1\right[, \\[1.0ex]
	\bullet \quad & \int_{-1}^{1} \tilde{t}\,{\rm d}x=0. 
\end{aligned}
\]

\begin{prop}
	\label{thm:eqPb}
	We assume that \(P>0\), \(f\in W^{1,\infty}(-1,1)\) and \(g\in W^{1,\infty}(-1,1)\) is convex.  Then, problem~\(\mathscr{P}_{\!\rm o}\) and  problem~\(\mathscr{P}_{\!\rm a}\) are equivalent in the following sense.
	\begin{itemize}
		\item If \((\tilde{t},u)\in \cup_{p>1}L^p(-1,1)\times \cup_{p>1} W^{1,p}(-1,1)\) denotes an arbitrary solution of problem~\(\mathscr{P}_{\!\rm o}\), then, setting:
		\[
		v(x) := \int_{-1}^x \frac{u'e^{-\bar{\tau}}}{\scriptstyle\sqrt{1+f^2}}\,{\rm d}s - \sup_{x\in\left]-1,1\right[}\int_{-1}^x \frac{(u'-g')e^{-\bar{\tau}}}{\scriptstyle\sqrt{1+f^2}}\,{\rm d}s,
		\]
		the pair \((\tilde{t},v)\in \cup_{p>1}L^p(-1,1)\times \cup_{p>1} W^{1,p}(-1,1)\) yields a solution of problem~\(\mathscr{P}_{\!\rm a}\).
		\item If \((\tilde{t},v)\in \cup_{p>1}L^p(-1,1)\times \cup_{p>1} W^{1,p}(-1,1)\) denotes an arbitrary solution of problem~\(\mathscr{P}_{\!\rm a}\), then, setting:
		\[
		u(x) := \int_{-1}^x v'e^{\bar{\tau}}{\scriptstyle\sqrt{1+f^2}}\,{\rm d}s - \sup_{x\in\left]-1,1\right[}\biggl(-g(x)+\int_{-1}^x v'e^{\bar{\tau}}{\scriptstyle\sqrt{1+f^2}}\,{\rm d}s\biggr),
		\]
		the pair \((\tilde{t},u)\in \cup_{p>1}L^p(-1,1)\times \cup_{p>1} W^{1,p}(-1,1)\) yields a solution of problem~\(\mathscr{P}_{\!\rm o}\).
	\end{itemize}
	In addition, if \((\tilde{t},u)\in \cup_{p>1}L^p(-1,1)\times \cup_{p>1} W^{1,p}(-1,1)\) denotes an arbitrary solution of problem~\(\mathscr{P}_{\!\rm o}\), then there exist \(a,b\in [-1,1]\) such that \(\text{supp}(\tilde{t}-Pt_0) = [a,b]\) and:
	\begin{itemize}
		\item \(\displaystyle u<g, \qquad g'< u' <0, \qquad \text{on }\left]-1,a\right[,\)
		\item \(\displaystyle u = g,\qquad\text{on }\left]a,b\right[,\)
		\item \(\displaystyle u < g, \qquad 0 < u' < g',\qquad \text{on }\left]b,1\right[.\)
	\end{itemize}
	In particular, \(\| u'\|_{L^\infty(-1,1)} = \| g'\|_{L^\infty(-1,1)}\).
\end{prop}

\noindent\textbf{Proof.}
Let \((\tilde{t},u)\in \cup_{p>1}L^p(-1,1)\times \cup_{p>1} W^{1,p}(-1,1)\) be an arbitrary solution of problem~\(\mathscr{P}_{\!\rm o}\).  The pair \((\tilde{t},v)\) obviously fulfils all the statements of problem~\(\mathscr{P}_{\!\rm a}\) except maybe:
\begin{equation}
	\label{eq:anotherExpectedConcl}
	\biggl(v(x)-\int_{-1}^x \frac{g'e^{-\bar{\tau}}}{\scriptstyle\sqrt{1+f^2}}\,{\rm d}s\biggr)\,\Bigl(\tilde{t}(x)-Pt_0(x)\Bigr)= 0 .
\end{equation}
Noting that the proof of Corollary~\ref{thm:convexIndentor} uses only the fact that \((\tilde{t},u)\) belongs to \(\cup_{p>1}L^p(-1,1)\times \cup_{p>1} W^{1,p}(-1,1)\) and not that the friction coefficient \(f\) is constant, we can apply its conclusion.  In particular, \(u'\in L^\infty(-1,1)\), \(\text{supp}\,(\tilde{t}-Pt_0)\) is an interval and its bounds \(a,b \in [-1,1]\) are such that \(u'-g'\geq 0\) on \(\left]-1,a\right[\), \(u'-g'=0\) on \(\left]a,b\right[\) and \(u'-g'\leq 0\) on \(\left]b,1\right[\).  This entails that the function:
\[
v(x)-\int_{-1}^x \frac{g'e^{-\bar{\tau}}}{\scriptstyle\sqrt{1+f^2}}\,{\rm d}s,
\]
is non-decreasing on \(\left]-1,a\right[\), is constant on \(\left]a,b\right[\) and is non-increasing on \(\left]b,1\right[\).  As its supremum is zero, by construction, this function actually vanishes on \(\left]a,b\right[\).  This is sufficient to say that identity~\eqref{eq:anotherExpectedConcl} holds and that the pair \((\tilde{t},v)\) solves problem~\(\mathscr{P}_{\!\rm a}\).

Reciprocally, let \((\tilde{t},v)\in \cup_{p>1}L^p(-1,1)\times \cup_{p>1} W^{1,p}(-1,1)\) be an arbitrary solution of problem~\(\mathscr{P}_{\!\rm a}\).  We have \(u' = v'e^{\bar{\tau}}{\scriptstyle\sqrt{1+f^2}} \in \cup_{p>1}L^p(-1,1)\).  The pair \((\tilde{t},u)\) obviously fulfils all the statements of problem~\(\mathscr{P}_{\!\rm o}\) except maybe:
\begin{equation}
	\label{eq:expectedConcl}
	\Bigl(u(x)-g(x)\Bigr)\,\Bigl(\tilde{t}(x)-Pt_0(x)\Bigr)= 0 .
\end{equation}
To prove this, we are going to prove that the set of those \(x\in \left]-1,1\right[\) such that \(v(x)=\varphi(x):={\displaystyle\int_{-1}^x} \frac{g'e^{-\bar{\tau}}}{\scriptstyle\sqrt{1+f^2}}\,{\rm d}s\) is an interval.  Consider \(x_0\in\left]-1,1\right[\) such that \(v(x_0)-\varphi(x_0)<0\) and let \(\left]a,b\right[\) be the largest interval containing \(x_0\) (\(x_0\in\left]a,b\right[\subset\left]-1,1\right[\)) in which this continuous function is negative.  As \(\tilde{t}-Pt_0\) vanishes on \(\left]a,b\right[\), we have, using Theorem~\ref{thm:Riesz1} of Appendix~A:
\[
\forall x\in \left]a,b\right[,\qquad u'(x) = \frac{1}{\pi}\int_{-1}^a \frac{\tilde{t}(s)-Pt_0(s)}{x-s}\,{\rm d}s + \frac{1}{\pi}\int_b^1 \frac{\tilde{t}(s)-Pt_0(s)}{x-s}\,{\rm d}s,
\]
which entails that \(u'\) is decreasing on \(\left]a,b\right[\), so that \(u'-g'\) is also decreasing on \(\left]a,b\right[\) (as \(g\) is convex).  If \(-1<a\) and \(b<1\), we would have, on the one hand \((v-\varphi)(a)=0=(v-\varphi)(b)\).  On the other hand, if \((u'-g')(a+)>0\), then the function \(v'-\varphi'\) would be positive on a right neighbourhood \(\left]a,a+\epsilon\right[\) of \(a\), which is impossible since \(v-\varphi\leq 0\), so we must have \((u'-g')(a+)\leq 0\).  Similarly \((u'-g')(b-)\geq 0\).  But, this implies that \(u'-g'\) vanishes identically on \(\left]a,b\right[\) and the same must be therefore true of \(v'-\varphi'\), which yields a contradiction.  Hence, either \(a=-1\) or \(b=1\), which implies that \((v-\varphi)^{-1}(\{0\})\) is an interval.  As previously, the function \(\tilde{t}-Pt_0\) vanishes outside this interval and the function \(u'-g'\) is decreasing outside this interval.  Since \(u'-g'=0\) on this interval, the constant value taken by \(u-g\) on that interval is also its maximum, which is zero by construction.  Finally, we have proved identity~\eqref{eq:expectedConcl}. \qed

\bigskip
The motivation for defining problem~\(\mathscr{P}_{\!\rm a}\) lies in the fact that the underlying linear operator enjoys the same good properties as those of the underlying linear operator of the case of homogeneous friction (which the underlying linear operator of problem~\(\mathscr{P}_{\!\rm o}\) does not enjoy, as was seen in the preceding section).  These good properties are summarized in the following theorem, whose proof is postponed to Section~\ref{sec:proofConvexHetero}.

\begin{theo}
	\label{thm:heartConvex}
	Let \(f\in W^{1,\infty}(-1,1)\).  Then, for all \(\tilde{t}\in H_0^*:=\{\hat{t}\in H^{-1/2}(-1,1)|\langle\hat{t},1\rangle=0\}\), the identity:
	\[
	\frac{e^{-\bar{\tau}}}{\scriptstyle\sqrt{1+f^2}}\biggl(- \frac{1}{\pi}\,\text{\rm pv}\frac{1}{x} * \tilde{t} + f\,\tilde{t}\biggr) = v',\qquad \text{in }\left]-1,1\right[,
	\]
	defines a unique \(\mathring{v}=\bar{A}\tilde{t}\) in \(H_0:=H^{1/2}(-1,1)/\mathbb{R}\).  The mapping \(\bar{A}:H_0^*\rightarrow H_0\) is continuous and coercive, so that it is actually an isomorphism.  Its inverse mapping \(\bar{A}^{-1}:H_0 \rightarrow H_0^*\) is also continuous, coercive, and enjoys, in addition, the T-monotonicity property:
	\[
	\forall v\in H^{1/2}(-1,1),\qquad \bigl\langle \bar{A}^{-1}\mathring{v},v^+\bigr\rangle \geq 0,
	\]
	where the equality is achieved if and only if \(v^+\) is constant.
\end{theo}

The proof of theorem~\ref{thm:heartConvex} is postponed to Section~\ref{sec:proofConvexHetero}.  We are now going to state its consequences, which are roughly, that the contact problem with heterogeneous friction and convex indentor behaves like the contact problem with homogeneous friction.

\begin{coro}
	\label{thm:exUniqConvex}
	Let \(f\in W^{1,\infty}\) be an arbitrary Lipschitz-continuous friction coefficient, \(g\in W^{1,\infty}(-1,1)\) a \emph{convex} indentor shape and \(P>0\).  Then, there exists a unique \((\tilde{t},v)\in \cup_{p>1}L^p(-1,1)\times \cup_{p>1}W^{1,p}(-1,1)\) that solves problem~\(\mathscr{P}_{\!\rm a}\).
\end{coro}

\noindent\textbf{Proof.} As \(\bar{A}:H_0^* \rightarrow H_0\) is continuous and coercive, thanks to Theorem~\ref{thm:heartConvex}, the Lions-Stampacchia Theorem yields a unique pair \((\tilde{t},v)\in H^{-1/2}(-1,1)\times H^{1/2}(-1,1)\) satisfying all the statements of problem~\(\mathscr{P}_{\!\rm a}\) (with \(\int_{-1}^1\tilde{t}\,{\rm d}x=0\) replaced by \(\langle\tilde{t},1\rangle=0\)).  In addition, by the same reasoning as in the proof of Theorem~\ref{thm:LewyStampacchia}, the T-monotonicity of \(\bar{A}^{-1}\) entails the Lewy-Stampacchia inequality:
\begin{equation}
	\label{eq:LewyStampacchiaInequality}
	\min\biggl\{Pt_0,\bar{A}^{-1}\int_{-1}^x g'e^{-\bar{\tau}}/{\scriptstyle\sqrt{1+f^2}}\,{\rm d}s\biggr\}\leq \tilde{t} \leq Pt_0.
\end{equation}

Note that \(\bar{A}^{-1}\int_{-1}^x g'e^{-\bar{\tau}}/{\scriptstyle\sqrt{1+f^2}}\,{\rm d}s\) denotes the unique solution \(t\in H_0^*\) of the equation:
\[
\frac{e^{-\bar{\tau}}}{\scriptstyle\sqrt{1+f^2}}\biggl(- \frac{1}{\pi}\,\text{\rm pv}\frac{1}{x} * t + f\,t\biggr) = \frac{g'e^{-\bar{\tau}}}{\scriptstyle\sqrt{1+f^2}},\qquad \text{in }\mathscr{D}'(\left]-1,1\right[),
\]
that is, of the non-homogeneous Carleman equation with Lipschitz coefficient:
\[
- \frac{1}{\pi}\,\text{\rm pv}\frac{1}{x} * t + f\,t = g', \qquad \text{in }\mathscr{D}'(\left]-1,1\right[).
\]
It is given by:
\[
t(x) = \frac{f\,g'(x)}{1+f^2}+t_0(x)\, \Biggl\{\text{pv}\frac{1}{x} * \biggl[\sqrt{1-x^2}e^{-\tau(x)}g'(x)/\sqrt{1+f^2}\biggr]\Biggr\},
\]
where \(t_0\) is the function defined in Proposition~\ref{thm:contactZoneIncreasing} (for a proof, see~\cite[Section 4-4]{tricomi} or \cite[theorem 13]{bibiJiri}).  As \(e^{-\tau(x)}g'(x)\in L^\infty(-1,1)\), the Hilbert transform in the above formula is in \(L^p\), for all \(p\in \left]1,+\infty\right[\) by Theorem~\ref{thm:Riesz1} in Appendix~A. As \(t_0\in \cup_{p>1} L^p(-1,1)\), the same is true of \(t=\bar{A}^{-1}\int_{-1}^x g'e^{-\bar{\tau}}/{\scriptstyle\sqrt{1+f^2}}\,{\rm d}s\).  Going back to the Lewy-Stampacchia inequality~\eqref{eq:LewyStampacchiaInequality}, we can conclude \(\tilde{t}\in \cup_{p>1} L^p(-1,1)\).  As the same conclusion applies to:
\[
v' = \frac{e^{-\bar{\tau}}}{\scriptstyle\sqrt{1+f^2}}\biggl(- \frac{1}{\pi}\,\text{\rm pv}\frac{1}{x} * \tilde{t} + f\,\tilde{t}\biggr),
\]
the proof is complete.\qed

\bigskip
Note that the T-monotonicity property of the operator \(\bar{A}^{-1}\) yields incidentally that the contact zone \(\mathscr{C}_P\) (which is an interval as the indentor is assumed to be convex) is a non-decreasing function of \(P\), as in Corollary~\ref{thm:contactZoneIncreasing}.

We are now able to generalize the homogenization result (Theorem~\ref{thm:Homog}) already obtained for the particular case of the flat indentor to the case of an arbitrary convex Lipschitz-continuous (\(g\in W^{1,\infty}(-1,1)\)) indentor.  The following theorem requires that the oscillating friction coefficient \(f_n\) is Lipschitz-continuous.

\begin{theo}
	\label{thm:homogConvex}
	Let \(P>0\) and \(g\in W^{1,\infty}(-1,1)\) be convex.  We are given a period \(p\in W^{1,\infty}(-1,1)\), such that \(p(-1)=p(1)\) which is extended by $2$-periodicity to the whole line.  Let \(f_n:\left]-1,1\right[\rightarrow\mathbb{R}\) be the \(2/n\)-periodic piecewise Lipschitz-continuous function defined on \(\left]-1,1\right[\) by:
	\[
	f_n(x) := p(nx),
	\]
	and let \(\tilde{t}_n\in \cup_{p>1}L^p(-1,1)\) be the unique solution of the contact problem (either problem~\(\mathscr{P}_{\!\rm o}\) or problem~\(\mathscr{P}_{\!\rm a}\))	provided by Corollary~\ref{thm:exUniqConvex}.  Let also \(f_\textrm{\rm eff}\) be the constant:
	\[
	f_\textrm{\rm eff} := \tan \bigl\langle \arctan f_1 \bigr\rangle,
	\]
	where \(\langle \cdot\rangle\) denotes the average, that is, \(\langle h\rangle := (1/2)\int_{-1}^1 h\,{\rm d}x\), and let \(\tilde{t}_\textrm{\rm eff}\) be the unique solution of the contact problem associated with the constant \(f_\textrm{\rm eff}\).
	
	Then, the sequences \((\tilde{t}_n)\) and \((f_n\tilde{t}_n)\) converge weakly-* in \(\mathscr{M}([-1,1])\), respectively towards \(\tilde{t}_\textrm{\rm eff}\) and \(f_\textrm{\rm eff}\tilde{t}_\textrm{\rm eff}\).	In particular, the total tangential contact force \(-\int_{-1}^{1}f_n(\tilde{t}_n-Pt_n^0)\,{\rm d}x\) converges towards \(f_\textrm{\rm eff} P\).
\end{theo}

The proof of Theorem~\ref{thm:homogConvex} is postponed to Section~\ref{sec:proofConvexHetero}.  Note that the type of convergence which is proved for the sequences \((\tilde{t}_n)\) and \((f_n\tilde{t}_n)\) is slightly weaker than the convergence that was proved in the case where the indentor is a rigid flat punch: weak-* convergence in \(\mathscr{M}([-1,1])\) instead of weak convergence in \(L^p(-1,1)\) with \(p\) in a range of values larger than 1.

\section{Proofs of the stated results}

\label{sec:proofs}

\subsection{Proof of the results for the case of the rigid flat punch}

\label{sec:proofFlat}

First, we give a proof of the fact that, in the case of the flat punch, any solution of the contact problem achieves active contact everywhere.

\begin{prop}
	\label{thm:UnilLin}
	Let \(P>0\) and \(f:\left]-1,1\right[\rightarrow\mathbb{R}\) be a given piecewise Lipschitz-continuous function.  Let \((t,u)\in \cup_{p>1}L^p(-1,1)\times \cup_{p>1}W^{1,p}(-1,1)\) be such that:
	\[
	\begin{aligned}
		\bullet \quad & - \frac{1}{\pi}\,\text{\rm pv}\frac{1}{x} * t + f\,t = u',\qquad \text{in }\mathscr{D}'(\left]-1,1\right[), \\[1.2ex]
		\bullet \quad & u \leq 0,\qquad t\leq 0,\qquad u\,t= 0 ,\qquad \text{a.e. in }\left]-1,1\right[, \\[1.2ex]
		\bullet \quad & \int_{-1}^{1} t\,{\rm d}x=-P<0. 
	\end{aligned}
	\]
	Then, \(u\) vanishes identically on \(\left]-1,1\right[\).
\end{prop}

\noindent\textbf{Proof.} Assume that there exists \(x_0\in\left]-1,1\right[\) such that \(u(x_0)<0\).  As \(u\) is continuous, we can consider the largest interval containing \(x_0\) (\(x_0\in\left]a,b\right[\subset\left]-1,1\right[\)) such that \(u\) is negative on \(\left]a,b\right[\).  As \(u\) cannot be negative all over \(\left]-1,1\right[\) because otherwise \(t\) would vanish identically, contradicting \(\int_{-1}^{1} t(x)\,{\rm d}x=-P<0\), only the three following cases can happen:
\begin{enumerate}
	\item \(a=-1\) and \(-1<b<1\), in which case \(u(b)=0\),
	\item \(b=1\) and \(-1<a<1\), in which case \(u(a)=0\),
	\item \(-1<a<b<1\), in which case \(u(a)=u(b)=0\).
\end{enumerate}
In any case, \(t\) must vanish on \(\left]a,b\right[\), so that using Theorem~\ref{thm:Riesz1} in Appendix~A:
\[
\forall x\in \left]a,b\right[, \qquad u'(x)=\frac{1}{\pi}\int_{-1}^a \frac{t(s)}{s-x}\,{\rm d}s + \frac{1}{\pi}\int_b^1 \frac{t(s)}{s-x}\,{\rm d}s.
\]
In case 1, the above identity implies that \(u'\) is negative on \(\left]-1,b\right[\), which contradicts \(u\leq0\).  In case 2, it implies that \(u'\) is positive on \(\left]a,1\right[\), which also contradicts \(u\leq 0\).  In case 3, it implies that \(u'\) must be non-increasing all over \(\left]a,b\right[\), so that \(u\) must be concave on \(\left]a,b\right[\).  As it is also nonpositive and \(u(a)=u(b)=0\), it must vanish identically which also yields contradiction.  Finally, there is no \(x_0\in\left]-1,1\right[\) such that \(u(x_0)<0\). \qed

\bigskip
The next proposition provides a solution of constant sign for the homogeneous Carleman equation.  This solution is classical in the case where \(f\) is Lipschitz-continuous (see, for example, \cite[Section 4-4]{tricomi}).  The proof is extended here to the case where \(f\) is piecewise Lipschitz-continuous and a Green's formula is established in addition.  This will turn out to be crucial for the homogenization analysis.

\begin{prop}
	\label{thm:defTzero}
	Let \(f:\left]-1,1\right[\rightarrow\mathbb{R}\) be a piecewise Lipschitz-continuous function that we extend by zero on \(\mathbb{R}\setminus\left]-1,1\right[\).  We denote by \(\tau:=-\frac{1}{\pi}\,\text{\rm pv}\frac{1}{x}*\arctan f\), the Hilbert transform of the function \(\arctan f\).  Let \(t_0\) be the function defined by:
	\[
	t_0(x):= \frac{e^{\tau(x)}}{\pi\sqrt{1-x^2}\sqrt{1+f^2(x)}},\qquad \text{for x}\in \left]-1,1\right[,
	\]
	and extended by zero on \(\mathbb{R}\setminus\left]-1,1\right[\).  Then, \(t_0\in \cup_{p>1}L^p(-1,1)\), is obviously positive on \(\left]-1,1\right[\), has total mass \(\int_{-1}^1t_0\,{\rm d}x=1\), and its Hilbert transform is given by:
	\[
	-\frac{1}{\pi}\,\text{\rm pv}\frac{1}{x} * t_0 = \left|
	\begin{array}{ll}
		\displaystyle - \frac{f\,e^{\tau}}{\pi\sqrt{1-x^2}\sqrt{1+f^2}} , & \text{ in }\left]-1,1\right[,\\[3.0ex]
		\displaystyle - \frac{\text{\rm sgn}(x)\,e^{\tau}}{\pi\sqrt{x^2-1}} , & \text{ in }\mathbb{R}\setminus\left]-1,1\right[.
	\end{array}
	\right.
	\]
	In particular, \(t_0\) solves the homogeneous Carleman equation:
	\[
	- \frac{1}{\pi}\,\text{\rm pv}\frac{1}{x} * t_0 + f\,t_0 = 0,\qquad \text{in }\left]-1,1\right[.
	\]
	In addition, \(t_0\) can be represented by use of the Green's formula:
	\begin{multline}
		\forall \varphi\in C_c^\infty(\overline{\Pi^+}),\qquad \int_{-1}^{1}t_0(x)\,\varphi(x,0)\,{\rm d}x = \frac{1}{2\pi}\int_{\Pi^+} \Re \Bigl(e^{\phi^+(x+iy)}-e^{\phi^-(x+iy)}\Bigr)\,\frac{\partial\varphi}{\partial x}(x,y) 
		\\
		-\Im \Bigl(e^{\phi^+(x+iy)}-e^{\phi^-(x+iy)}\Bigr)\,\frac{\partial\varphi}{\partial y}(x,y)\,{\rm d}x\,{\rm d}y,
		\label{eq:expectedGreenFormula}
	\end{multline}
	where \(\Pi^+\) denotes the open upper Euclidean plane, \(\overline{\Pi^+}\) its closure and \(\phi^\pm\) are the holomorphic complex functions defined on \(\Pi^+\) by:
	\[
	\phi^\pm(x+iy) := \frac{1}{\pi} \int_{-1}^1 \frac{\pm\pi/2+\arctan f(s)}{s-(x+iy)}\,{\rm d}s.
	\]
	Finally, the following identity holds true:
	\[
	\forall (x,y)\in\Pi^+, \qquad \frac{1}{\pi}\int_{-1}^1 \frac{t_0(s)}{s-(x+iy)}\,{\rm d}s = \frac{1}{2\pi} \Bigl(e^{\phi^+(x+iy)}-e^{\phi^-(x+iy)}\Bigr).
	\]
\end{prop}

\noindent\textbf{Proof.} The strategy of proof to solve the Carleman equation is as follows.  We know that if \(\phi(z)\) is a holomorphic function on \(\Pi^+\) admitting the trace \(\phi(x+i0)=\rho(x)+i\theta(x)\) on \(\partial\Pi^+\), then \(e^{\phi(z)}\) is holomorphic on \(\Pi^+\) with trace \(e^{\rho(x)}\,\cos\theta(x) + ie^{\rho(x)}\, \sin\theta(x)\).  If, in addition, this trace belongs to \(L^p(\mathbb{R})\), then \(e^{\rho(x)}\,\cos\theta(x)\) will be the Hilbert transform of \(e^{\rho(x)}\,\sin\theta(x)\) (see Appendix~A).  Looking for a solution of the Carleman equation of the form \(t=e^{\rho(x)}\,\sin\theta(x)\), we are led to solve \(f=-\mathcal{H}[t]/t=-1/\tan \theta\), obtaining formally the equation \(\tan \theta(x)=-1/f(x)\) for some unknown function \(\theta\), with support in \([-1,1]\).  In order to make this argument rigorous, we set:
\[
\phi^\pm(z) := \frac{1}{\pi} \int_{-1}^1 \frac{\pm\pi/2+\arctan f(s)}{s-z}\,{\rm d}s,
\]
which by use of Theorem~\ref{thm:Riesz2} in Appendix~A satisfies:
\[
\text{for a.a. }x\in\mathbb{R},\qquad \lim_{y\rightarrow 0+}\phi^\pm(x+iy) = \log\sqrt{\biggl|\frac{1\mp x}{1\pm x}\biggr|}+\tau(x) +i\biggl(\pm \frac{\pi}{2}\chi_{\left]-1,1\right[}(x)+\arctan f(x)\biggr) 
\]
where \(\chi_{\left]-1,1\right[}\) is the characteristic function of \(\left]-1,1\right[\) whose Hilbert transform is readily calculated as \((1/\pi)\log |(1-x)/(1+x)|\), using Theorem~\ref{thm:Riesz1} in Appendix~A.  Hence, the function:
\[
\psi^\pm(z) := e^{\phi^\pm(z)} - 1,
\]
is holomorphic in \(\Pi^+\) and satisfies:
\[
\text{for a.a. }x\in\left]-1,1\right[,\quad \lim_{y\rightarrow 0+}\psi^\pm(x+iy) = \mp f(x)\,\sqrt{\frac{1\mp x}{1\pm x}}\,\frac{e^{\tau(x)}}{\sqrt{1+f^2(x)}} - 1 \pm i\,\sqrt{\frac{1\mp x}{1\pm x}}\,\frac{e^{\tau(x)}}{\sqrt{1+f^2(x)}},
\]
and:
\[
\text{for a.a. }x\in\mathbb{R}\setminus[-1,1],\qquad \lim_{y\rightarrow 0+}\psi^\pm(x+iy) = \sqrt{\biggl|\frac{1\mp x}{1\pm x}\biggr|}\,e^{\tau(x)}-1.
\]
We are now going to check that this trace \(\psi^\pm(x+i0)\) belongs to \(L^p(\mathbb{R})\) for some \(p\in\left]1,+\infty\right[\).  First, note that:
\[
\forall x\in \mathbb{R}\setminus[-1,1],\qquad \tau(x) = \frac{1}{\pi}\int_{-1}^1 \frac{\arctan f(s)}{s-x}\,{\rm d}s, 
\]
so that \(\tau\) is \(C^\infty\) on \(\mathbb{R}\setminus[-1,1]\) and \(\tau(x)=O(1/x)\) as \(|x|\rightarrow\infty\).  This entails that \(|\psi^\pm(x+i0)|^p\) is integrable on a neighbourhood of infinity for all \(p>1\).  As \(f\) is piecewise Lipschitz-continuous with support in \([-1,1]\), the function \(\arctan f\) is the sum of a Lipschitz-continuous function on \(\mathbb{R}\) and a piecewise constant function.  This entails that its Hilbert transform \(\tau\) is continuous on \(\mathbb{R}\) except possibly at a finite number of points \(x_i\) (namely the discontinuity points of \(f\) in \([-1,1]\)) where \(\tau\) has the singularity:
\begin{equation}
	\label{eq:estSingul}
	\tau(x) \sim \frac{\arctan f(x_i-) - \arctan f(x_i+)}{\pi}\,\log |x-x_i|, \qquad \text{as }x\rightarrow x_i.
\end{equation}
All in all, the function \(\sqrt{|(1\mp x)/(1\pm x)|}e^{\tau(x)}\) is continuous on \(\mathbb{R}\) except possibly at a finite number of points where it may have a power singularity which belongs to \(L^p\) for some \(p>1\).  Hence, \(\psi^\pm\) is holomorphic on \(\Pi^+\) and \(\psi^\pm(x+i0)\in L^p(\mathbb{R})\) for some \(p>1\).  In addition, \(\phi^\pm\) goes to zero at infinity and the same is therefore true of \(\psi^\pm\).  Hence, Corollary~\ref{thm:noteTricomi} in Appendix~A yields:
\[
-\frac{1}{\pi}\,\text{\rm pv}\frac{1}{x}*\Im \psi^\pm(x+i0) = \Re \psi^\pm(x+i0), \qquad \text{on }\mathbb{R}.
\]
Setting:
\[
t^\pm(x) := \pm\Im \psi^\pm(x+i0) = \sqrt{\biggl|\frac{1\mp x}{1\pm x}\biggr|}\,\frac{e^{\tau(x)}}{\sqrt{1+f^2(x)}}\,\chi_{\left]-1,1\right[}(x),
\]
we get:
\begin{equation}
	\label{eq:eqCarlpm}
	-\frac{1}{\pi}\,\text{pv}\frac{1}{x}*t^\pm + f\,t^\pm = \mp 1, \qquad \text{in }\left]-1,1\right[,
\end{equation}
so that \(t_0:=(t^++t^-)/(2\pi)\) solves the homogeneous Carleman equation on \(\left]-1,1\right[\).  The last identity of the Proposition is now a direct consequence of the application of Theorem~\ref{thm:Riesz2} in Appendix~A to the holomorphic functions \(\psi^+\) and \(\psi^-\).

We now turn to the proof of the representation of \(t_0\) by means of the Green's formula.  The functions \(h^\pm\) defined on \(\Pi^+\) by:
\[
h^\pm(x,y) := \frac{1}{\pi}\int_{-1}^1 t^\pm(s) \,\log\sqrt{(x-s)^2+y^2}\,{\rm d}s,
\]
are clearly harmonic on \(\Pi^+\) (as \(\log\sqrt{x^2+y^2}\) is) and, by Theorem~\ref{thm:Riesz2} of Appendix~A, we have:
\[
\pm\psi^\pm(x+iy) = \frac{1}{\pi}\int_{-1}^1 \frac{t^\pm(s)\,{\rm d}s}{s-(x+iy)} = -\frac{\partial h^\pm}{\partial x}(x,y) + i \frac{\partial h^\pm}{\partial y}(x,y),
\]
on \(\Pi^+\).  By this and Green's formula:
\[
\forall \varphi\in C_c^\infty(\overline{\Pi^+}),\qquad \int_{-1}^1 \frac{\partial h^\pm}{\partial y}(x,0)\,\varphi(x,0)\,{\rm d}x = -\int_{\Pi^+}\nabla h^\pm\cdot\nabla\varphi\:{\rm d}x\,{\rm d}y,
\]
we finally get formula~\eqref{eq:expectedGreenFormula}.

Finally, the only thing which remains to prove is the identity:
\[
\int_{-1}^1 t_0(x)\,{\rm d}x = 1.
\]
Set:
\[
I := \int_{-1}^1 \arctan f(x)\,{\rm d}x,
\]
and:
\[
\tilde{\psi}^\pm(z):=ze^{\phi^\pm(z)} - z \pm 1 + I,
\]
so that \(\tilde{\psi}^\pm\) is holomorphic in \(\Pi^+\) and:
\[
\begin{aligned}
	&\text{for a.a. }x\in\left]-1,1\right[,\qquad \lim_{y\rightarrow 0+}\tilde{\psi}^\pm(x+iy) = \mp x\,\sqrt{\frac{1\mp x}{1\pm x}}\,\frac{e^{\tau(x)}}{\sqrt{1+f^2(x)}} \,\Bigl(f(x) - i\Bigr) -x \pm 1 + I,\\
	&\text{for a.a. }x\in\mathbb{R}\setminus[-1,1],\qquad \lim_{y\rightarrow 0+}\tilde{\psi}^\pm(x+iy) = x\sqrt{\biggl|\frac{1\mp x}{1\pm x}\biggr|}\,e^{\tau(x)}-x \pm 1 + I.
\end{aligned}
\]
As:
\[
\sqrt{\biggl|\frac{1\mp x}{1\pm x}\biggr|}\,e^{\tau(x)} = 1 - \frac{I\pm 1}{x} + O(1/x^2),\qquad \text{as }|x|\rightarrow +\infty,
\]
we have \(\tilde{\psi}^\pm(x+i0)\in L^p(\mathbb{R})\) for some \(p>1\).  Using once more Corollary~\ref{thm:noteTricomi} in appendix~A (since \(\psi^\pm\) is bounded at infinity), we have:
\[
-\frac{1}{\pi}\,\text{\rm pv}\frac{1}{x}*\Im \tilde{\psi}^\pm(x+i0) = \Re \tilde{\psi}^\pm(x+i0), \qquad \text{on }\mathbb{R}.
\]
As the restriction of \(\tilde{\psi}^\pm(x+i0)\) to a neighbourhood of \(x=0\) is actually in \(W^{1,p}\) for some \(p>1\) (thanks to the estimate~\eqref{eq:estSingul}), both members of the above identity are continuous functions on a neighbourhood of \(x=0\), so that we can consider this identity at \(x=0\).  The right-hand side at \(x=0\) reduces to \(\pm 1 + I\) and the left-hand side reduces to:
\[
\pm \frac{1}{\pi}\int_{-1}^1\sqrt{\frac{1\mp x}{1\pm x}}\,\frac{e^{\tau(x)}}{\sqrt{1+f^2(x)}}\,{\rm d}x = \pm 1 + I,
\]
so that:
\[
\frac{1}{\pi}\int_{-1}^1 t^\pm(x)\,{\rm d}x = 1 \pm I \qquad \Rightarrow\qquad \int_{-1}^1 t_0(x)\,{\rm d}x = 1. 
\]
\qed

In the case where \(f\) is Lipschitz-continuous, all the solutions in \(\cup_{p>1}L^p(-1,1)\) of the homogeneous Carleman equation are proportional to \(t_0\) (see \cite[Section 4-4]{tricomi}).  In the case where \(f\) is \emph{piecewise} Lipschitz-continuous, these solutions form a finite-dimensional linear space whose dimension may be larger than one.  This fact was already observed in \cite[theorem 14]{bibiJiri}.  The result and the proof is adapted here to the notation of this paper for the sake of completeness.

\begin{prop}
	\label{thm:defTi}
	Let \(f:\left]-1,1\right[\rightarrow\mathbb{R}\) be a piecewise Lipschitz-continuous function.  Let \(n\in \mathbb{N}\) be the total number of those discontinuity points \(x_i\in\left]-1,1\right[\) of \(f\) such that:
	\[
	f(x_i-) > f(x_i+).
	\]
	Then, all the solutions \(t\in \cup_{p>1}L^p(-1,1)\) of the homogeneous Carleman equation:
	\[
	- \frac{1}{\pi}\,\text{\rm pv}\frac{1}{x} * t + f\,t = 0,\qquad \text{in }\left]-1,1\right[,
	\]
	(where the convolution is meant in terms of the extension of \(t\) by zero on \(\mathbb{R}\)), are given by:
	\[
	t(x) = C_0\,t_0(x) + \sum_{i=1}^n \frac{C_i\, t_0(x)}{x-x_i},
	\]
	where the function \(t_0\) was defined in Proposition~\ref{thm:defTzero} and the \(C_i\)'s are arbitrary real constants.  In addition, we have:
	\[
	\forall i\in\{1,\ldots,n\},\qquad\int_{-1}^1 \frac{t_0(x)}{x-x_i}\,{\rm d}x = 0.
	\]
\end{prop}

\noindent\textbf{Proof.} As already noted in the proof of Proposition~\ref{thm:defTzero}, the Hilbert transform \(\tau\) of \(\arctan f\) admits the estimate:
\[
\tau(x) \sim \frac{\arctan f(x_i-) - \arctan f(x_i+)}{\pi}\,\log |x-x_i|, \qquad \text{as }x\rightarrow x_i,
\]
at every discontinuity point \(x_i\) of \(\arctan f\).  Therefore, if \(x_i\) is a discontinuity point such that \(f(x_i-)>f(x_i+)\), then the function \(t_0(x)/(x-x_i)\) belongs to \(\cup_{p>1}L^p(-1,1)\).  Let \(x_i\) be such a discontinuity point.  In view of Proposition~\ref{thm:defTzero}, we have:
\[
f\,t_0 = \frac{1}{\pi}\,\text{pv}\frac{1}{x}*t_0, \qquad \text{a.e. in }\left]-1,1\right[,
\]
from which we get, using Theorem~\ref{thm:Riesz1} in Appendix~A:
\[
\text{for a.a. }x\in\left]-1,1\right[,\qquad f(x)\,t_0(x) = \frac{x-x_i}{\pi}\,\biggl(\text{pv}\frac{1}{x}*\frac{t_0}{x-x_i}\biggr)(x) -\frac{1}{\pi}\int_{-1}^1 \frac{t_0(x)}{x-x_i}\,{\rm d}x,
\]
that is:
\[
- \frac{1}{\pi}\,\text{\rm pv}\frac{1}{x} * \frac{t_0}{x-x_i} + f\,\frac{t_0}{x-x_i} = -\frac{1}{\pi(x-x_i)}\int_{-1}^1 \frac{t_0(x)}{x-x_i}\,{\rm d}x,\qquad \text{a.e. in }\left]-1,1\right[.
\]
By Theorem~\ref{thm:Riesz1} in Appendix~A, the left-hand side of this identity is in \(\cup_{p>1}L^p(-1,1)\) and, in particular, is integrable.  Therefore:
\[
\int_{-1}^1 \frac{t_0(x)}{x-x_i}\,{\rm d}x = 0,
\]
and the function \(t_0(x)/(x-x_i)\) solves the homogeneous Carleman equation.

There remains only to prove that any solution in \(\cup_{p>1}L^p(-1,1)\) of the homogeneous Carleman equation is a linear combination of \(t_0\) and the \(t_0/(x-x_i)\) (for \(i\in \{1,\ldots n\}\)).  Introducing the holomorphic functions on \(\Pi^+\) defined by:
\[
\phi(z):=-\frac{1}{\pi}\int_{-1}^1 \frac{\pi/2+\arctan f(s)}{s-z}\,{\rm d}s,\qquad \psi(z):=e^{\phi(z)}-1,
\]
and mimicking the proof of Proposition~\ref{thm:defTzero} (or simply applying formula~\eqref{eq:eqCarlpm} after changing \(f\) into \(-f\)), we obtain:
\begin{equation}
	\label{eq:Flav}
	-\frac{1}{\pi}\,\text{pv}\frac{1}{x}*\biggl\{\sqrt{\biggl|\frac{1+x}{1-x}\biggr|}\frac{e^{-\tau(x)}}{\sqrt{1+f^2(x)}}\,\chi_{\left]-1,1\right[}(x)\biggr\} = f(x)\,\sqrt{\frac{1+x}{1-x}}\frac{e^{-\tau(x)}}{\sqrt{1+f^2(x)}} + 1,\qquad \text{a.e. in }\left]-1,1\right[.
\end{equation}
But, by Theorem~\ref{thm:Riesz1} in Appendix~A, for \(g\in \cup_{p>1}L^p(-1,1)\) extended by zero outside \(\left]-1,1\right[\) and \(c\in\mathbb{R}\), we have:
\[
\text{pv}\frac{1}{x}*\Bigl\{(x-c)\,g(x)\Bigr\} = (x-c)\,\Bigl\{\text{pv}\frac{1}{x}*g\Bigr\} - \int_{-1}^1 g(x)\,{\rm d}x.
\]
Applying this inductively with the choices \(c=1,x_i\) together with formula~\eqref{eq:Flav}, we get:
\begin{multline}
	-\frac{1}{\pi}\,\text{pv}\frac{1}{x}*\biggl\{\frac{\sqrt{|1-x^2|}e^{-\tau(x)}\Pi_{i=1}^n(x-x_i)}{\sqrt{1+f^2(x)}}\,\chi_{\left]-1,1\right[}(x)\biggr\} = \mbox{}\\
	\mbox{} = f(x)\,\frac{\sqrt{|1-x^2|}e^{-\tau(x)}\Pi_{i=1}^n(x-x_i)}{\sqrt{1+f^2(x)}} + Q_{n+1}(x),\qquad \text{a.e. in }\left]-1,1\right[,
	\label{eq:lHilbert}
\end{multline}
where \(Q_{n+1}\) denotes a real polynomial of degree \(n+1\).  Now, consider an arbitrary solution \(t\in\cup_{p>1}L^p(-1,1)\) of the homogeneous Carleman equation:
\begin{equation}
	\label{eq:homCarl}
	\mathcal{H}[t] + f\,t = 0,\qquad \text{a.e. in }\left]-1,1\right[,
\end{equation}
where we have used the notation:
\[
\mathcal{H}[t] = - \frac{1}{\pi}\,\text{\rm pv}\frac{1}{x} * t.
\]
Noting that:
\[
l:=\frac{\sqrt{|1-x^2|}e^{-\tau(x)}\Pi_{i=1}^n(x-x_i)}{\sqrt{1+f^2(x)}}\,\chi_{\left]-1,1\right[}(x)\in L^\infty(-1,1),
\]
formula~\eqref{eq:homCarl} entails that:
\[
\mathcal{H}\Bigl[l\,\mathcal{H}[t]\Bigr] + \mathcal{H}[f\,l\,t] = 0, \qquad \text{a.e. in }\mathbb{R}.
\]
Making use of the Poincaré-Bertrand-Tricomi theorem (Corollary~\ref{thm:pbt} in Appendix~A), we get:
\[
\mathcal{H}[l]\,\mathcal{H}[t]-l\,t-\mathcal{H}\Bigl[\mathcal{H}[l]\,t\Bigr] + \mathcal{H}[f\,l\,t] = 0,
\]
which, in view of formula~\eqref{eq:lHilbert}, yields:
\[
f\,l\,\mathcal{H}[t]+Q_{n+1}\,\mathcal{H}[t]-l\,t - \mathcal{H}\bigl[Q_{n+1}\,t\bigr]=0,\qquad \text{a.e. in }\left]-1,1\right[.
\]
Therefore, formula~\eqref{eq:homCarl} yields:
\[
(1+f^2)\,l\,t = Q_{n+1}\,\mathcal{H}[t] - \mathcal{H}\bigl[Q_{n+1}\,t\bigr],\qquad \text{a.e. in }\left]-1,1\right[,
\]
that is:
\[
(1+f^2(x))\,l(x)\,t(x) = -\frac{1}{\pi}\int_{-1}^1 t(s) \,\frac{Q_{n+1}(x)-Q_{n+1}(s)}{x-s}\,{\rm d}s,\qquad \text{a.e. in }\left]-1,1\right[.
\]
But the right-hand side of this identity is a polynomial \(P\) of degree at most \(n\).  Hence, we have proved:
\[
t(x) = \frac{P(x)\,e^{\tau(x)}}{\Pi_{i=1}^n(x-x_i)\,\sqrt{1-x^2}\sqrt{1+f^2(x)}},\qquad \text{a.e. in }\left]-1,1\right[,
\]
which is necessarily a linear combination of \(t_0(x)\) and the \(t_0(x)/(x-x_i)\).
\qed

\bigskip
\noindent\textbf{Proof of Proposition~\ref{thm:ResPb}.} By Proposition~\ref{thm:UnilLin}, any solution of the contact problem is such that \(u=0\) on \(\left]-1,1\right[\), so that \(t\) must actually solve a homogeneous Carleman equation.  Denoting by \(\tau\) the Hilbert transform of \(\arctan f\) (extended by zero outside \(\left]-1,1\right[\)), defining:
\[
t_0(x) := \frac{e^{\tau(x)}}{\pi\sqrt{1-x^2}\sqrt{1+f^2(x)}},\qquad \text{for }x\in \left]-1,1\right[,
\]
and denoting by \(x_i\) those discontinuity points of \(f\) in \(\left]-1,1\right[\) such that \(f(x_i-)>f(x_i+)\), we get:
\[
t(x) = -P\,t_0(x) + \sum_{i=1}^n C_i\,\frac{t_0(x)}{x-x_i},
\]
for some real constants \(C_i\in \mathbb{R}\), thanks to Propositions~\ref{thm:defTzero} and~\ref{thm:defTi}.  Splitting \(\arctan f\) as the sum of a Lipschitz-continuous function on \(\left]-1,1\right[\) and a piecewise constant function, we have the following estimate for \(\tau\) at \(x_i\):
\[
\tau(x) \sim \frac{\arctan f(x_i-) - \arctan f(x_i+)}{\pi}\,\log |x-x_i|, \qquad \text{as }x\rightarrow x_i,
\]
which entails, on the one hand, that \(t_0\) is bounded on a neighbourhood of \(x_i\), and on the other hand, that \(t_0(x)/(x-x_i)\) goes to \(-\infty\) on \(x_i-\) and to \(+\infty\) on \(x_i+\).  Therefore, the condition \(t\leq 0\) entails that all the constants \(C_i\) above must vanish and \(t=-Pt_0\).  In the particular case where \(f\) is a piecewise constant function, then the Hilbert transform \(\tau(x)\) is explicitly obtained from Theorem~\ref{thm:Riesz1} in Appendix~A, yielding the claimed explicit formula for \(t\).
\qed

\bigskip
The next lemma is a minor and usual reformulation of the Riemann-Lebesgue lemma.

\begin{lemm}
	\label{thm:lemmHomog}
	Let \(p\in L^\infty(-1,1)\) be a given bounded function on \(\left]-1,1\right[\) that is extended to \(\mathbb{R}\) by \(2\)-periodicity.  Let \(f_n:\left]-1,1\right[\rightarrow\mathbb{R}\) be the \(2/n\)-periodic function defined on \(\left]-1,1\right[\) by:
	\[
	f_n(x) = p(nx).
	\]
	Then, the sequence \((f_n)\) converges in \(L^\infty(-1,1)\) weak-* towards the constant function \(\langle p \rangle := \frac{1}{2}\int_{-1}^1 p(x)\,{\rm d}x\). 
\end{lemm}

\noindent\textbf{Proof.} As \(\|f_n\|_{L^\infty(-1,1)}=\|p\|_{L^\infty(-1,1)}\), the sequence \((f_n)\) is obviously bounded in \(L^\infty(-1,1)\) and we can extract a subsequence (still denoted by \((f_n)\)) that converges in \(L^\infty(-1,1)\) weak-* towards some limit \(\bar{f}\).  It is readily checked that, for all \(a,b\in \left[-1,1\right]\):
\[
\lim_{n\rightarrow\infty}\int_a^b f_n(x)\,{\rm d}x = \frac{(b-a)}{2}\int_{-1}^1p(x)\,{\rm d}x = (b-a)\,\langle p \rangle.
\]
Therefore, for any piecewise constant function \(c:\left]-1,1\right[\rightarrow\mathbb{R}\), we have:
\[
\lim_{n\rightarrow\infty}\int_{-1}^1c(x)\, f_n(x)\,{\rm d}x = \langle p \rangle\,\int_{-1}^1c(x)\,{\rm d}x .
\]
But, as the piecewise constant functions are dense in \(L^1\), \(\bar{f}\) must be the constant function \(\langle p \rangle\).  As any converging subsequence has the same limit, the whole sequence \((f_n)\) must converge in \(L^\infty(-1,1)\) weak-* towards the constant function \(\langle p\rangle\).\qed

\bigskip
The next technical lemma contains a \(L^p\)-estimate which is crucial for the proof of the homogenization result.

\begin{lemm}
	\label{thm:LpEst}
	Let \((f_n)\) be a sequence of piecewise Lipschitz-continuous functions that is bounded in \(L^\infty(-1,1)\).  We denote by \(\tau_n\) the Hilbert transform of \(\arctan f_n\) and by \((t_n)\) the sequence of functions in \(\cup_{p>1}L^p(-1,1)\) defined by:
	\[
	t_n(x):= \frac{e^{\tau_n(x)}}{\pi\sqrt{1-x^2}\sqrt{1+f_n^2(x)}},\qquad \text{for x}\in \left]-1,1\right[,
	\]
	which is a positive function of total mass \(1\), thanks to Proposition~\ref{thm:defTzero}.  We set:
	\[
	\beta := \frac{1}{\pi} \arctan \Bigl(\sup_{n}\bigl\|f_n\bigr\|_{L^\infty(-1,1)}\Bigr) < \frac{1}{2}.
	\]
	Then, the sequence \((t_n)\) is bounded in \(L^p(-1,1)\), for all \(p\in \left]1,(1/2+\beta)^{-1}\right[\).
\end{lemm}

\noindent\textbf{Proof.} Fix \(p\in \left]1,(1/2+\beta)^{-1}\right[\) arbitrarily.  Note that \(p<2\) and \(2\beta p<1\).  Our goal is to find a majorant of:
\[
\int_{-1}^1 \frac{e^{p\tau_n(x)}}{(1-x^2)^{p/2}}\,{\rm d}x.
\]
We apply Hölder inequality:
\begin{equation}
	\label{eq:HolderIneq}
	\int_{-1}^1 \frac{e^{p\tau_n(x)}}{(1-x^2)^{p/2}}\,{\rm d}x \leq \biggl(\int_{-1}^1\frac{{\rm d}x}{(1-x^2)^{(ps-s+1)/2}}\biggr)^{1/s}\Biggl(\int_{-1}^1\frac{e^{\frac{ps}{s-1}\tau_n(x)}}{\sqrt{1-x^2}}\,{\rm d}x\Biggr)^{(s-1)/s}
\end{equation}
with the choice of an arbitrary \(s\) in the interval:
\[
s\in \biggl] \frac{1}{1-2\beta p},\frac{1}{p-1} \biggr[.
\]
Note that all the \(s\) in that interval are larger than \(1\) and that the interval is non-empty, thanks to the condition \(1<p<(1/2+\beta)^{-1}\).  Such a choice for \(s\) ensures that:
\[
0 <\frac{ps-s+1}{2} < 1,
\]
so that the first integral in the right-hand side of the inequality is finite.  It also ensures that:
\[
q:= \frac{ps}{s-1} \in \left] p,1/(2\beta)\right[.
\]
Finally, inequality~\eqref{eq:HolderIneq} shows that the sequence \((t_n)\) is bounded in \(L^p(-1,1)\), provided that the integral
\[
\int_{-1}^1 \frac{e^{q\tau_n(x)}}{\sqrt{1-x^2}}\,{\rm d}x,
\]
is bounded by above by some constant independent of \(n\).

Define:
\[
\tilde{f}_n(x) := \tan \bigl(q\,\arctan f_n(x)\bigr),
\]
so that \(\tilde{f}_n\) is piecewise Lipschitz continuous and the sequence \((\tilde{f}_n)\) is bounded in \(L^\infty(-1,1)\), thanks to the condition \(q<1/(2\beta)\).  Denoting by \(\tilde{\tau}_n:=q\tau_n\) the Hilbert transform of \(\arctan \tilde{f}_n\),  Proposition~\ref{thm:defTzero} ensures:
\[
\int_{-1}^1 \frac{e^{\tilde{\tau}_n(x)}}{\pi\sqrt{1-x^2}\sqrt{1+\tilde{f}_n^2(x)}}\,{\rm d}x = 1,
\]
which yields:
\[
\int_{-1}^1 \frac{e^{q\tau_n(x)}}{\sqrt{1-x^2}}\,{\rm d}x \leq \pi \sqrt{1+\|\tilde{f}_n^2\|_{L^\infty}} ,
\]
and therefore the claimed estimate. \qed

\bigskip
\noindent\textbf{Proof of Theorem~\ref{thm:Homog}.}  We are given \(p:\left]-1,1\right[\rightarrow\mathbb{R}\) a piecewise Lipschitz-continuous function on \(\left]-1,1\right[\) that is extended to the whole real line by \(2\)-periodicity.  The function \(f_n\) is defined on \(\left]-1,1\right[\) by:
\[
f_n(x) := p(nx),
\]
and extended by zero outside \(\left]-1,1\right[\), so that it is supported in \([-1,1]\) and its restriction to that interval is \(2/n\)-periodic.  Let \(t_n\in \cup_{p>1}L^p(-1,1)\) be the unique nonpositive solution of total mass \(-P\) of the equation:
\begin{equation}
	\label{eq:defTn}
	- \frac{1}{\pi}\,\text{\rm pv}\frac{1}{x} * t_n + f_n\,t_n = 0,\qquad \text{in }\left]-1,1\right[,
\end{equation}
provided by Proposition~\ref{thm:ResPb}.

We set:
\[
\beta := \frac{1}{\pi}\arctan \bigl\|f_n\bigr\|_{L^\infty(-1,1)} = \frac{1}{\pi}\arctan \bigl\|f_1\bigr\|_{L^\infty(-1,1)} < \frac{1}{2},
\]
and pick an arbitrary \(p\in \left]1,(1/2+\beta)^{-1}\right[\).  By Proposition~\ref{thm:defTzero} and lemma~\ref{thm:LpEst}, the sequence \((t_n)\) is bounded in \(L^p(-1,1)\).  As \((f_n)\) is bounded in \(L^\infty(-1,1)\), the sequence \((f_nt_n)\) is also bounded in \(L^p(-1,1)\).  

Extracting a subsequence, if necessary, the sequence \((t_n)\) converges weakly in \(L^p(-1,1)\) towards some limit \(\bar{t}\in L^p(-1,1)\).  We are going to prove that \(\bar{t}=t_{\rm eff}\), where \(t_{\rm eff}\) is defined in the statement of Theorem~\ref{thm:Homog}.  The proof will be based on the representation of \(t_n\) by means of a Green's formula as established in Proposition~\ref{thm:defTzero}:
\begin{multline}
	\forall \varphi\in C_c^\infty(\overline{\Pi^+}),\qquad \int_{-1}^{1}t_n(x)\,\varphi(x,0)\,{\rm d}x = -\frac{P}{2\pi}\int_{\Pi^+} \Re \Bigl(e^{\phi_n^+(x+iy)}-e^{\phi_n^-(x+iy)}\Bigr)\,\frac{\partial\varphi}{\partial x}(x,y)
	\\
	-\Im \Bigl(e^{\phi_n^+(x+iy)}-e^{\phi_n^-(x+iy)}\Bigr)\,\frac{\partial\varphi}{\partial y}(x,y)\,{\rm d}x\,{\rm d}y,
	\label{eq:origGreen}
\end{multline}
where:
\[
\phi_n^\pm(z) := \frac{1}{\pi} \int_{-1}^1 \frac{\pm\pi/2+\arctan f_n(s)}{s-z}\,{\rm d}s,\qquad z=x+iy\in \Pi^+.
\]
The left-hand side of formula~\eqref{eq:origGreen} converges towards \(\langle \bar{t},\varphi\rangle\).  In the right-hand side, the sequences \((\phi_n^\pm(x+iy))\) converge pointwisely in \(\Pi^+\) towards \(\phi_{\rm eff}^\pm(x+iy)\) defined by:
\[
\phi_{\rm eff}^\pm(z) := \frac{1}{\pi} \int_{-1}^1 \frac{\pm\pi/2+\arctan f_{\rm eff}(s)}{s-z}\,{\rm d}s,\qquad z=x+iy\in \Pi^+,
\]
thanks to lemma~\ref{thm:lemmHomog}.  Furthermore, the restriction of \(\phi_n^\pm(x+iy)\) to any compact subset \(K\subset \Pi^+\) of the \emph{open} half-plane \(\Pi^+\) is bounded by a constant depending only on \(K\).  Therefore, we can pass to the limit by dominated convergence in the right-hand side of formula~\eqref{eq:origGreen}, replacing the index \(n\) by the index `eff', \emph{provided} that the support of \(\varphi\) is contained in \(\Pi^+\).  We need to extend the conclusion to the case where the support of \(\varphi\) can be any compact subset of \(\overline{\Pi^+}\).  So let \(K\subset \overline{\Pi^+}\) be an arbitrary compact subset of \(\overline{\Pi^+}\).  Take \(Y>0\) such that \(K\subset \mathbb{R}\times[0,Y]\).  By Proposition~\ref{thm:defTzero}, we have:
\[
\forall (x,y)\in\Pi^+, \qquad \frac{1}{\pi}\int_{-1}^1 \frac{t_n(s)}{s-(x+iy)}\,{\rm d}s = \frac{1}{2\pi}\Bigl(e^{\phi_n^+(x+iy)}-e^{\phi_n^-(x+iy)}\Bigr).
\]
Picking an arbitrary \(y_0>0\), we obtain:
\[
\frac{1}{2\pi}\Bigl\| \Im (e^{\phi_n^+(x+iy_0)}- e^{\phi_n^-(x+iy_0)})\Bigr\|_{L^p(\mathbb{R})}  =  \biggl\| \frac{1}{\pi}\frac{y_0}{x^2+y_0^2}\stackrel{x}{*}t_n(x)\biggr\|_{L^p(\mathbb{R})} \leq \bigl\| t_n \bigr\|_{L^p(\mathbb{R})}.
\]
But, Theorem~\ref{thm:Riesz2} of Appendix~A also yields that \(\Re(e^{\phi_n^+(x+iy_0)}- e^{\phi_n^-(x+iy_0)})\) is the Hilbert transform of \(\Im(e^{\phi_n^+(x+iy_0)}- e^{\phi_n^-(x+iy_0)})\).  As the Hilbert transform is linear continuous in \(L^p(\mathbb{R})\), we have:
\[
\Bigl\| e^{\phi_n^+(x+iy_0)}- e^{\phi_n^-(x+iy_0)}\Bigr\|_{L^p(\mathbb{R})} \leq C \bigl\| t_n \bigr\|_{L^p(\mathbb{R})},
\]
for some real constant \(C\) that depends on \(p\), but not of \(y_0\) and \(n\).  Using Lemma~\ref{thm:LpEst}, we can conclude that there exists a constant \(K>0\), independent of \(y_0\) and \(n\) such that:
\[
\Bigl\| e^{\phi_n^+(x+iy_0)}- e^{\phi_n^-(x+iy_0)}\Bigr\|_{L^p(\mathbb{R})} \leq K.
\]
In particular, the sequence \((e^{\phi_n^+(x+iy)}- e^{\phi_n^-(x+iy)})\) is bounded in \(L^p(\mathbb{R}\times [0,Y];\mathbb{C})\).  Therefore, it has a weakly convergent subsequence.  But the above analysis shows that this weak limit must be \(e^{\phi_{\rm eff}^+(x+iy)}- e^{\phi_{\rm eff}^-(x+iy)}\).

Now, we are able to pass to the limit in formula~\eqref{eq:origGreen} to get, for all \(\varphi\in C_c^\infty(\overline{\Pi^+})\):
\begin{align*}
	\int_{-1}^{1}\bar{t}(x)\,\varphi(x,0)\,{\rm d}x 
	& = -\frac{P}{2\pi}\int_{\Pi^+} \Re \Bigl(e^{\phi_{\rm eff}^+(x+iy)}-e^{\phi_{\rm eff}^-(x+iy)}\Bigr)\,\frac{\partial\varphi}{\partial x}(x,y)
	\\
	& \phantom{\mbox{}=\mbox{}}\qquad\mbox{}-\Im \Bigl(e^{\phi_{\rm eff}^+(x+iy)}-e^{\phi_{\rm eff}^-(x+iy)}\Bigr)\,\frac{\partial\varphi}{\partial y}(x,y)\,{\rm d}x\,{\rm d}y,\\
	& = \int_{-1}^{1}t_{\rm eff}(x)\,\varphi(x,0)\,{\rm d}x.
\end{align*}
All in all, we have proved that the sequence \((t_n)\) converges weakly in \(L^p(-1,1)\) towards \(t_{\rm eff}\).

As the Hilbert transform maps continuously \(L^p(\mathbb{R})\) onto itself, by Theorem~\ref{thm:Riesz1} in Appendix~A, \((1/\pi)\text{pv}\,1/x * t_n\) converges weakly in \(L^p(-1,1)\) towards \((1/\pi)\text{pv}\,1/x * t_{\rm eff}=f_{\rm eff}t_{\rm eff}\).  Hence, the sequence \((f_nt_n)\) converges weakly in \(L^p(-1,1)\) towards the function \(f_{\rm eff}t_{\rm eff}\).\qed

\subsection{Proof of the results for the case of homogeneous friction}

\label{sec:homogeneousFriction}

In the case of the analysis of the contact problem with homogeneous friction, the only result whose proof was postponed here is Theorem~\ref{thm:corodelamort}.  Given \(f\in\mathbb{R}\) and \(\mathring{u}\in H_0:=H^{1/2}(-1,1)/\mathbb{R}\), let \(\tilde{t}=A^{-1}\mathring{u}\) be the unique distribution in \(H_0^*:=\{\hat{t}\in H^{-1/2}(-1,1)\;|\;\langle\hat{t},1\rangle=0\}\) such that:
\[
- \frac{1}{\pi}\,\text{pv}\frac{1}{x} * \tilde{t} + f\,\tilde{t} = \mathring{u}',\qquad \text{in }\mathscr{D}'(\left]-1,1\right[).
\]
Then, the linear operator \(A^{-1}:H_0 \rightarrow H_0^*\) is an isomorphism.  Theorem~\ref{thm:corodelamort} states that it has the following \emph{T-monotonicity} property:
\[
\forall u\in H^{1/2}(-1,1),\qquad \bigl\langle A^{-1}\mathring{u} , u^+\bigr\rangle \geq 0, 
\]
(where \(u^+(x):=\max\{u(x),0\}\)) with equality if and only if \(u\) is a constant function on \(\left]-1,1\right[\).  

\begin{lemm}
	\label{thm:coroR}
	Let \(u\in H^{1/2}(\mathbb{R})\) and \(C\in\mathbb{R}\) be arbitrary.  Then, \({(u+C)}^+-C^+\in H^{1/2}(\mathbb{R})\) and:
	\[
	\biggl\langle \text{\rm pv}\frac{1}{x}*u' \;,\; {(u+C)}^+ - C^+ \biggr\rangle \geq 0,
	\]
	where the equality is achieved if and only if \({(u+C)}^+-C^+=0\).
\end{lemm}

\noindent\textbf{Proof.} We recall that the space \(H^{1/2}(\mathbb{R})\) can be defined as the space of those distributions \(u\) such that:
\[
\|u\|_{L^2}^2 + \int_{\mathbb{R}^2} \Bigl| \frac{u(x)-u(y)}{x-y} \Bigr|^2 \,{\rm d}x\,{\rm d}y < \infty.
\]
Taking \(u\in H^{1/2}(\mathbb{R})\) and \(C\in \mathbb{R}\), we have:
\[
\Bigl|{\bigl(u(x)+C\bigr)}^+-C^+\Bigr| \leq |u(x)|, \qquad \text{and} \qquad
\Bigl| {\bigl(u(x)+C\bigr)}^+ - {\bigl(u(y)+C\bigr)}^+ \Bigr| \leq | u(x) - u(y) |,
\]
which shows that \({(u+C)}^+-C^+\in H^{1/2}(\mathbb{R})\).

Next, we denote by \(\Phi_u\) the Poisson integral of \(u\) (see Proposition~\ref{thm:prolAnal} in Appendix~C) and we define:
\[
\varphi_n(x,y) := \min\{1,n/y\} \Bigl[ {\bigl(\Phi_u(x,y) + C\bigr)}^+ - C^+ \Bigr],
\]
for \((x,y)\in\Pi^+\), so that \(\varphi_n\in H^1(\Pi^+)\) and the sequence \((\nabla \varphi_n)\) converges strongly in \(L^2(\Pi^+)\) towards the limit \(\nabla {\bigl(\Phi_u(x,y) + C\bigr)}^+\) (we recall that \(\varphi\in H^1\Rightarrow\varphi^+\in H^1\), by Stampacchia's theorem \cite[Theorem~1.56, p~79]{Troianiello}).  By~Proposition~\ref{thm:prolAnal}, this entails:
\begin{align*}
	\biggl\langle \text{\rm pv}\frac{1}{x}*u' \;,\; {(u+C)}^+ - C^+ \biggr\rangle 
	& = \int_{\Pi^+}\biggl( \frac{\partial\Phi_u}{\partial x}\frac{\partial}{\partial x}{\bigl(\Phi_u+C\bigr)}^+ + \frac{\partial\Phi_u}{\partial y}\frac{\partial}{\partial y}{\bigl(\Phi_u+C\bigr)}^+
	\biggr){\rm d}x\,{\rm d}y, \\
	& = \int_{\Pi^+}\biggl\{ {\biggl[ \frac{\partial}{\partial x}{\bigl(\Phi_u+C\bigr)}^+ \biggr]}^2 + {\biggl[\frac{\partial}{\partial y}{\bigl(\Phi_u+C\bigr)}^+ \biggr]}^2\biggr\}{\rm d}x\,{\rm d}y,
\end{align*}
thanks to Stampacchia's theorem.  Hence, the duality product is nonnegative.  It only vanishes when \({(\Phi_u+C)}^+\) is a constant, implying that \({(u+C)}^+\) is a constant.  But, as \({(u+C)}^+-C^+\in L^2(\mathbb{R})\), this constant must be equal to \(C^+\). \qed

\bigskip
\noindent\textbf{Proof of Theorem~\ref{thm:corodelamort}.} Pick \(u\in H^{1/2}(-1,1)\) and define \(\tilde{t}:=A^{-1}\mathring{u}\in H_0^*\).  The distribution \(\tilde{t}\in H^{-1/2}(-1,1)\) can be seen as an element of \(H^{-1/2}(\mathbb{R})\) with support in \([-1,1]\).  Its Fourier transform \(\hat{\tilde{t}}\) is therefore a \(C^\infty\) function and, as \(\langle \tilde{t},1\rangle=0\), it satisfies \(\hat{\tilde{t}}(0)=0\).  Next, we define:
\[
\tilde{u} := -\frac{1}{\pi}\,\log|x|*\tilde{t} + \frac{f}{2}\,\text{sgn}(x)*\tilde{t}.
\]
By Proposition~\ref{thm:FourierTransforms} in Appendix~C, we obtain:
\[
\hat{\tilde{u}}(\xi) = \bigl(\text{sgn}(\xi) - if\bigr)\,\frac{\hat{\tilde{t}}(\xi)}{\xi}\qquad \Rightarrow\qquad \bigl|\hat{\tilde{u}}(\xi)\bigr| = \sqrt{1+f^2}\,\frac{|\hat{\tilde{t}}(\xi)|}{|\xi|},
\]
where the `hat' stands for the Fourier transform.  As \(\hat{\tilde{t}}\) is a \(C^\infty\) function which satisfies:
\[
\hat{\tilde{t}}(0) = 0, \qquad \text{and}\qquad \int_{-\infty}^{+\infty} \frac{|\hat{\tilde{t}}(\xi)|^2}{\sqrt{1+\xi^2}}\,{\rm d}\xi < \infty,
\]
we have:
\[
\int_{-\infty}^{+\infty} |\hat{\tilde{u}}(\xi)|^2\,\sqrt{1+\xi^2}\,{\rm d}\xi < \infty,
\]
that is, \(\tilde{u}\in H^{1/2}(\mathbb{R})\), by Proposition~\ref{thm:defSobolevFourier} in Appendix~B.

We have:
\[
-\frac{1}{\pi}\,\text{pv}\frac{1}{x}*\tilde{t} + f\,\tilde{t} = \tilde{u}',\qquad \text{on }\mathbb{R},
\]
which is easily inverted by taking the Hilbert transform and use of Theorem~\ref{thm:Riesz1} in Appendix~A, as:
\[
\tilde{t} = \frac{f}{1+f^2} \,\tilde{u}' + \frac{1}{1+f^2}\,\frac{1}{\pi}\,\text{pv}\frac{1}{x}*\tilde{u}'.
\]
As \(\tilde{u}'=u'\) on \(\left]-1,1\right[\), there exists \(C\in\mathbb{R}\) such that \(u=\tilde{u}+C\) on \(\left]-1,1\right[\).  Hence, using this and the previous identity, it is easily calculated that:
\begin{align}
	\bigl\langle A^{-1}\mathring{u}\;,\; u^+ \bigr\rangle 
	& = \Bigl\langle\tilde{t}\;,\; {(\tilde{u}+C)}^+ - C^+\Bigr\rangle,\nonumber\\
	& = \frac{f}{1+f^2}\Bigl\langle\tilde{u}'\;,\; {(\tilde{u}+C)}^+ - C^+\Bigr\rangle + \frac{1}{\pi(1+f^2)}\Bigl\langle\text{pv}1/x*\tilde{u}'\;,\; {(\tilde{u}+C)}^+ - C^+\Bigr\rangle.
	\label{eq:flaviana}
\end{align}
We are now going to prove that the first term in the right-hand side of identity~\eqref{eq:flaviana} vanishes.  In the particular case where \(w\in C_c^\infty(\mathbb{R})\), we have:
\[
\Bigl\langle w'\;,\; {(w+C)}^+-C^+\Bigr\rangle = \Bigl\langle w'\;,\; {(w+C)}^+\Bigr\rangle =  
\Bigl\langle (w+C)'\;,\; {(w+C)}^+\Bigr\rangle = 0. 
\]
Taking a sequence \((w_n)\) in \(C_c^\infty(\mathbb{R})\) converging strongly in \(H^{1/2}(\mathbb{R})\) towards \(\tilde{u}\), the sequence  \((w_n')\) converges strongly in \(H^{-1/2}(\mathbb{R})\) towards \(\tilde{u}'\) and the sequence \(({(w_n+C)}^+-C^+)\) converges strongly in \(H^{1/2}(\mathbb{R})\) towards \(({(\tilde{u}+C)}^+-C^+)\).  Hence, we have proved:
\[
\Bigl\langle\tilde{u}'\;,\; {(\tilde{u}+C)}^+ - C^+\Bigr\rangle = 0,
\]
and therefore, by identity~\eqref{eq:flaviana}:
\[
\bigl\langle A^{-1}\mathring{u}\;,\; u^+ \bigr\rangle  = \frac{1}{\pi(1+f^2)}\Bigl\langle\text{pv}1/x*\tilde{u}'\;,\; {(\tilde{u}+C)}^+ - C^+\Bigr\rangle,
\]
and the claim is now a straightforward consequence of Lemma~\ref{thm:coroR}. \qed

\subsection{Proof of the existence result for heterogeneous friction and arbitrary indentor}

\label{sec:proofGeneralHetero}

In the case where \(f\) is a piecewise Lipschitz-continuous function on \([-1,1]\) (extended by zero on \(\mathbb{R}\setminus [-1,1]\)) and denoting \(\tau:=-\frac{1}{\pi}\,\text{\rm pv}\frac{1}{x}*\arctan f\), the Hilbert transform of \(\arctan f\), it was proved in Proposition~\ref{thm:defTzero} that the nonnegative function:
\[
t_0(x):= \frac{e^{\tau(x)}}{\pi\sqrt{1-x^2}\sqrt{1+f^2(x)}},\qquad \text{for }x\in \left]-1,1\right[,
\]
and extended by zero on \(\mathbb{R}\setminus\left]-1,1\right[\), belongs to \(\cup_{p>1}L^p(-1,1)\), has total mass \(\int_{-1}^1t_0\,{\rm d}x=1\) and solves the homogeneous Carleman equation:
\[
- \frac{1}{\pi}\,\text{\rm pv}\frac{1}{x} * t_0 + f\,t_0 = 0,\qquad \text{a.e. in }\left]-1,1\right[.
\]
We are now going to prove that all these facts remain true when \(f\) is only supposed to be in \(L^\infty(-1,1)\).

\bigskip
\noindent\textbf{Proof of Proposition~\ref{thm:defTzeroLinfty}.} Pick \(f\in L^\infty(-1,1)\) and an approximation of identity \(\rho_n(x):=n\,\rho(nx)\) based on some nonnegative even \(C^\infty\) function \(\rho\) with compact support and total mass \(1\).  Then, the sequence of functions \((f_n)\) defined by:
\[
f_n(x) := (\rho_n*f)(x)\,\chi_{\left]-1,1\right[}(x),
\]
is a sequence of Lipschitz-continuous functions on \(\left]-1,1\right[\), which is bounded in \(L^\infty(-1,1)\) (\(\|f_n\|_{L^\infty}\leq \|f\|_{L^\infty}\)), and converges towards \(f\) in \(L^\infty(-1,1)\) weak-*.  Set:
\[
t_n^0(x) : = \frac{e^{\tau_n(x)}}{\pi\sqrt{1-x^2}\sqrt{1+f_n^2(x)}},\qquad \text{for }x\in \left]-1,1\right[,
\]
extended by zero on \(\mathbb{R}\setminus\left]-1,1\right[\), where \(\tau_n\) denotes the Hilbert transform of \(\arctan f_n\).  By Proposition~\ref{thm:defTzero}, the function \(t_n^0\) is in \(\cup_{p>1}L^p(-1,1)\) and has total mass \(1\).  By Lemma~\ref{thm:LpEst}, the sequence \((t_n^0)\) is bounded in \(L^p(-1,1)\), for all	\(p\in \left]1,(1/2+\beta)^{-1}\right[\), where:
\[
\beta := \frac{1}{\pi} \arctan \bigl\|f\bigr\|_{L^\infty(-1,1)} < \frac{1}{2}.
\]
Hence, upon extracting a subsequence, the sequence \((t_n^0)\) converges weakly in \(L^p(-1,1)\), for all	\(p\in \left]1,(1/2+\beta)^{-1}\right[\), towards some limit \(\bar{t}\) which has total mass \(1\).  Writing the Green's formula provided by Proposition~\ref{thm:defTzero} for \(t_n^0\) and going to the limit \(n \rightarrow+\infty\) along the lines of the proof of Theorem~\ref{thm:Homog}, we get:
\begin{multline*}
	\forall \varphi\in C_c^\infty(\overline{\Pi^+}),\qquad \int_{-1}^{1}\bar{t}(x)\,\varphi(x,0)\,{\rm d}x = \frac{1}{2\pi}\int_{\Pi^+} \Re \Bigl(e^{\phi^+(x+iy)}-e^{\phi^-(x+iy)}\Bigr)\,\frac{\partial\varphi}{\partial x}(x,y)
	\\
	-\Im \Bigl(e^{\phi^+(x+iy)}-e^{\phi^-(x+iy)}\Bigr)\,\frac{\partial\varphi}{\partial y}(x,y)\,{\rm d}x\,{\rm d}y,
\end{multline*}
where:
\[
\phi^\pm(z) := \frac{1}{\pi} \int_{-1}^1 \frac{\pm\pi/2+\arctan f(s)}{s-z}\,{\rm d}s.
\]
As \(\Im(e^{\phi^+(x+i0)}-e^{\phi^-(x+i0)})/(2\pi)=\bar{t}\in \cup_{p>1} L^p(-1,1)\), this entails by Theorem~\ref{thm:Riesz2} of Appendix~A that:
\[
\forall z\in \Pi^+,\qquad \frac{1}{\pi}\int_{-1}^1  \frac{\bar{t}(s)\,{\rm d}s}{s-z} = \frac{e^{\phi^+(z)}-e^{\phi^-(z)}}{2\pi},
\]
and that \((e^{\phi^+(x+iy)}-e^{\phi^-(x+iy)})/(2\pi)\) converges, as \(y \rightarrow 0+\), for almost all \(x\), towards \(\mathcal{H}[\bar{t}](x)+i\bar{t}(x)\).  But, the pointwise limit of \((e^{\phi^+(x+iy)}-e^{\phi^-(x+iy)})/(2\pi)\) can also be computed from the definition of \(\phi^\pm(x+iy)\) using again Theorem~\ref{thm:Riesz2} of Appendix~A. Its imaginary part is found to be:
\[
\frac{e^{\tau(x)}}{\pi\sqrt{|1-x^2|}\sqrt{1+f^2(x)}}\,\chi_{\left]-1,1\right[}(x) = t_0(x),
\]
and the restriction of its real part to \(\left]-1,1\right[\) is found to equal \(-ft_0\).  Therefore, \(\bar{t}=t_0\) and it solves the homogeneous Carleman equation.\qed

\bigskip
We now turn to the proof of Theorem~\ref{thm:pseudomonot}.  Assume that \(f\in BV([-1,1])\).  The Banach space \(H_{0}^*\cap\mathscr{M}([-1,1])\) is the dual of the Banach space \(H_0+C^0([-1,1])\).  Given an arbitrary \(\tilde{t}\in H_0^*\cap \mathscr{M}([-1,1])\), we have defined \(\mathring{u}=A\tilde{t}\) by the identity:
\[
- \frac{1}{\pi}\,\text{pv}\frac{1}{x} * \tilde{t} + f\,\tilde{t} = \mathring{u}',\qquad \text{in }\mathscr{D}'(\left]-1,1\right[).
\]
As \(\log|x|*\tilde{t}\) is in \(H^{1/2}(-1,1)\) and \(x\mapsto \int_0^x f\tilde{t}\) is a continuous function with bounded variation, the linear operator \(A\) maps \(H_{0}^*\cap\mathscr{M}([-1,1])\) into \(H_0+C^0([-1,1])\cap BV([-1,1])\subset H_0+C^0([-1,1])\).  It was shown in Proposition~\ref{thm:boundedA} that the linear operator \(A:H_{0}^*\cap\mathscr{M}([-1,1])\rightarrow H_0+C^0([-1,1])\) is bounded.

The following lemma contains the result of a calculation which will prove in the sequel to be cornerstone for the analysis of the monotonicity and pseudomonotonicity of the operator \(A\).

\begin{lemm}
	\label{thm:monotLemma}
	Let \(f\in BV([-1,1])\).  For all \(t\in H_{0}^*\cap\mathscr{M}([-1,1])\), we have:
	\[
	\Bigl\langle t,{\textstyle \int_{-1}^x ft}\Bigr\rangle = \frac{1}{2}\int_{-1}^1 f' \biggl(\int_{-1}^x t\biggr)^2,
	\]
	where the integral in the right-hand side stands for the integral of the continuous function \((\int_{-1}^xt)^2\) (which vanishes at $x=-1,1$) with respect to the measure \(f'\).
\end{lemm}

\noindent\textbf{Proof.}
First note that in the case where \(f\) is a constant function, the duality product at stake is zero.  Therefore, in the general case, we have:
\[
\Bigl\langle t,{\textstyle \int_{-1}^x ft}\Bigr\rangle = \Bigl\langle t,{\textstyle \int_{-1}^x (f+m)t}\Bigr\rangle,
\]
where \(m\) is an arbitrary constant.  From now on, we take \(m> -\inf_{[-1,1]}f\) (we recall that any function in \(\mathit{BV}([-1,1])\) is bounded), so that the function \(f+m\) is positive and bounded away from \(0\).  In that case, the function \(1/(f+m)\) is in \(\mathit{BV}([-1,1])\) and its derivative is the measure:
\[
-\frac{f'}{(f(x-)+m)(f(x+)+m)}.
\]
We therefore have:
\begin{align*}
	\Bigl\langle t,{\textstyle \int_{-1}^x (f+m)t}\Bigr\rangle & = \int_{-1}^1 t \biggl(\int_{-1}^x (f+m)t\biggr), \\
	& = \frac{1}{2}\int_{-1}^1 \frac{1}{f(x)+m} \frac{{\rm d}}{{\rm d}x}\biggl(\int_{-1}^x (f+m)t\biggr)^2, \\
	& = \frac{1}{2(f(1-)+m)}\biggl(\int_{-1}^1 (f+m)t\biggr)^2 \\
	& \phantom{=} \mbox{}+ \frac{1}{2}\int_{\left]-1,1\right[}\frac{f'}{(f(x-)+m)(f(x+)+m)}\biggl(\int_{-1}^x (f+m)t\biggr)^2.
\end{align*}
where the second identity has made use of the fact that the countable set of discontinuities of \(f\) is negligible for the measure \(t\) as \(t\) has no atom.  The last identity has been obtained through the integration by part formula for functions with bounded variation.  The first term of the right-hand side of the last identity equals:
\[
\frac{1}{2(f(1-)+m)}\biggl(\int_{-1}^1 ft\biggr)^2
\]
and therefore tends to \(0\) as \(m\) grows to \(+\infty\).  In addition, the function:
\[
x \mapsto \frac{1}{(f(x-)+m)(f(x+)+m)}\biggl(\int_{-1}^x (f+m)t\biggr)^2,
\]
is uniformly bounded with respect to \(m\) as \(m\) goes to infinity and converges pointwisely towards:
\[
x \mapsto \biggl(\int_{-1}^x t\biggr)^2.
\]
By the dominated convergence theorem, we have:
\[
\lim_{m \rightarrow+\infty}\int_{\left]-1,1\right[}\frac{f'}{(f(x-)+m)(f(x+)+m)}\biggl(\int_{-1}^x (f+m)t\biggr)^2 = \int_{-1}^1f'\biggl(\int_{-1}^x t\biggr)^2,
\]
which gives the announced identity. \qed

\begin{coro}
	\label{thm:monot}
	We assume that the function \(f\) is nondecreasing.  The bounded linear operator \(A:H_{0}^*\cap\mathscr{M}([-1,1])\rightarrow  H_0+C^0([-1,1])\) is strictly monotone (see Definition~\ref{thm:defMonot}).
\end{coro}

\noindent\textbf{Proof.}
If \(f\) is nondecreasing, then the measure \(f'\) is nonnegative.  The previous lemma entails, for all nonzero \(t\in H_{0}^*\cap\mathscr{M}([-1,1])\):
\[
\bigl\langle At,t\bigr\rangle = -\frac{1}{\pi} \Bigl\langle \log|x|*t,t\Bigr\rangle + \Bigl\langle t,{\textstyle \int_{-1}^x ft}\Bigr\rangle \geq  -\frac{1}{\pi} \Bigl\langle \log|x|*t,t\Bigr\rangle >0,
\]
where the last inequality relies, once again, on Lemma~\ref{thm:lemmHminusHalf} of Appendix~C.\qed

\bigskip
\noindent\textbf{Proof of Theorem~\ref{thm:pseudomonot}.}  Let \((t_n)\) be a sequence in \(H_{0}^*\cap\mathscr{M}([-1,1])\) that converges weakly-* in \(H_{0}^*\cap\mathscr{M}([-1,1])\) towards a limit \(t\), such that:
\[
\limsup_{n\rightarrow+\infty} \bigl\langle A t_n,t_n-t\bigr\rangle \leq 0.
\]
From the definition of pseudomonotonicity in Appendix~D, we have to prove:
\[
\forall \tau\in H_{0}^*\cap\mathscr{M}([-1,1]),\qquad \bigl\langle A t,t-\tau\bigr\rangle \leq \liminf_{n\rightarrow+\infty} \bigl\langle A t_n,t_n-\tau\bigr\rangle.
\]

\noindent\textbf{Step 1.} \textit{The sequence \((At_n)\) converges weakly in \(H_0+C^0([-1,1])\) towards \(At\).}

We pick an arbitrary \(\tau\in H_{0}^*\cap\mathscr{M}([-1,1])\).  We have to prove that \(\lim_{n \rightarrow+\infty} \langle At_n,\tau\rangle = \langle At,\tau\rangle\).  But:
\[
\bigl\langle At_n,\tau\bigr\rangle =  -\frac{1}{\pi} \Bigl\langle \log|x|*t_n,\tau\Bigr\rangle + \int_{-1}^1\tau\int_{-1}^x ft_n.
\]
The first term of the right-hand side is the \(H^{-1/2}\)-scalar product of \(t_n\) and \(\tau\) (thanks to~\cite[Theorem 3]{bibiJiri}.  As \((t_n)\) converges weakly in \(H^{-1/2}\) towards \(t\), it converges towards \((-1/\pi)\langle  \log|x|*t,\tau\rangle\).  There remains only to prove that \(\int_{-1}^1\tau\int_{-1}^x ft_n\) converges towards \(\int_{-1}^1\tau\int_{-1}^x ft\).  As  \((t_n)\) converges weakly-* in \(\mathscr{M}([-1,1])\) towards \(t\), the sequence \((t_n)\) is bounded in \(\mathscr{M}([-1,1])\).  Therefore, the functions \(x\mapsto \int_{-1}^x ft_n\) are bounded uniformly with respect to \(n\).  Hence, it is sufficient to prove that the sequence of functions \(x\mapsto \int_{-1}^x ft_n\) converges pointwisely towards \(x\mapsto \int_{-1}^x ft\), to obtain the expected conclusion \(\lim_{n \rightarrow\infty}\int_{-1}^1\tau\int_{-1}^x ft_n=\int_{-1}^1\tau\int_{-1}^x ft\) from the dominated convergence theorem.  So, let us prove that the sequence of functions \(x\mapsto \int_{-1}^x ft_n\) converges pointwisely towards \(x\mapsto \int_{-1}^x ft\).  In the particular case where all the \(t_n\), and therefore \(t\), are nonnegative measures (forgetting for a while the condition \(\int_{[-1,1]}t=0\)) and \(x\in \left]-1,1\right[\) is fixed, given an arbitrary \(\varepsilon>0\), we can build a continuous function \(\varphi_\varepsilon:[-1,1]\rightarrow[0,1]\), supported in \([-1,x]\), such that:
\[
(1-\varepsilon)\int_{-1}^{x}t \leq \int_{-1}^{1}\varphi_\varepsilon t \leq \int_{-1}^{x}t,
\]
and a function \(\psi_\varepsilon:[-1,1]\rightarrow[0,1]\) that takes the value \(1\) all over \([-1,x]\), such that:
\[
\int_{-1}^{x}t \leq \int_{-1}^{1}\psi_\varepsilon t \leq (1+\varepsilon) \int_{-1}^{x}t,
\]
where the fact that the measure \(t\) has no atom at \(x\) has been used.  As:
\[
\forall n\in \mathbb{N},\qquad \int_{-1}^{1}\varphi_\varepsilon t_n \leq \int_{-1}^x t_n \leq \int_{-1}^{1}\psi_\varepsilon t_n,
\]
we obtain that the sequence of functions \(x\mapsto \int_{-1}^x t_n\) converges pointwisely towards \(x\mapsto \int_{-1}^x t\).  Hence, provided that the sequence \((t_n)\) is nonnegative, and \(f\) is a stepwise constant function, the sequence of functions \(x\mapsto \int_{-1}^x ft_n\) converges pointwisely towards \(x\mapsto \int_{-1}^x ft\).  As any function with bounded variation in \([-1,1]\) is a uniform limit of a sequence of stepwise constant functions, the conclusion extends to the case where \((t_n)\) is nonnegative, and \(f\in BV([-1,1])\).  Finally, let us drop the hypothesis that \((t_n)\) is nonnegative.  An arbitrary measure \(t_n\) is the difference of its positive and negative parts: \(t_n=t_n^+-t_n^-\).  As \(t_n\) has no atom, the same is true of \(t_n^+\) and \(t_n^-\).  As the sequence \((t_n)\) is weakly-* convergent in \(\mathscr{M}([-1,1])\), it is bounded in \(\mathscr{M}([-1,1])\) and therefore, so are the sequences \((t_n^+)\) and \((t_n^-)\).  Hence, extracting, if necessary, a subsequence, the sequences \((t_n^+)\) and \((t_n^-)\) converge weakly-* in \(\mathscr{M}([-1,1])\) towards limits \(t_+\) and \(t_-\) such that \(t=t_+-t_-\).  Now, applying the previous conclusion to the sequences \((t_n^+)\) and \((t_n^-)\), we have proved the desired conclusion that the sequence of functions \(x\mapsto \int_{-1}^x ft_n\) converges pointwisely towards \(x\mapsto \int_{-1}^x ft\).

\noindent\textbf{Step 2.} \textit{Conclusion.}

We have, by hypothesis, \(\limsup_{n\rightarrow+\infty} \bigl\langle A t_n,t_n-t\bigr\rangle \leq 0\).  Since \((At_n)\) converges weakly in \(H_0+C^0([-1,1])\) towards \(At\), this entails that:
\[
\limsup_{n\rightarrow+\infty} \bigl\langle A t_n,t_n\bigr\rangle \leq \bigl\langle A t,t\bigr\rangle.
\]
But, by lemma~\ref{thm:monotLemma}:
\[
\bigl\langle A t_n,t_n\bigr\rangle =  -\frac{1}{\pi} \Bigl\langle \log|x|*t_n,t_n\Bigr\rangle + \frac{1}{2}\int_{-1}^1 f' \biggl(\int_{-1}^x t_n\biggr)^2.
\]
By dominated convergence, we have:
\[
\lim_{n \rightarrow +\infty} \int_{-1}^1 f' \biggl(\int_{-1}^x t_n\biggr)^2 = \int_{-1}^1 f' \biggl(\int_{-1}^x t\biggr)^2,
\]
so that:
\[
\limsup_{n\rightarrow+\infty}  -\frac{1}{\pi} \Bigl\langle \log|x|*t_n,t_n\Bigr\rangle \leq  -\frac{1}{\pi} \Bigl\langle \log|x|*t,t\Bigr\rangle.
\]
Hence, \((t_n)\) is a weakly convergent sequence in the Hilbert space \(H^{-1/2}(-1,1)\) such that \(\limsup_{n\rightarrow+\infty} \|t_n\|_{H^{-1/2}} \leq  \|t\|_{H^{-1/2}}\).  By the weak lower semicontinuity of the norm in a Hilbert space, this classically entails that \(\lim_{n\rightarrow+\infty} \|t_n\|_{H^{-1/2}} =\|t\|_{H^{-1/2}}\).  Finally, we have proved that:
\[
\forall \tau\in H_{0}^*\cap\mathscr{M}([-1,1]),\qquad \lim_{n \rightarrow + \infty} \bigl\langle A t_n,t_n-\tau\bigr\rangle = \bigl\langle A t,t-\tau\bigr\rangle,
\]
which contains the expected conclusion. \qed

\subsection{Proof of the uniqueness, regularity, and homogenization results for heterogeneous friction and convex indentor}

\label{sec:proofConvexHetero}

The first part of this section is devoted to the proof of Theorem~\ref{thm:heartConvex}, that is, essentially to prove that the linear mapping \(\bar{A}:H_0^* \rightarrow H_0\) defined by the identity:
\[
\frac{e^{-\bar{\tau}}}{\scriptstyle\sqrt{1+f^2}}\biggl(- \frac{1}{\pi}\,\text{\rm pv}\frac{1}{x} * \tilde{t} + f\,\tilde{t}\biggr) = v',\qquad \text{in }\left]-1,1\right[,
\]
is continuous and coercive on \(H_0^*:=\{t\in H^{-1/2}(-1,1)|\langle t,1\rangle = 0\}\).  As we shall see, the fundamental reason is, once again, the link made by Theorem~\ref{thm:Riesz2} of Appendix~A between the Hilbert transform and a class of holomorphic functions on the complex upper half-plane.

\begin{lemm}
	\label{thm:propPhi}
	Let \(\bar{f}:\mathbb{R}\rightarrow\mathbb{R}\) be an arbitrary Lipschitz-continuous function with compact support.  We denote by \(\bar{\tau}:=-\frac{1}{\pi}\text{\rm pv}\frac{1}{x}*\arctan \bar{f}\) the Hilbert transform of \(\arctan \bar{f}\).  Let \(\phi(z)\) and \(\psi(z)\) be the holomorphic complex functions defined on the complex upper half-plane \(\Pi^+\) by:
	\[
	\phi(z) := \frac{1}{\pi}\int_{-\infty}^{\infty} \frac{\arctan\bar{f}(s)\,{\rm d}s}{s-z}, \qquad \psi(z):=e^{\phi(z)}-1.
	\]
	Then, the functions \(\phi(x+iy)\) and \(\psi(x+iy)\) converge as \(y \rightarrow 0+\), for almost all \(x\in\mathbb{R}\), towards the limits:
	\begin{align*}
		\phi(x+i0) & = \bar{\tau}(x) + i \,\arctan\bar{f}(x),\\
		\psi(x+i0) & = \biggl[\frac{e^{\bar{\tau}(x)}}{\scriptstyle\sqrt{1+\bar{f}^2(x)}}-1\biggr] + i\,\frac{\bar{f}(x)\,e^{\bar{\tau}(x)}}{\scriptstyle\sqrt{1+\bar{f}^2(x)}}.
	\end{align*}
	In addition, \(\phi(x+i0),\psi(x+i0)\in L^p(\mathbb{R};\mathbb{C})\), for all \(p\in \left]1,+\infty\right[\),  and we have:
	\[
	-\frac{1}{\pi}\,\text{\rm pv}\frac{1}{x}*\biggl[\frac{e^{\bar{\tau}}}{\scriptstyle\sqrt{1+\bar{f}^2}}-1\biggr] + \bar{f}\, \frac{e^{\bar{\tau}}}{\scriptstyle\sqrt{1+\bar{f}^2}} = 0,\qquad \text{on }\mathbb{R}.
	\]
	Finally, the function \(\Re e^{\phi(x+iy)}>0\) is positive on \(\Pi^+\), is bounded away from zero and infinity:
	\[
	\exists C_1,C_2\in \mathbb{R},\quad \forall (x,y)\in\overline{\Pi^+},\qquad 0<C_1\leq \Re e^{\phi(x+iy)}\leq C_2,
	\]
	and:
	\[
	\forall (x,y)\in\Pi^+,\qquad \inf \bar{f} \leq \frac{\Im e^{\phi(x+iy)}}{\Re e^{\phi(x+iy)}} \leq \sup \bar{f}.
	\]
\end{lemm}

\noindent\textbf{Proof.}
As \(\arctan\bar{f}\in L^p(\mathbb{R})\), for all \(p\in\left]1,+\infty\right[\), Theorem~\ref{thm:Riesz2} yields that existence of the almost everywhere pointwise limit \(\phi(x+i0)\in L^p(\mathbb{R};\mathbb{C})\), for all \(p\in\left]1,+\infty\right[\).  Using:
\[
\cos\arctan\bar{f} = \frac{1}{\sqrt{1+\bar{f}^2}},\qquad \sin\arctan\bar{f} = \frac{\bar{f}}{\sqrt{1+\bar{f}^2}},
\]
we get the existence of the almost everywhere pointwise limit \(\psi(x+i0)\) and its value.  As \(\bar{f}\) is Lipschitz-continuous with compact support, its Hilbert transform \(\bar{\tau}\) is continuous on \(\mathbb{R}\).  Since \(\bar{f}\) is compactly supported, \(\bar{\tau}(x)=O(1/x)\) and \(e^{\bar{\tau}(x)}-1=O(1/x)\), as \(|x|\rightarrow+\infty\), so that \(\bar{\tau}(x)\) and \(e^{\bar{\tau}(x)}-1\) have integrable \(p\)-power in the neighbourhood of infinity, for all \(p\in\left]1,+\infty\right[\).  Finally, we have proved that \(\phi(x+i0),\psi(x+i0)\in L^p(\mathbb{R};\mathbb{C})\), for all \(p\in\left]1,+\infty\right[\).  As \(\phi\) and \(\psi\) converge to \(0\) at infinity, Corollary~\ref{thm:noteTricomi} shows that \(\Re\psi(x+i0)\) is the Hilbert transform of \(\Im\psi(x+i0)\), which entails that \(-\Im\psi(x+i0)\) is the Hilbert transform of \(\Re\psi(x+i0)\).

Finally, we have:
\[
\phi(x+iy) = \frac{1}{\pi}\int_{-\infty}^\infty \frac{(s-x)\,\arctan\bar{f}(s)\,{\rm d}s}{(s-x)^2+y^2} + \frac{i}{\pi}\int_{-\infty}^\infty \frac{y\,\arctan\bar{f}(s)\,{\rm d}s}{(s-x)^2+y^2},
\]
so that, for all \((x,y)\in\Pi^+\):
\[
-\frac{\pi}{2}< \arctan (\inf\bar{f}) \leq \Im \phi(x+iy) \leq \arctan (\sup\bar{f}) < \frac{\pi}{2},
\]
which entails that:
\[
\forall (x,y)\in \Pi^+,\qquad \Re e^{\phi(x+iy)} > 0.
\]
Hence:
\[
\forall (x,y)\in \Pi^+,\qquad \frac{\Im e^{\phi(x+iy)}}{\Re e^{\phi(x+iy)}} = \tan \biggl(\frac{1}{\pi}\int_{-\infty}^\infty \frac{y\,\arctan\bar{f}(s)\,{\rm d}s}{(s-x)^2+y^2}\biggr)  \in \bigl[\inf\bar{f},\sup\bar{f}\bigr].
\]
As \(\bar{f}\) is Lipschitz-continuous and compactly supported, \(\bar{\tau}\) is in \(W^{1,p}(\mathbb{R})\), for all \(p\in \left]1,\infty\right[\), by Theorem~\ref{thm:Riesz1} of Appendix~A, which entails that \(\phi(x+iy)\) and \(e^{\phi(x+iy)}\) are both continuous functions on \(\overline{\Pi^+}\).  The additional fact that \(\phi(x+iy)\) goes to zero at infinity on \(\overline{\Pi^+}\) yields the fact that \(\Re e^{\phi(x+iy)}\) is bounded away from zero and infinity.
\qed

\bigskip
\begin{lemm}
	\label{thm:multHhalf}
	Let \(p\in \left]2,+\infty\right]\) and \(\theta\in L^\infty(\mathbb{R})\) such that \(\theta'\in L^p(\mathbb{R})\).  Then, the linear mappings:
	\[
	\left\{
	\begin{array}{rcl}
		H^{1/2}(\mathbb{R}) & \rightarrow & H^{1/2}(\mathbb{R})\\
		u & \mapsto & \theta u
	\end{array}
	\right.\qquad\qquad
	\left\{
	\begin{array}{rcl}
		H^{-1/2}(\mathbb{R}) & \rightarrow & H^{-1/2}(\mathbb{R})\\
		t & \mapsto & \theta t
	\end{array}
	\right.
	\]
	are continuous.
\end{lemm}

\noindent\textbf{Proof.} The space \(H^{1/2}(\mathbb{R})\) is endowed with the Sobolev-Slobodetskii norm:
\[
\bigl\|u\bigr\|_{H^{1/2}}^2 := \bigl\|u\bigr\|_{L^{2}}^2 + \int_{\mathbb{R}\times\mathbb{R}} \biggl(\frac{u(x)-u(y)}{x-y}\biggr)^2{\rm d}x\,{\rm d}y.
\]
For \(u\in H^{1/2}(\mathbb{R})\), we have:
\[
\bigl\| \theta u\bigr\|_{L^2}^2 \leq \bigl\| \theta \bigr\|_{L^\infty}^2\bigl\|  u\bigr\|_{L^2}^2,
\]
and:
\begin{multline*}
	\int_{\mathbb{R}\times\mathbb{R}} \biggl(\frac{\theta(x)u(x)-\theta(y)u(y)}{x-y}\biggr)^2{\rm d}x\,{\rm d}y \leq 2\int_{\mathbb{R}\times\mathbb{R}} \theta^2(x)\,\biggl(\frac{u(x)-u(y)}{x-y}\biggr)^2{\rm d}x\,{\rm d}y +\mbox{}\\
	\mbox{}+ 2 \int_{\mathbb{R}\times\mathbb{R}} u^2(y)\,\biggl(\frac{\theta(x)-\theta(y)}{x-y}\biggr)^2{\rm d}x\,{\rm d}y.
\end{multline*}
Fixing \(y\in\mathbb{R}\), we have:
\begin{align*}
	\int_{-\infty}^\infty \biggl(\frac{\theta(x)-\theta(y)}{x-y}\biggr)^2{\rm d}x	
	& \leq \int_{\mathbb{R}\setminus[y-1,y+1]} \biggl(\frac{\theta(x)-\theta(y)}{x-y}\biggr)^2{\rm d}x
	+ \int_{y-1}^{y+1}\biggl(\frac{\theta(x)-\theta(y)}{x-y}\biggr)^2{\rm d}x	
	,\\
	& \leq  8\bigl\|\theta\bigr\|_{L^\infty}^2 + \int_{y-1}^{y+1}\biggl(\frac{1}{|x-y|}\int_{x}^y|\theta'(s)|\,{\rm d}s\biggr)^2{\rm d}x,\\
	& \leq  8\bigl\|\theta\bigr\|_{L^\infty}^2 + \int_{y-1}^{y+1}\frac{1}{|x-y|^{2/p}}\biggl(
	\int_{x}^y|\theta'(s)|^p\,{\rm d}s\biggr)^{2/p}{\rm d}x,\\
	& \leq  8\bigl\|\theta\bigr\|_{L^\infty}^2 + \frac{2p}{p-2}\bigl\|\theta'\bigr\|_{L^p}^2.
\end{align*}
Finally, we have proved:
\[
\bigl\|\theta u\bigr\|_{H^{1/2}}^2 \leq \Bigl[19\bigl\|\theta\bigr\|_{L^\infty}^2 + \frac{4p}{p-2}\bigl\|\theta'\bigr\|_{L^p}^2\Bigr]\bigl\|u\bigr\|_{H^{1/2}}^2,
\]
which yields the claimed result for the first mapping.  For the second mapping, we recall that \(H^{-1/2}(\mathbb{R})\) is the dual space of \(H^{1/2}(\mathbb{R})\), and that we have:
\[
\bigl\langle \theta t ,u\bigr\rangle = \bigl\langle t ,\theta u\bigr\rangle,
\]
by definition.\qed

\bigskip
The following corollary will be the cornerstone of the proof of Theorem~\ref{thm:heartConvex}, but also of the homogenization analysis.

\begin{prop}
	\label{thm:fundCoroConvex}
	Let \(\bar{f}\in W^{1,\infty}(\mathbb{R})\) be a compactly supported Lipschitz-continuous function, and \(\eta\in H^{1/2}(\mathbb{R})\).  Then, 
	\[
	\frac{e^{-\bar{\tau}}}{\scriptstyle\sqrt{1+\bar{f}^2}}\biggl(-\frac{1}{\pi}\,\text{\rm pv}\frac{1}{x} * \eta' + \bar{f}\,\eta'\biggr) \in H^{-1/2}(\mathbb{R}),
	\]
	and we have the Green's formula, for all \(\varphi\in H^1(\Pi^+)\):
	\begin{align}
		\biggl\langle 	\frac{e^{-\bar{\tau}}}{\scriptstyle\sqrt{1+\bar{f}^2}}\biggl(\frac{1}{\pi}\,\text{\rm pv}\frac{1}{x} * \eta' - \bar{f}\,\eta'\biggr) \;,\;\varphi(x,0)\biggr\rangle & = \int_{\Pi^+} \Bigl[\Re e^{-\phi(x+iy)}\frac{\partial \Phi_\eta}{\partial x} + \Im  e^{-\phi(x+iy)}\frac{\partial \Phi_\eta}{\partial y}\Bigr]\frac{\partial \varphi}{\partial x} +\mbox{}\nonumber\\
		& \qquad\mbox{} + 
		\Bigl[-\Im e^{-\phi(x+iy)}\frac{\partial \Phi_\eta}{\partial x} + \Re  e^{-\phi(x+iy)}\frac{\partial \Phi_\eta}{\partial y}\Bigr]\frac{\partial \varphi}{\partial y},
		\label{eq:genGreenFormula}
	\end{align}
	where:
	\[
	\phi(z) :=  \frac{1}{\pi}\int_{-\infty}^{\infty} \frac{\arctan\bar{f}(s)\,{\rm d}s}{s-z},\qquad \text{for }z\in \Pi^+,
	\]
	and \(\Phi_{\eta}\) denotes the Poisson integral of \(\eta\) (see Proposition~\ref{thm:prolAnal} in Appendix~C):
	\[
	\Phi_\eta(x,y):= \frac{1}{\pi}\int_{-\infty}^\infty \frac{y\,\eta(s)\,{\rm d}s}{(s-x)^2+y^2} = \frac{1}{\pi}\frac{y}{x^2+y^2}\stackrel{x}{*}\eta = \frac{1}{\pi}\arctan \frac{x}{y}\stackrel{x}{*}\eta',
	\]
	and where the terms between \([\,]\) in formula~\eqref{eq:genGreenFormula} are the two components of the gradient of a harmonic function in \(\Pi^+\).  This entails the following formula valid for arbitrary \(\eta_1,\eta_2\in H^{1/2}(\mathbb{R})\):
	\begin{align*}
		\biggl\langle 	\frac{e^{-\bar{\tau}}}{\scriptstyle\sqrt{1+\bar{f}^2}}\biggl(\frac{1}{\pi}\,\text{\rm pv}\frac{1}{x} * \eta_1' - \bar{f}\,\eta_1'\biggr) \;,\;\eta_2\biggr\rangle & = \int_{\Pi^+} \Bigl[\Re e^{-\phi(x+iy)}\frac{\partial \Phi_{\eta_1}}{\partial x} + \Im  e^{-\phi(x+iy)}\frac{\partial \Phi_{\eta_1}}{\partial y}\Bigr]\frac{\partial \Phi_{\eta_2}}{\partial x} +\mbox{}\\
		& \qquad\mbox{} + 
		\Bigl[-\Im e^{-\phi(x+iy)}\frac{\partial \Phi_{\eta_1}}{\partial x} + \Re  e^{-\phi(x+iy)}\frac{\partial \Phi_{\eta_1}}{\partial y}\Bigr]\frac{\partial \Phi_{\eta_2}}{\partial y}.
	\end{align*}
\end{prop}

\noindent\textbf{Proof.}
Since \(\bar{f}\in W^{1,\infty}(\mathbb{R})\) is a compactly supported Lipschitz-continuous function, \(\bar{\tau}\in W^{1,p}(\mathbb{R})\), for all \(p\in \left]1,\infty\right[\).  This entails that the function:
\[
\theta := \frac{e^{-\bar{\tau}}}{\scriptstyle\sqrt{1+\bar{f}^2}},
\]
satisfies \(\theta \in L^\infty(\mathbb{R})\) and \(\theta'\in L^p(\mathbb{R})\), for all \(p\in \left]1,\infty\right[\).  Also, the definition of the spaces \(H^{-1/2}(\mathbb{R})\) and \(H^{1/2}(\mathbb{R})\) in terms of the Fourier transform (Proposition~\ref{thm:defSobolevFourier} in Appendix~B) shows that \(\eta\in H^{1/2}(\mathbb{R}) \Rightarrow \eta'\in H^{-1/2}(\mathbb{R})\).  Hence, 
\[
\frac{e^{-\bar{\tau}}}{\scriptstyle\sqrt{1+\bar{f}^2}}\biggl(-\frac{1}{\pi}\,\text{\rm pv}\frac{1}{x} * \eta' + \bar{f}\,\eta'\biggr) \in H^{-1/2}(\mathbb{R}), 
\]
is now an immediate consequence of Lemma~\ref{thm:multHhalf}. 

Pick \(\eta\in C_c^\infty(\mathbb{R})\) and consider the following holomorphic complex function 
\[
\tilde{\psi}(x+iy) := ie^{-\phi(x+iy)}\,\frac{1}{\pi}\int_{-\infty}^\infty \frac{\eta'(s)\,{\rm d}s}{s-(x+iy)} = ie^{-\phi(x+iy)}\,\biggl(\frac{\partial \Phi_\eta}{\partial y}(x,y) + i\frac{\partial \Phi_\eta}{\partial x}(x,y)\biggr),
\]
where the last identity is established using the last definition of \(\Phi_\eta\) in formula~\eqref{eq:defPoissonIntegral}.  By Lemma~\ref{thm:propPhi}, the function \(e^{-\phi(x+iy)}\) is bounded on \(\overline{\Pi^+}\), so that Proposition~\ref{thm:prolAnal} in Appendix~C entails that \(\tilde{\psi}\in L^2(\Pi^+;\mathbb{C})\).  In addition, using Theorem~\ref{thm:Riesz2} in Appendix~A, \(\tilde{\psi}(x+iy)\) converges as \(y \rightarrow 0+\), for almost all \(x\in\mathbb{R}\), but also in \(L^p(\mathbb{R})\), for all \(p\in\left]1,+\infty\right[\), towards:
\[
\tilde{\psi}(x+i0) = \frac{e^{-\bar{\tau}}}{\scriptstyle\sqrt{1+\bar{f}^2}}\bigl(\bar{f}+i\bigr)\Bigl(-\frac{1}{\pi}\text{pv}\frac{1}{x}*\eta' + i \eta'\Bigr).
\]
Finally, the harmonic function \(\tilde{h}(x,y):=(1/\pi)\log\sqrt{x^2+y^2}\stackrel{x}{*}\Im\tilde{\psi}(x+i0)\) is linked to \(\tilde{\psi}\) by:
\[
\forall (x,y)\in\Pi^+,\qquad 	\tilde{\psi}(x+iy) = -\frac{\partial \tilde{h}}{\partial x}(x,y) + i \frac{\partial \tilde{h}}{\partial y}(x,y),
\]
and satisfies the Green's formula:
\[
\forall\varphi\in C_c^\infty(\overline{\Pi^+}),\qquad -\int_{-\infty}^\infty \Im\tilde{\psi}(x+i0)\,\varphi(x,0)\,{\rm d}x = \int_{\Pi^+} \biggl(\frac{\partial \tilde{h}}{\partial x}\frac{\partial \varphi}{\partial x} + \frac{\partial \tilde{h}}{\partial y}\frac{\partial \varphi}{\partial y}\biggr){\rm d}x\,{\rm d}y.
\]
This is readily seen to be the expected Green's formula in the particular case \(\eta\in C_c^\infty(\mathbb{R})\) and \(\varphi\in C_c^\infty(\overline{\Pi^+})\).  By Proposition~\ref{thm:prolAnal} in Appendix~C and the definition of \(H^{1/2}(\mathbb{R})\) in terms of the Fourier transform (Proposition~\ref{thm:defSobolevFourier} in Appendix~B), we have:
\[
\forall \eta\in H^{1/2}(\mathbb{R}),\qquad \bigl\| \nabla \Phi_\eta\bigr\|_{L^2(\Pi^+)}^2 = \frac{1}{\pi}\bigl\langle\text{\rm pv}1/x*\eta',\eta\bigr\rangle \leq \bigl\| \eta\bigr\|_{H^{1/2}(\mathbb{R})}^2.
\]
By the density of \(C_c^\infty(\mathbb{R})\) in \(H^{1/2}(\mathbb{R})\), Lemma~\ref{thm:multHhalf} (to go to the limit on the boundary term) and Lemma~\ref{thm:propPhi} (which yields that \(e^{-\phi(x+iy)}\) is bounded on \(\Pi^+\) and enables to go to the limit in the integral over \(\Pi^+\)), we obtain the Green's formula for all \(\eta\in H^{1/2}(\mathbb{R})\) and \(\varphi\in C_c^\infty(\overline{\Pi^+})\).  The generalization to \(\varphi\in H^1(\Pi^+)\) now follows from the density of \(C_c^\infty(\overline{\Pi^+})\) in \(H^1(\Pi^+)\).

Finally, take arbitrary \(\eta_1,\eta_2\in  H^{1/2}(\mathbb{R})\).  As \(\Phi_{\eta_2}\notin H^1(\Pi^+)\), in general, but is only in \(H^1(\mathbb{R}\times \left]0,Y\right[)\), for all \(Y>0\), by Proposition~\ref{thm:prolAnal}, we define:
\[
\varphi_n(x,y) := \min\{1,n/y\}\,\Phi_{\eta_2}(x,y),
\]
so that \(\varphi_n\in H^1(\Pi^+)\) and \((\nabla \varphi_n)\) converges strongly in \(L^2(\Pi^+)\) towards \(\nabla \Phi_{\eta_2}\).  Using \(\varphi_n\) as a test-function in the Green's formula written for \(\Phi_{\eta_1}\) and going to the limit yields the final claim of the proposition.\qed

\bigskip
\noindent\textbf{Proof of Theorem~\ref{thm:heartConvex}.}
We assume that \(f\in W^{1,\infty}(-1,1)\), and we consider any compactly supported extension \(\bar{f}\in W^{1,\infty}(\mathbb{R})\) of \(f\).  

\noindent\textbf{Step 1.} \textit{Let \(\tilde{t}\in H_0^*:=\{t\in H^{-1/2}(-1,1)\;|\;\langle t,1\rangle=0\}\).  Then there exists a unique \(v \in H^{1/2}(\mathbb{R})\) such that:
\[
\frac{e^{-\bar{\tau}}}{\scriptstyle\sqrt{1+\bar{f}^2}}\biggl(- \frac{1}{\pi}\,\text{\rm pv}\frac{1}{x} * \tilde{t} + \bar{f}\,\tilde{t}\biggr) = v',\qquad \text{in } \mathscr{D}'(\mathbb{R}).
\]
}
We know that \(\tilde{t}*\log|x|\in H^{1/2}(\mathbb{R})\) by using Propositions~\ref{thm:FourierTransforms} and~\ref{thm:defSobolevFourier}, so that we only have to prove that there exists a unique \(\tilde{v}\in H^{1/2}(\mathbb{R})\) such that:
\begin{equation}
	\label{eq:distOverLine}
	\biggl(\frac{e^{-\bar{\tau}}}{\scriptstyle\sqrt{1+\bar{f}^2}}-1\biggr)\biggl(- \frac{1}{\pi}\,\text{\rm pv}\frac{1}{x} * \tilde{t}\biggr) + \bar{f}\,\frac{e^{-\bar{\tau}}}{\scriptstyle\sqrt{1+\bar{f}^2}}\,\tilde{t}  = \tilde{v}',\qquad \text{in } \mathscr{D}'(\mathbb{R}),
\end{equation}
because \(v:=\tilde{v}-(1/\pi)\tilde{t}*\log|x|\) would then satisfy the claim.  The function \(e^{-\bar{\tau}}/{\scriptstyle\sqrt{1+\bar{f}^2}}-1\) is readily checked to be in \(H^1(\mathbb{R})\subset H^{1/2}(\mathbb{R})\).  For \((u,t)\in H^{1/2}(\mathbb{R})\times H^{-1/2}(\mathbb{R})\), we have \(\langle u, \text{pv}1/x * t\rangle = - \langle  \text{pv}1/x * u, t\rangle\).  Hence, Lemma~\ref{thm:propPhi} for \(-\bar{f}\) yields:
\begin{equation}
	\label{eq:constInfinity}
	\biggl\langle \frac{e^{-\bar{\tau}}}{\scriptstyle\sqrt{1+\bar{f}^2}}-1\;,\; - \frac{1}{\pi}\,\text{\rm pv}\frac{1}{x} * \tilde{t}\biggr\rangle + \biggl\langle \bar{f}\,\frac{e^{-\bar{\tau}}}{\scriptstyle\sqrt{1+\bar{f}^2}}\;,\;\tilde{t}\biggr\rangle = 0.
\end{equation}
As \(\tilde{t}\in H^{-1/2}(-1,1)\) is a compactly supported distribution, the left-hand side of~\eqref{eq:distOverLine} (which is in \(H^{-1/2}(\mathbb{R})\), thanks to Lemma~\ref{thm:multHhalf}) can be written as a sum \(T+L\) where \(T\in H^{-1/2}(\mathbb{R})\) with compact support and \(L:\mathbb{R}\rightarrow\mathbb{R}\) is a continuous integrable function such that \(L(x)=O(1/x^2)\) at infinity.  Denoting by \(H\) the Heaviside function, \(H*T\) is therefore in \(L_{\rm loc}^2\).  Setting \(\tilde{v}:= H*(T+L)\), the function \(\tilde{v}\) is therefore in \(L_{\rm loc}^2\) and fulfils identity~\eqref{eq:distOverLine}.  In addition, we have \(\tilde{v}(x) = O(1/x)\) at \(-\infty\) and \(\tilde{v}(x) = C_{\infty} + O(1/x)\) at \(+\infty\), where \(C_{\infty}\) is checked to equal the left-hand side of identity~\eqref{eq:constInfinity} and therefore vanishes.  Hence, \(\tilde{v}\in L^2(\mathbb{R})\).  By Proposition~\ref{thm:defSobolevFourier} in Appendix~B, we see that a function \(\tilde{v}\in L^2(\mathbb{R})\) such that \(\tilde{v}'\in H^{-1/2}(\mathbb{R})\) is actually in \(H^{1/2}(\mathbb{R})\).  Hence, \(\tilde{v}\in H^{1/2}(\mathbb{R})\), and it is clearly the only solution of~\eqref{eq:distOverLine} with this regularity.  The proof of Step~1 is complete.

\smallskip
\noindent\textbf{Step 2.} \textit{The mapping \(\bar{A}:H_0^* \rightarrow H_0\), which associates with \(\tilde{t}\in H_0^*:=\{t\in H^{-1/2}(-1,1)\;|\;\langle t,1\rangle=0\}\) the unique \(\mathring{v}=\bar{A}\tilde{t}\) in \(H_0:=H^{1/2}(-1,1)/\mathbb{R}\) such that:
\[
\frac{e^{-\bar{\tau}}}{\scriptstyle\sqrt{1+f^2}}\biggl(- \frac{1}{\pi}\,\text{\rm pv}\frac{1}{x} * \tilde{t} + \bar{f}\,\tilde{t}\biggr) = \mathring{v}',\qquad \text{in }\mathscr{D}'(\left]-1,1\right[),
\]
is continuous and coercive.}

The fact that there exists a unique \(v \in H_0\) satisfying the above identity is a direct consequence of step~1.

Consider arbitrary \(\tilde{t}_1,\tilde{t}_2\in H_0^*\).  Their extension by zero to the real line are elements of \(H^{-1/2}(\mathbb{R})\) whose Fourier transforms are \(C^\infty\) functions vanishing at \(0\).  By the definition of the spaces \(H^{1/2}(\mathbb{R})\) and \(H^{1/2}(\mathbb{R})\) in terms of the Fourier transform (Proposition~\ref{thm:defSobolevFourier} in Appendix~B), this entails that the distributions \(\eta_i:=(1/2)\,\text{\rm sgn}*\tilde{t}_i\) are in \(H^{1/2}(\mathbb{R})\).  By Proposition~\ref{thm:fundCoroConvex},
we obtain:
\[
\Biggl|  \biggl\langle 	\frac{e^{-\bar{\tau}}}{\scriptstyle\sqrt{1+\bar{f}^2}}\biggl(\frac{1}{\pi}\,\text{\rm pv}\frac{1}{x} * \tilde{t}_1 - \bar{f}\,\tilde{t}_1\biggr) \;,\;\frac{1}{2}\,\text{\rm sgn}*\tilde{t}_2\biggr\rangle \Biggr|  \leq C \bigl\| \nabla \Phi_{\eta_1} \bigr\|_{L^2(\Pi^+)}\, \bigl\| \nabla \Phi_{\eta_2} \bigr\|_{L^2(\Pi^+)},
\]
for some real constant \(C\) independent of \(\tilde{t}_i\), since the function \(e^{-\bar{\tau}(x+iy)}\) is bounded on \(\Pi^+\) by Lemma~\ref{thm:propPhi}.  Combining Proposition~\ref{thm:prolAnal} and~\ref{thm:lemmHminusHalf}, we have:
\[
\forall \tilde{t}\in H_0^*,\qquad \bigl\| \nabla \Phi_\eta\bigr\|_{L^2(\Pi^+)}^2 = -\frac{1}{\pi}\bigl\langle \log|x|*\tilde{t},\tilde{t}\bigr\rangle = \bigl\| \tilde{t}\bigr\|_{H^{-1/2}(-1,1)}^2,
\]
where \(\eta:=(1/2)\text{\rm sgn}*\tilde{t}\in H^{1/2}(\mathbb{R})\), so that, we have proved that the bilinear form:
\[
\tilde{t}_1,\tilde{t}_2 \mapsto \biggl\langle 	\frac{e^{-\bar{\tau}}}{\scriptstyle\sqrt{1+\bar{f}^2}}\biggl(\frac{1}{\pi}\,\text{\rm pv}\frac{1}{x} * \tilde{t}_1 - \bar{f}\,\tilde{t}_1\biggr) \;,\;\frac{1}{2}\,\text{\rm sgn}*\tilde{t}_2\biggr\rangle
\]
is continuous on \(H_0^*\).  Hence, the formula:
\[
\bigl\langle \bar{A}\tilde{t}_1,\tilde{t}_2\bigr\rangle = \biggl\langle 	\frac{e^{-\bar{\tau}}}{\scriptstyle\sqrt{1+\bar{f}^2}}\biggl(\frac{1}{\pi}\,\text{\rm pv}\frac{1}{x} * \tilde{t}_1 - \bar{f}\,\tilde{t}_1\biggr) \;,\;\frac{1}{2}\,\text{\rm sgn}*\tilde{t}_2\biggr\rangle
\]
defines a continuous linear mapping \(\bar{A}:H_0^*\rightarrow H_0\).  Also, using again Lemma~\ref{thm:propPhi} and Proposition~\ref{thm:fundCoroConvex}, we obtain:
\[
\exists C>0 , \quad \forall\tilde{t}\in H_0^*,\qquad \bigl\langle \bar{A}\tilde{t},\tilde{t}\bigr\rangle \geq C\, \bigl\| \nabla \Phi_{\eta} \bigr\|_{L^2(\Pi^+)}^2 = C \bigl\| \tilde{t}\bigr\|_{H^{-1/2}(-1,1)}^2,
\]
which shows that \(\bar{A}:H_0^*\rightarrow H_0\) is coercive, so that it is an isomorphism, by the Lax-Milgram theorem.

\smallskip
\noindent\textbf{Step 3.} \textit{Let \(\tilde{t}\in H_0^*\) and the element \(v\in H^{1/2}(\mathbb{R})\) satisfying:
\[
v' = \frac{e^{-\bar{\tau}}}{\scriptstyle\sqrt{1+\bar{f}^2}}\biggl(- \frac{1}{\pi}\,\text{\rm pv}\frac{1}{x} * \tilde{t} + \bar{f}\,\tilde{t}\biggr) ,\qquad \text{in } \mathscr{D}'(\mathbb{R}).
\]
Then, we have the inversion formula:
\[
\tilde{t} = \frac{e^{\bar{\tau}}}{\scriptstyle\sqrt{1+\bar{f}^2}}\biggl( \frac{1}{\pi}\,\text{\rm pv}\frac{1}{x} * v' + \bar{f}\,v'\biggr),\qquad \text{in } \mathscr{D}'(\mathbb{R}).
\]
}

Pick an arbitrary \(\tilde{t}\in H_0^*\).  Assume in addition that \(t\in L^p(-1,1)\) for some \(p\in \left]1,\infty\right[\).  Set \(\tilde{v} := v - (1/\pi) \log|x| * \tilde{t}\), so that \(\tilde{v}\in H^{1/2}(\mathbb{R})\), \(\tilde{v}'\in L^p(\mathbb{R})\) and:
\[
\tilde{v}' = \biggl(\frac{e^{-\bar{\tau}}}{\scriptstyle\sqrt{1+\bar{f}^2}}-1\biggr)\biggl(- \frac{1}{\pi}\,\text{\rm pv}\frac{1}{x} * \tilde{t}\biggr) + \bar{f}\,\frac{e^{-\bar{\tau}}}{\scriptstyle\sqrt{1+\bar{f}^2}}\,\tilde{t} ,\qquad \text{in } \mathscr{D}'(\mathbb{R}).
\]
Recalling that the function:
\[
\frac{e^{-\bar{\tau}}}{\scriptstyle\sqrt{1+\bar{f}^2}} - 1,
\]
is in \(L^q(\mathbb{R})\), for all \(q\in\left]1,\infty\right[\), and using Lemma~\ref{thm:propPhi} and Corollary~\ref{thm:pbt} of Appendix~A, we obtain:
\[
-\frac{1}{\pi} \text{pv}\frac{1}{x}*\tilde{v}'=\bar{f}\,\frac{e^{-\bar{\tau}}}{\scriptstyle\sqrt{1+\bar{f}^2}}\biggl(- \frac{1}{\pi}\,\text{\rm pv}\frac{1}{x} * \tilde{t}\biggr) - \biggl(\frac{e^{-\bar{\tau}}}{\scriptstyle\sqrt{1+\bar{f}^2}}-1\biggr)\,\tilde{t}.
\]
Combining:
\begin{align*}
	v' & = \frac{e^{-\bar{\tau}}}{\scriptstyle\sqrt{1+\bar{f}^2}}\biggl(- \frac{1}{\pi}\,\text{\rm pv}\frac{1}{x} * \tilde{t} + \bar{f}\,\tilde{t}\biggr),\\
	-\frac{1}{\pi} \text{pv}\frac{1}{x}*v' & = \frac{e^{-\bar{\tau}}}{\scriptstyle\sqrt{1+\bar{f}^2}}\biggl(- \frac{\bar{f}}{\pi}\,\text{\rm pv}\frac{1}{x} * \tilde{t} - \tilde{t}\biggr),
\end{align*}
we obtain the expected inversion formula for the particular case \(\tilde{t}\in L^p(-1,1)\;\Rightarrow\; v'\in L^p(\mathbb{R})\).  Using the density of \(H_0^*\cap L^p(-1,1)\) into \(H_0^*\) and making use of Lemma~\ref{thm:multHhalf}, we obtain the general case. 

\smallskip
\noindent\textbf{Step 4.} \textit{The operator \(\bar{A}^{-1}\) has the T-monotonicity property.}

The above inversion formula is all what is needed to prove the T-monotonicity property of the operator \(\bar{A}^{-1}\) along the lines of the proof already given in the case of homogeneous friction in Section~\ref{sec:homogeneousFriction}.  Take an arbitrary \(v\in H^{1/2}(-1,1)\).  Set \(\tilde{t}=\bar{A}^{-1}\mathring{v}\in H_0^*\) and \(\bar{v}\) as the unique function in \(H^{1/2}(\mathbb{R})\) such that:
\[
\bar{v}' = \frac{e^{-\bar{\tau}}}{\scriptstyle\sqrt{1+\bar{f}^2}}\biggl(- \frac{1}{\pi}\,\text{\rm pv}\frac{1}{x} * \tilde{t} + \bar{f}\,\tilde{t}\biggr),\qquad \text{in } \mathscr{D}'(\mathbb{R}).
\]
As \(v\) and \(\bar{v}\) have the same derivative on \(\left]-1,1\right[\), there exists \(C\in \mathbb{R}\) such that \(v=\bar{v}+C\) on \(\left]-1,1\right[\), so that:
\[
\bigl\langle \bar{A}^{-1}\mathring{v},v^+\bigr\rangle = \biggl\langle \frac{e^{\bar{\tau}}}{\scriptstyle\sqrt{1+\bar{f}^2}}\biggl(\frac{1}{\pi} \text{pv}\frac{1}{x}*\bar{v}' + \bar{f}\,\bar{v}'\biggr)\;, (\bar{v}+C)^+-C^+\biggr\rangle,
\]
where \((\bar{v}+C)^+-C^+\in H^{1/2}(\mathbb{R})\) (see Lemma~\ref{thm:coroR}).  Applying the Green's formula~\eqref{eq:genGreenFormula} of Proposition~\ref{thm:fundCoroConvex} upon replacing \(\bar{f}\) by \(-\bar{f}\), \(\eta\) by \(\bar{v}\), \(\varphi\) by \((\Phi_{\bar{v}}+C)^+-C^+\) and mimicking the proof of Lemma~\ref{thm:coroR}, we get:
\begin{align*}
	\lefteqn{\biggl\langle 	\frac{e^{\bar{\tau}}}{\scriptstyle\sqrt{1+\bar{f}^2}}\biggl(\frac{1}{\pi}\,\text{\rm pv}\frac{1}{x} * \bar{v}' + \bar{f}\,\bar{v}'\biggr) \;,\;(\bar{v}+C)^+-C^+\biggr\rangle} \\
	& \qquad\qquad \mbox{} = \int_{\Pi^+} \Bigl[\Re e^{\phi(x+iy)}\frac{\partial \Phi_{\bar{v}}}{\partial x} + \Im  e^{\phi(x+iy)}\frac{\partial \Phi_{\bar{v}}}{\partial y}\Bigr]\frac{\partial}{\partial x}\Bigl(\Phi_{\bar{v}}+C\Bigr)^+ +\mbox{}\\
	&\qquad\qquad\qquad\mbox{} + 
	\Bigl[-\Im e^{\phi(x+iy)}\frac{\partial \Phi_{\bar{v}}}{\partial x} + \Re  e^{\phi(x+iy)}\frac{\partial \Phi_{\bar{v}}}{\partial y}\Bigr]\frac{\partial}{\partial y}\Bigl(\Phi_{\bar{v}}+C\Bigr)^+, \\
	&
	\qquad\qquad \mbox{}= \int_{\Pi^+}\Re e^{\phi(x+iy)}\biggl\{ \biggl[\frac{\partial}{\partial x}\Bigl(\Phi_{\bar{v}}+C\Bigr)^+\biggr]^2+\biggl[\frac{\partial}{\partial y}\Bigl(\Phi_{\bar{v}}+C\Bigr)^+\biggr]^2 \biggr\}\geq 0,
\end{align*}
where Proposition~\ref{thm:prolAnal} and Stampacchia's theorem \cite[Theorem~1.56, p~79]{Troianiello} have been used.  If the equality is achieved, then \((\Phi_{\bar{v}}+C)^+\) must be a constant, implying that \((\bar{v}+C)^+\) is also constant. \qed

\bigskip
We now turn to the detailed proof of the homogenization result in the case of an arbitrary convex indentor (Theorem~\ref{thm:homogConvex}).  We are given \(P>0\) and a \emph{convex} indentor shape \(g\in W^{1,\infty}(-1,1)\).  The microscopic friction coefficient \(f_n\) is built from a given period \(p\in W^{1,\infty}(-1,1)\), such that \(p(-1)=p(1)\) which is extended by $2$-periodicity to the whole line.  The Lipschitz-continuous function \(f_n\) is defined on \(\left]-1,1\right[\) by:
\[
f_n(x) := p(nx).
\]
Picking an arbitrary function \(\varphi\in C_c^\infty(\mathbb{R};[0,1])\) such that \(\varphi(x)=1\), for all \(x\in [-1,1]\), the function \(\bar{f}_n(x):=\varphi(x)\,p(nx)\) is a compactly supported Lipschitz-continuous extension of \(f_n\) to the whole line.  Denoting by \(f_\textrm{\rm eff}\) the constant:
\[
f_\textrm{\rm eff} := \tan \bigl\langle \arctan f_1 \bigr\rangle,
\]
Lemma~\ref{thm:lemmHomog} yields that the sequence \((f_n)\) converges weakly-* in \(L^\infty(-1,1)\) towards the constant function \(f_\textrm{\rm eff}\) and the sequence \((\bar{f}_n)\) converges weakly-* in \(L^\infty(\mathbb{R})\) towards \(f_\textrm{\rm eff}\,\varphi(x)\).

Let  \(\tilde{t}_n\in \cup_{p>1}L^p(-1,1)\) be the unique solution of the contact problem (either problem~\(\mathscr{P}_{\!\rm o}\) or problem~\(\mathscr{P}_{\!\rm a}\))	provided by Corollary~\ref{thm:exUniqConvex}.  Let also \(\tilde{t}_\textrm{\rm eff}\) be the unique solution of the contact problem associated with the constant \(f_\textrm{\rm eff}\).  We also denote by \(t_n^0\) the function:
\[
t_n^0(x) :=  \frac{e^{\tau_n(x)}}{\pi\sqrt{1-x^2}\sqrt{1+f_n^2(x)}}, \qquad \text{on }\left]-1,1\right[.
\]
The function \(t_n^0\) has integral over \(\left]-1,1\right[\) equal to \(1\) and solves the homogeneous Carleman equation.  The function \(-Pt_n^0\) solves the contact problem with microscopic friction coefficient \(f_n\) in the case of the flat indentor \(g=0\).  From the analysis in Section~\ref{sec:constantFriction}, we know that the sequence \((t_n^0)\) converges weakly in \(L^p(-1,1)\), for all \(p\in \left]1,(1/2+\beta)^{-1}\right[\), where:
\begin{equation}
	\label{eq:beta}
	\beta := \frac{1}{\pi}\arctan \bigl\|f_1\bigr\|_{L^\infty(-1,1)} < \frac{1}{2},
\end{equation}
towards:
\[
t_{\rm eff}^0(x) :=  \frac{e^{\tau_{\rm eff}(x)}}{\pi\sqrt{1-x^2}\sqrt{1+f_{\rm eff}^2(x)}}, \qquad \text{on }\left]-1,1\right[.
\]
In particular, as the embedding of \(L^p(-1,1)\) (\(p\in \left]1,\infty\right[\)) into \(H^{-1/2}(-1,1)\) is compact, the sequence \((t_n^0)\) converges \emph{strongly} in \(H^{-1/2}(-1,1)\) towards \(t_{\rm eff}^0\).

The function \(\tilde{t}_n\in \cup_{p>1}L^p(-1,1)\) is characterized by the following assertions:
\begin{itemize}
	\item \(\displaystyle\int_{-1}^1\tilde{t}_n\,{\rm d}x=0\) and \(\tilde{t}_n\leq Pt_n^0\).
	\item There exists (uniquely) \(C_n\in\mathbb{R}\) such that:
	\begin{align}
		&\int_{-1}^x \frac{e^{-\bar{\tau}_n}}{\scriptstyle \sqrt{1+f_n^2}}\biggl(-\frac{1}{\pi}\text{pv}\frac{1}{x}*\tilde{t}_n+f_n\tilde{t}_n-g'\biggr){\rm d}s\leq C_n,\qquad \text{and}\qquad \nonumber\\
		& \biggl\langle \tilde{t}_n-Pt_n^0\;,\;\int_{-1}^x \frac{e^{-\bar{\tau}_n}}{\scriptstyle \sqrt{1+f_n^2}}\biggl(-\frac{1}{\pi}\text{pv}\frac{1}{x}*\tilde{t}_n+f_n\tilde{t}_n-g'\biggr){\rm d}s-C_n\biggr\rangle = 0.
		\label{eq:defCn}
	\end{align}
\end{itemize}

The following lemma contains the technical estimate on which the proof of Theorem~\ref{thm:homogConvex} is going to rely.

\begin{lemm}
	\label{thm:homogEstimate}
	The sequence:
	\[
	\Biggl(\frac{e^{-\bar{\tau}_n}}{\scriptstyle \sqrt{1+\bar{f}_n^2}}\biggl(-\frac{1}{\pi}\text{\rm pv}\frac{1}{x}*\tilde{t}_n+f_n\tilde{t}_n\biggr)\Biggr)_{n\in\mathbb{N}\setminus\{0\}}
	\]
	is bounded in \(L^p(\mathbb{R})\), for all \(p\in \left]1,(1/2 + \beta)^{-1}\right[\), where \(\beta\) is given by formula~\eqref{eq:beta}.
\end{lemm}

\noindent\textbf{Proof.}

\noindent\textbf{Step 1.} \textit{Setting:
\[
u_n' := -\frac{1}{\pi}\,\text{\rm pv}\frac{1}{x}*\tilde{t}_n +f_n\tilde{t}_n,\qquad \text{in }\left]-1,1\right[,
\]
we have \(\bigl\| u_n' \bigr\|_{L^\infty(-1,1)} = \bigl\| g' \bigr\|_{L^\infty(-1,1)}\), for all \(n\).}

As \((\tilde{t}_n,u_n)\in \cup_{p>1}L^p(-1,1)\times \cup_{p>1}W^{1,p}(-1,1)\) is the unique solution of problem~\(\mathscr{P}_{\!\rm o}\) associated with friction coefficient \(f_n\), the claim is given by Proposition~\ref{thm:eqPb}.

\smallskip
\noindent\textbf{Step 2.} \textit{We have:
\[
\tilde{t}_n = \frac{f_n\,u_n'\,\chi_{\left]-1,1\right[}}{1+f_n^2} + \underbrace{\frac{e^{\tau_n}\,\chi_{\left]-1,1\right[}}{\pi\sqrt{1-x^2}\sqrt{1+f_n^2}}}_{\displaystyle t_n^0}\;\text{\rm pv} \frac{1}{x} * \biggl[\frac{e^{-\tau_n}\,\sqrt{1-x^2}\,u_n'\,\chi_{\left]-1,1\right[}}{\sqrt{1+f_n^2}}\biggr],
\]
(where \(\chi_{\left]-1,1\right[}\) is the indicatrix function) which entails:
\[
-\frac{1}{\pi}\,\text{\rm pv}\frac{1}{x}*\tilde{t}_n +f_n\tilde{t}_n = \left|
\begin{array}{ll}
	\displaystyle u_n', & \text{in }\left]-1,1\right[,\\[1.0ex]
	\displaystyle -\frac{\text{\rm sgn}(x)\,e^{\tau_n}}{\pi\sqrt{x^2-1}}\;\text{\rm pv} \frac{1}{x} * \biggl[\frac{e^{-\tau_n}\,\sqrt{1-x^2}\,u_n'\,\chi_{\left]-1,1\right[}}{\sqrt{1+f_n^2}}\biggr], & \text{in }\mathbb{R}\setminus [-1,1].
\end{array}
\right.
\]
}

Set:
\[
T_n := \frac{f_n\,u_n'\,\chi_{\left]-1,1\right[}}{1+f_n^2} + \underbrace{\frac{e^{\tau_n}\,\chi_{\left]-1,1\right[}}{\pi\sqrt{1-x^2}\sqrt{1+f_n^2}}}_{\displaystyle t_n^0}\;\text{\rm pv} \frac{1}{x} * \biggl[\frac{e^{-\tau_n}\,\sqrt{1-x^2}\,u_n'\,\chi_{\left]-1,1\right[}}{\sqrt{1+f_n^2}}\biggr].
\]
As:
\[
\frac{e^{-\tau_n}\,\sqrt{1-x^2}\,u_n'\,\chi_{\left]-1,1\right[}}{\sqrt{1+f_n^2}}\in L^\infty(-1,1),\quad\text{ and }\quad \exists r\in\left]1,2\right[,\quad \frac{e^{\tau_n}\,\chi_{\left]-1,1\right[}}{\pi\sqrt{1-x^2}\sqrt{1+f_n^2}}\in L^r(-1,1),
\]
we obtain \(T_n\in L^r(-1,1)\), for some \(r\in\left]1,2\right[\), thanks to Theorem~\ref{thm:Riesz1} in Appendix~A.  In addition: 
\begin{align*}
	\int_{-1}^1 T_n\,{\rm d}x & = \int_{-1}^1 \frac{f_n\,u_n'\,\chi_{\left]-1,1\right[}}{1+f_n^2}\,{\rm d}x - \int_{-1}^1 \biggl[\frac{e^{-\tau_n}\,\sqrt{1-x^2}\,u_n'\,\chi_{\left]-1,1\right[}}{\sqrt{1+f_n^2}}\biggr]\biggl[\text{\rm pv} \frac{1}{x} * t_n^0\biggr]{\rm d}x,\\
	& = \int_{-1}^1 \frac{f_n\,u_n'\,\chi_{\left]-1,1\right[}}{1+f_n^2}\,{\rm d}x - \int_{-1}^1 \frac{f_n\,u_n'\,\chi_{\left]-1,1\right[}}{1+f_n^2}\,{\rm d}x = 0.
\end{align*}
Also, the Hilbert transform of \(T_n\) can be calculated, thanks to Corollary~\ref{thm:pbt}:
\begin{equation}
	\label{eq:labEqFla}
	-\frac{1}{\pi}\,\text{\rm pv}\frac{1}{x} * T_n = \frac{u_n'\,\chi_{\left]-1,1\right[}}{1+f_n^2} - \frac{1}{\pi}\Biggl(\text{\rm pv}\frac{1}{x} * \biggl[ \frac{e^{\tau_n}\,\chi_{\left]-1,1\right[}}{\pi\sqrt{1-x^2}\sqrt{1+f_n^2}} \biggr]\Biggr) \Biggl(\text{\rm pv}\frac{1}{x} * \biggl[ \frac{e^{-\tau_n}\,\sqrt{1-x^2}\,u_n'\,\chi_{\left]-1,1\right[}}{\sqrt{1+f_n^2}} \biggr]\Biggr).
\end{equation}
This yields, in particular:
\[
-\frac{1}{\pi}\,\text{\rm pv}\frac{1}{x} * T_n + f_nT_n = u_n',\qquad \text{a.e. in }\left]-1,1\right[, 
\]
so that \(T_n-\tilde{t}_n\) solves the homogeneous Carleman equation.  But, as \(f_n\) is Lipschitz-continuous, all the solutions in \(\cup_{p>1}L^p(-1,1)\) of the homogeneous Carleman equation are proportional to \(t_n^0\) (see \cite[Section 4-4]{tricomi} or Proposition~\ref{thm:defTi}).  Since both \(T_n\) and \(\tilde{t}_n\) have zero integral over \(\left]-1,1\right[\), we must have \(T_n=\tilde{t}_n\) and the claimed expression for \(\tilde{t}_n\) is therefore established.

Finally, by Proposition~\ref{thm:defTzero}, the Hilbert transform of \(t_n^0\) is:
\[
-\frac{1}{\pi}\,\text{\rm pv}\frac{1}{x} * \biggl[ \frac{e^{\tau_n}\,\chi_{\left]-1,1\right[}}{\pi\sqrt{1-x^2}\sqrt{1+f_n^2}} \biggr] = \left|
\begin{array}{ll}
	\displaystyle - \frac{f_n\,e^{\tau_n}}{\pi\sqrt{1-x^2}\sqrt{1+f_n^2}} , & \text{ in }\left]-1,1\right[,\\[3.0ex]
	\displaystyle - \frac{\text{\rm sgn}(x)\,e^{\tau_n}}{\pi\sqrt{x^2-1}} , & \text{ in }\mathbb{R}\setminus\left]-1,1\right[,
\end{array}
\right.
\]
which, together with identity~\eqref{eq:labEqFla}, is sufficient to obtain the claim in Step~2.

\smallskip
\noindent\textbf{Step 3.} \textit{We have:
\[
\biggl| -\frac{1}{\pi}\,\text{\rm pv}\frac{1}{x}*\tilde{t}_n +f_n\tilde{t}_n \biggr| \leq \left|
\begin{array}{ll}
	\displaystyle \bigl\|g'\bigr\|_{L^\infty(-1,1)} & \text{ in }\left]-1,1\right[,\\[3.0ex]
	\displaystyle \bigl\|g'\bigr\|_{L^\infty(-1,1)}\biggl[ 1 + \frac{e^{\tau_n(x)}}{\sqrt{x^2-1}}\bigl(|x|+2\pi\beta\bigr) \biggr] & \text{ in }\mathbb{R}\setminus\left]-1,1\right[,
\end{array}
\right.
\]
where \(\beta := (1/\pi)\arctan\|f_1\|_{L^\infty(-1,1)}\).}

The estimate in \(\left]-1,1\right[\) was already proved in Step~1.  Let us first estimate the following function for \(|x|>1\):
\begin{align*}
	\Biggl| \text{\rm pv} \frac{1}{x} * \biggl[\frac{e^{-\tau_n}\,\sqrt{1-x^2}\,u_n'\,\chi_{\left]-1,1\right[}}{\sqrt{1+f_n^2}}\biggr] \Biggr| & = \Biggl| \int_{-1}^1\frac{e^{-\tau_n(s)}\,\sqrt{1-s^2}\,u_n'(s)}{\sqrt{1+f_n^2(s)(x-s)}}{\rm d}s\Biggr| ,\\
	& \leq \bigl\|g'\bigr\|_{L^\infty} \biggl| \int_{-1}^1\frac{e^{-\tau_n(s)}\,\sqrt{1-s^2}}{\sqrt{1+f_n^2(s)(x-s)}}{\rm d}s\Biggr|.
\end{align*}
We have to estimate the integral in the right-hand side:
\begin{multline*}
	\int_{-1}^1\frac{e^{-\tau_n(s)}}{\sqrt{1+f_n^2(s)}}\,\sqrt{\frac{1-s}{1+s}}\,\frac{1+s}{x-s}\,{\rm d}s = (1+x)\int_{-1}^1\frac{e^{-\tau_n(s)}}{\sqrt{1+f_n^2(s)}}\,\sqrt{\frac{1-s}{1+s}}\,\frac{{\rm d}s}{x-s} - \mbox{} \\
	\mbox{} - \int_{-1}^1 \sqrt{\frac{1-s}{1+s}}\frac{e^{-\tau_n(s)}\,{\rm d}s}{\sqrt{1+f_n^2(s)}},
\end{multline*}
where the two integrals appearing in the right-hand side have actually already been calculated in the proof of Proposition~\ref{thm:defTzero}:
\begin{align*}
	\int_{-1}^1\frac{e^{-\tau_n(s)}}{\sqrt{1+f_n^2(s)}}\,\sqrt{\frac{1-s}{1+s}}\,\frac{1+s}{x-s}\,{\rm d}s 
	& = (1+x)\pi\biggl(1-e^{-\tau_n(x)}\sqrt{\biggl|\frac{1-x}{1+x}\biggr|}\biggr) - \pi \biggl(1-\int_{-1}^1\!\!\!\!\!\arctan f_n(s)\,{\rm d}s\biggr)\\
	& = -\text{sgn}(x)\,\pi\,\sqrt{x^2-1}\,e^{-\tau_n(x)} + \pi \biggl(x+\int_{-1}^1\arctan f_n(s)\,{\rm d}s\biggr).
\end{align*}
Using this and Step~2, we obtain the claim in Step~3.

\smallskip
\noindent\textbf{Step 4.} \textit{Conclusion.}

Set:
\[
V_n := \frac{e^{-\bar{\tau}_n}}{{\scriptstyle\sqrt{1+\bar{f}_n^2}}}\biggl( -\frac{1}{\pi}\,\text{\rm pv}\frac{1}{x}*\tilde{t}_n +f_n\tilde{t}_n \biggr).
\]
The set \(\text{supp}\bar{f}_n\) is bounded uniformly with respect to \(n\).  To fix ideas, we can suppose:
\[
\forall n\in \mathbb{N}\setminus\{0\},\qquad \text{supp}\bar{f}_n\subset [-2,2].
\]

Note that:
\[
\int_{-1}^1 |\tilde{t}_n(s)|\,{\rm d}s \leq \int_{-1}^1 \Bigl(|\tilde{t}_n(s)-Pt_n^0(s)|+Pt_n^0(s)\Bigr)\,{\rm d}s \leq 2P,
\]
and also that:
\[
\int_{-2}^2 |\arctan \bar{f}_n(s)| \,{\rm d}s \leq 4\pi\beta ,
\]
so that, for \(|x|>3\), we have:
\begin{align*}
	|V_n(x)| & = \frac{e^{-\bar{\tau}_n(x)}}{\pi}\biggl| \int_{-1}^1 \frac{\tilde{t}_n(s)\,{\rm d}s}{x-s} \biggr|,\\
	& \leq \frac{1}{\pi}\,\frac{\int_{-1}^1|\tilde{t}_n(s)|\,{\rm d}s}{|x|-1}\, e^{\frac{1}{\pi}\int_{-2}^2 \frac{\arctan\bar{f}_n(s)\,ds}{x-s}},\\
	& \leq \frac{2P}{\pi}\,\frac{1}{|x|-1}\,e^{\frac{4\beta}{|x|-2}}.
\end{align*}
The function appearing in the right-hand side is independent of \(n\) and is in \(L^p(-\infty,-3)\) and in \(L^p(3,\infty)\), for all \(p>1\).

Hence, to obtain the conclusion of the Lemma, it is now sufficient to prove that the sequence \(V_n\) is bounded in \(L^p(-3,3)\), for all \(p\in \left]1,(1/2+\beta)^{-1}\right[\).  To do this, we invoke the estimate of \(V_n\) provided by Step~3:
\[
\bigl|V_n(x)\bigr| \leq \left|
\begin{array}{ll}
	\displaystyle \bigl\|g'\bigr\|_{L^\infty(-1,1)} \frac{e^{-\bar{\tau}_n(x)}}{\scriptstyle\sqrt{1+\bar{f}_n^2(x)}} & \text{ in }\left]-1,1\right[,\\[3.0ex]
	\displaystyle \bigl\|g'\bigr\|_{L^\infty(-1,1)}\frac{e^{-\bar{\tau}_n(x)}}{\scriptstyle\sqrt{1+\bar{f}_n^2(x)}}\biggl[ 1 + \frac{e^{\tau_n(x)}}{\sqrt{x^2-1}}\bigl(|x|+2\pi\beta\bigr) \biggr] & \text{ in }\left]-3,-1\right[\cup\left]1,3\right[.
\end{array}
\right.
\]
By Lemma~\ref{thm:LpEst} in \(\left]-3,3\right[\) with \(f\) replaced by \(-\bar{f}\), the sequence \(\bigl(e^{-\bar{\tau}_n}/{\scriptstyle\sqrt{1+\bar{f}_n^2}}\bigr)\) is bounded in \(L^p(-3,3)\), for all \(p\in \left]1,(1/2+\beta)^{-1}\right[\).  Now, there remains only to prove that the sequence:
\begin{equation}
	\label{eq:boundedSeq}
	\biggl(\frac{e^{\tau_n(x)-\bar{\tau}_n(x)}}{{\scriptstyle\sqrt{1+\bar{f}_n^2(x)}}\sqrt{x^2-1}}\biggr)
\end{equation}
is bounded both in \(L^p(-3,-1)\) and \(L^p(1,3)\), for all \(p\in \left]1,(1/2+\beta)^{-1}\right[\), to reach the conclusion of the Lemma.  We are going to prove it for \(L^p(1,3)\), the other case being similar.  Set:
\[
\nu_n^{\text{\tiny\textcircled{+}}} := \bar{f}_n\,\chi_{\left]1,2\right[},\qquad \nu_n^{\text{\tiny\textcircled{-}}} := \bar{f}_n\,\chi_{\left]-2,-1\right[}
\]
so that \(\arctan f_n -\arctan\bar{f}_n = -\arctan\nu_n^{\text{\tiny\textcircled{-}}} - \arctan\nu_n^{\text{\tiny\textcircled{+}}}\) and \(\tau_n-\bar{\tau}_n=-\mathcal{H}[\arctan\nu_n^{\text{\tiny\textcircled{-}}}] - \mathcal{H}[\arctan\nu_n^{\text{\tiny\textcircled{+}}}]\) (where \(\mathcal{H}\) stands for the Hilbert transform).  As the sequence \(\arctan\nu_n^{\text{\tiny\textcircled{-}}}\) is supported in \([-2,-1]\) and is bounded in \(L^1(-2,-1)\), the sequence \(\mathcal{H}[\arctan\nu_n^{\text{\tiny\textcircled{-}}}]\) is bounded in \(C^0([1,3])\).  Finally, we can write:
\[
\text{for a.a. } x \in \left]1,3\right[,\qquad \frac{e^{\tau_n(x)-\bar{\tau}_n(x)}}{{\scriptstyle\sqrt{1+\bar{f}_n^2(x)}}\sqrt{x^2-1}} = C_n(x) \frac{e^{-\mathcal{H}[\arctan\nu_n^{\text{\tiny\textcircled{+}}}}](x)}{{\scriptstyle\sqrt{1+\nu_n^{{\text{\tiny\textcircled{+}}}2}(x)}}\sqrt{(3-x)(x-1)}},
\]
where \(C_n\) denotes a bounded sequence in \(C^0([1,3])\).  Invoking again Lemma~\ref{thm:LpEst}, the sequence \(e^{-\mathcal{H}[\arctan\nu_n^{\text{\tiny\textcircled{+}}}]}/({\scriptstyle\sqrt{1+\nu_n^{{\text{\tiny\textcircled{+}}}2}}}\sqrt{(3-x)(x-1)})\) is bounded in \(L^p(1,3)\), for all \(p\in \left]1,(1/2+\beta)^{-1}\right[\).  Hence we have proved that the sequence~\eqref{eq:boundedSeq} is bounded in \(L^p(1,3)\), for all \(p\in \left]1,(1/2+\beta)^{-1}\right[\).  By the same argument, it is also bounded in \(L^p(-3,-1)\).  All in all, we have proved that the sequence \((V_n)\) is bounded in \(L^p(-3,3)\), 
for all \(p\in \left]1,(1/2+\beta)^{-1}\right[\), and the Lemma is proved. \qed

\bigskip
\noindent\textbf{Proof of Theorem~\ref{thm:homogConvex}.}

\noindent\textbf{Step 1.} \textit{The sequence \((\tilde{t}_n)\) is bounded in \(\mathscr{M}([-1,1])\) and in \(H^{-1/2}(-1,1)\).  Therefore, upon extracting a subsequence, it converges towards some limit \(\bar{t}\) in \(\mathscr{M}([-1,1])\) weak-* and in \(H^{-1/2}(-1,1)\) weak.}

The integrable function \(\tilde{t}_n-Pt_n^0\) is nonpositive and has total mass \(P\).  This entails that \(\tilde{t}_n\) belongs to the closed unit ball of radius \(2P\) in \(\mathscr{M}([-1,1])\).  Also, as \(\tilde{t}_n\) solves problem~\(\mathscr{P}_{\!\rm o}\), it fulfils the variational inequality:
\[
\biggl\langle -\frac{1}{\pi}\,\log|x|*\tilde{t}_n + \frac{1}{2}\,\text{\rm sgn}(x)*\bigl(\bar{f}_n\tilde{t}_n\bigr) \;,\, \hat{t}-\tilde{t}_n\biggr\rangle \geq \bigl\langle g\;,\; \hat{t}-\tilde{t}_n\bigr\rangle,
\]
for all \(\hat{t}\in H_0^*\cap \mathscr{M}([-1,1])\) such that \(\hat{t}-Pt_n^0\leq 0\).  Taking \(\hat{t}=0\), we obtain:
\begin{align*}
	\biggl\langle -\frac{1}{\pi}\,\log|x|*\tilde{t}_n \;,\; \tilde{t}_n\biggr\rangle & \leq \bigl\| g\bigr\|_{H^{1/2}}\,\bigl\| \tilde{t}_n\bigr\|_{H^{-1/2}}  +  \bigl\| \tilde{t}_n\bigr\|_{\mathscr{M}}\,\bigl\| f_n\tilde{t}_n\bigr\|_{\mathscr{M}},
\end{align*}
where the second term in the right-hand side is bounded and the left-hand side can be identified with \(\|\tilde{t}_n\|_{H^{-1/2}}^2\), thanks to Lemma~\ref{thm:lemmHminusHalf} in Appendix~C.  This inequality therefore shows that \(\|\tilde{t}_n\|_{H^{-1/2}}\) is bounded.

\smallskip
\noindent\textbf{Step 2.} \textit{We have the following convergence:
\[
\frac{e^{-\bar{\tau}_n}}{\scriptstyle \sqrt{1+\bar{f}_n^2}}\biggl(-\frac{1}{\pi}\text{\rm pv}\frac{1}{x}*\tilde{t}_n+f_n\tilde{t}_n\biggr)  \rightharpoonup  \frac{e^{-\bar{\tau}_{\rm eff}}}{\scriptstyle \sqrt{1+\bar{f}_{\rm eff}^2}}\biggl(-\frac{1}{\pi}\text{\rm pv}\frac{1}{x}*\bar{t}+f_{\rm eff}\bar{t}\biggr)
\]
weakly in \(L^p(\mathbb{R})\), for all \(p\in \left]1,(1/2+\beta)^{-1}\right[\), where \(\beta\) is given by formula~\eqref{eq:beta}.}

From now on, we pick an arbitrary \(p\in\left]1,(1/2+\beta)^{-1}\right[\).  Given arbitrary \(t\in L^p(\mathbb{R})\) and \(y>0\), we have:
\[
\biggl\| \frac{1}{\pi}\int_{-\infty}^\infty \frac{y\,t(s)\,{\rm d}s}{(x-s)^2+y^2}\biggr\|_{L^p(\mathbb{R})} \leq \biggl\| \frac{1}{\pi}\frac{y}{x^2+y^2}\biggr\|_{L^1(\mathbb{R})} \|t\|_{L^p(\mathbb{R})}=\|t\|_{L^p(\mathbb{R})}.
\]
Applying Theorem~\ref{thm:Riesz2}\textit{(i)} in Appendix~A to the holomorphic function:
\[
\phi(z) := \frac{1}{\pi}\int_{-\infty}^\infty \frac{t(s)}{s-z}\,{\rm d}s,\qquad z\in \Pi^+,
\]
we obtain, for all \(y>0\):
\[
\frac{1}{\pi}\int_{-\infty}^{\infty}\frac{(s-x)\,t(s)}{(s-x)^2+y^2}\,{\rm d}s = -\frac{1}{\pi}\,\text{pv}\frac{1}{x} \stackrel{x}{*}\Biggl(\frac{1}{\pi}\int_{-\infty}^{\infty}\frac{y\,t(s)}{(s-x)^2+y^2}\,{\rm d}s\Biggr),
\]
so that, invoking Theorem~\ref{thm:Riesz1}\textit{(iii)} in Appendix~A, we have:
\[
\bigl\| \phi(x+iy) \bigr\|_{L^p(\mathbb{R})} \leq C \|t\|_{L^p(\mathbb{R})},
\]
for some universal constant \(C\), independent of \(t\), \(\phi\) and \(y>0\).  Applying Green's formula to the harmonic function \(h(x,y):= (1/\pi)\,\log\sqrt{x^2+y^2}\stackrel{x}{*}t\), we have, for all \(\varphi\in C_c^\infty(\overline{\Pi^+})\):
\[
\int_{-\infty}^{\infty} t(x)\,\varphi(x,0)\,{\rm d}x = -\frac{1}{\pi} \int_{\Pi^+}\!\biggl[\frac{\partial \varphi}{\partial x}(x,y)\int_{-\infty}^{\infty}\frac{(x-s)\,t(s)\,{\rm d}s}{(x-s)^2+y^2}+\frac{\partial \varphi}{\partial y}(x,y)\int_{-\infty}^{\infty}\frac{y\,t(s)\,{\rm d}s}{(x-s)^2+y^2}\biggr]{\rm d}x\,{\rm d}y,
\]
where, for all \(Y>0\), the following estimates hold:
\begin{equation}
	\label{eq:pairIdentities}
	\begin{aligned}
		\biggl\| \int_{-\infty}^{\infty}\frac{(x-s)\,t(s)\,{\rm d}s}{(x-s)^2+y^2}\biggr\|_{L^p(\mathbb{R}\times\left]0,Y\right[)} & \leq CY^{1/p}\,\| t\|_{L^p(\mathbb{R})},\\
		\biggl\| \int_{-\infty}^{\infty}\frac{y\,t(s)\,{\rm d}s}{(x-s)^2+y^2}\biggr\|_{L^p(\mathbb{R}\times\left]0,Y\right[)} & \leq CY^{1/p}\,\| t\|_{L^p(\mathbb{R})},
	\end{aligned}
\end{equation}
for some constant \(C\) independent of \(t\) and \(Y>0\).  When
\[
t =  \frac{e^{-\bar{\tau}_n}}{\scriptstyle \sqrt{1+\bar{f}_n^2}}\biggl(-\frac{1}{\pi}\text{\rm pv}\frac{1}{x}*\tilde{t}_n+f_n\tilde{t}_n\biggr),
\]
such a Green's formula is the one provided by Proposition~\ref{thm:fundCoroConvex}:
\begin{align*}
	\int_{-\infty}^\infty \frac{e^{-\bar{\tau}_n}}{\scriptstyle \sqrt{1+\bar{f}_n^2}}\biggl(-\frac{1}{\pi}\text{\rm pv}\frac{1}{x}*\tilde{t}_n+f_n\tilde{t}_n\biggr)\,\varphi(x,0)\,{\rm d}x 
	& =
	\int_{\Pi^+} \biggl[\Re e^{-\phi_n(x+iy)}\frac{\partial \Phi_{\tilde{\eta}_n}}{\partial x}+\Im e^{-\phi_n(x+iy)}\frac{\partial \Phi_{\tilde{\eta}_n}}{\partial y}\biggr]\frac{\partial\varphi}{\partial x}\\
	& \quad\mbox{} + \biggl[-\Im e^{-\phi_n(x+iy)}\frac{\partial \Phi_{\tilde{\eta}_n}}{\partial x}+\Re e^{-\phi_n(x+iy)}\frac{\partial \Phi_{\tilde{\eta}_n}}{\partial y}\biggr]\frac{\partial\varphi}{\partial y},
\end{align*}
where \(\tilde{\eta}_n:=(1/2)\,\text{sgn}(x)*\tilde{t}_n\) and:
\begin{align*}
	\phi_n(z) & := \frac{1}{\pi}\int_{-\infty}^\infty \frac{\arctan\bar{f}_n(s)\,{\rm d}s}{s-z},\\
	\Phi_{\tilde{\eta}_n}(x,y) & := \frac{1}{\pi}\int_{-\infty}^{\infty} \frac{y\,\tilde{\eta}_n\,{\rm d}s}{(x-s)^2+y^2}=\frac{1}{\pi}\, \arctan \frac{x}{y}\stackrel{x}{*}\tilde{t}_n.
\end{align*}
Using Lemma~\ref{thm:homogEstimate} and formula~\eqref{eq:pairIdentities}, this entails that the two sequences:
\[
\biggl(\Re e^{-\phi_n(x+iy)}\frac{\partial \Phi_{\tilde{\eta}_n}}{\partial x}+\Im e^{-\phi_n(x+iy)}\frac{\partial \Phi_{\tilde{\eta}_n}}{\partial y}\biggr)\qquad\text{and}\qquad \biggl(-\Im e^{-\phi_n(x+iy)}\frac{\partial \Phi_{\tilde{\eta}_n}}{\partial x}+\Re e^{-\phi_n(x+iy)}\frac{\partial \Phi_{\tilde{\eta}_n}}{\partial y}\biggr),
\]
are bounded in \(L^p(\mathbb{R}\times\left]0,Y\right[)\), for all \(Y>0\), and therefore have a weak limit (upon extracting subsequences) in \(L^p(\mathbb{R}\times\left]0,Y\right[)\).  As \((\arctan\bar{f}_n)\) converges towards \(\arctan \bar{f}_{\rm eff}\) in \(L^\infty(\mathbb{R})\) weak-*, \((\tilde{t}_n)\) converges towards \(\bar{t}\) in \(\mathscr{M}([-1,1])\) weak-* and \(H^{-1/2}(-1,1)\) weak,  and:
\begin{gather*}
	\phi_n(x+iy) = \frac{1}{\pi}\int_{-\infty}^\infty \frac{\arctan\bar{f}_n(s)\,{\rm d}s}{s-(x+iy)},\\
	\frac{\partial \Phi_{\tilde{\eta}_n}}{\partial x}(x,y)=\frac{1}{\pi}\int_{-1}^1 \frac{y\,\tilde{t}_n(s)\,{\rm d}s}{(x-s)^2+y^2}, \qquad \frac{\partial \Phi_{\tilde{\eta}_n}}{\partial y}(x,y)=\frac{1}{\pi}\int_{-1}^1 \frac{(s-x)\,\tilde{t}_n(s)\,{\rm d}s}{(x-s)^2+y^2},
\end{gather*}
it is readily checked that the sequences \((\phi_n(x+iy))\), \((\partial\Phi_{\tilde{\eta}_n}/\partial x)\) and \((\partial\Phi_{\tilde{\eta}_n}/\partial y)\) converge pointwisely in \(\Pi^+\) towards, respectively, \(\phi_{\rm eff}\), \(\partial\Phi_{\bar{\eta}}/\partial x\) and \(\partial \Phi_{\bar{\eta}}/\partial y\), where \(\bar{\eta}:=(1/2)\,\text{sgn}(x)*\bar{t}\) (this is the same argument already used in the case of the flat indentor).  Furthermore, these three sequences are bounded by a constant on each compact subset of \(\Pi^+\).  Using the dominated convergence theorem, this implies that the weak limits in \(L^p(\mathbb{R}\times\left]0,Y\right[)\) of the above two sequences must equal:
\[
\biggl(\Re e^{-\phi_{\rm eff}(x+iy)}\frac{\partial \Phi_{\bar{\eta}}}{\partial x}+\Im e^{-\phi_{\rm eff}(x+iy)}\frac{\partial \Phi_{\bar{\eta}}}{\partial y}\biggr)\qquad\text{and}\qquad \biggl(-\Im e^{-\phi_{\rm eff}(x+iy)}\frac{\partial \Phi_{\bar{\eta}}}{\partial x}+\Re e^{-\phi_{\rm eff}(x+iy)}\frac{\partial \Phi_{\bar{\eta}}}{\partial y}\biggr).
\]
This entails that the sequence:
\[
\Biggl(\frac{e^{-\bar{\tau}_n}}{\scriptstyle \sqrt{1+\bar{f}_n^2}}\biggl(-\frac{1}{\pi}\text{\rm pv}\frac{1}{x}*\tilde{t}_n+f_n\tilde{t}_n\biggr)\Biggr)_{n\in\mathbb{N}\setminus\{0\}}
\]
converges weakly in \(L^p(\mathbb{R})\) towards:
\[
\frac{e^{-\bar{\tau}_{\rm eff}}}{\scriptstyle \sqrt{1+\bar{f}_{\rm eff}^2}}\biggl(-\frac{1}{\pi}\text{\rm pv}\frac{1}{x}*\bar{t}+f_{\rm eff}\bar{t}\biggr).
\]

\smallskip
\noindent\textbf{Step 3.} \textit{We have \(\bar{t}=\tilde{t}_{\rm eff}\) and, in particular, the sequence \((\tilde{t}_n)\) converges towards \(\tilde{t}_{\rm eff}\) in \(\mathscr{M}([-1,1])\) weak-*.}

Step~2 and Lemma~\ref{thm:LpEst} entail that the sequence:
\[
\Biggl(\int_{-1}^x\frac{e^{-\bar{\tau}_n}}{\scriptstyle \sqrt{1+f_n^2}}\biggl(-\frac{1}{\pi}\text{\rm pv}\frac{1}{x}*\tilde{t}_n+f_n\tilde{t}_n - g' \biggr){\rm d}s\Biggr)_{n\in\mathbb{N}\setminus\{0\}}
\]
converges pointwisely towards:
\[
\int_{-1}^x\frac{e^{-\bar{\tau}_{\rm eff}}}{\scriptstyle \sqrt{1+f_{\rm eff}^2}}\biggl(-\frac{1}{\pi}\text{\rm pv}\frac{1}{x}*\bar{t}+f_{\rm eff}\bar{t} - g'\biggr){\rm d}s,
\]
and is bounded in \(W^{1,p}(-1,1)\), for all \(p\in \left]1,(1/2+\beta)^{-1}\right[\).  As \(W^{1,p}(-1,1)\) is continuously embedded in \(C^{0,(p-1)/p}([-1,1])\), Ascoli's theorem yields that the convergence is actually strong in \(C^0([-1,1])\).  Coming back to formula~\eqref{eq:defCn}, this entails that the sequence:
\[
C_n := \frac{1}{P} \,\biggl\langle Pt_n^0 - \tilde{t}_n\;,\; \int_{-1}^x\frac{e^{-\bar{\tau}_n}}{\scriptstyle \sqrt{1+f_n^2}}\biggl(-\frac{1}{\pi}\text{\rm pv}\frac{1}{x}*\tilde{t}_n+f_n\tilde{t}_n - g' \biggr)\biggr\rangle
\]
converges towards a limit \(C_\infty\), as the first term converges in \(\mathscr{M}([-1,1])\) weak-* and the second one converges strongly in \(C^0([-1,1])\).  Recalling that \(t^0_n\) converges weakly in \(L^p(-1,1)\) to \(t^0_{\rm eff}\), by Theorem~\ref{thm:Homog}, we therefore have:
\begin{itemize}
	\item \(\displaystyle \bigl\langle\bar{t},1\bigr\rangle =0\) and \(\bar{t}\leq Pt_{\rm eff}^0\).
	\item \(\displaystyle \int_{-1}^x \frac{e^{-\bar{\tau}_{\rm eff}}}{\scriptstyle \sqrt{1+f_{\rm eff}^2}}\biggl(-\frac{1}{\pi}\text{pv}\frac{1}{x}*\bar{t}+f_{\rm eff}\bar{t}-g'\biggr){\rm d}s\leq C_\infty\),
	\item \(\displaystyle \biggl\langle \bar{t}-Pt_{\rm eff}^0\;,\;\int_{-1}^x \frac{e^{-\bar{\tau}_{\rm eff}}}{\scriptstyle \sqrt{1+f_{\rm eff}^2}}\biggl(-\frac{1}{\pi}\text{pv}\frac{1}{x}*\bar{t}+f_{\rm eff}\tilde{t}_{\rm eff}-g'\biggr){\rm d}s-C_\infty\biggr\rangle = 0\).
\end{itemize}
But, there is one and only one element of \(H^{-1/2}(-1,1)\) satisfying all these conditions: \(\tilde{t}_{\rm eff}\).  Hence \(\bar{t}=\tilde{t}_{\rm eff}\) and we have proved that the sequence \((\tilde{t}_n)\) converges towards \(\tilde{t}_{\rm eff}\) in \(\mathscr{M}([-1,1])\) weak-* and in \(H^{-1/2}(-1,1)\) weak.

\smallskip
\noindent\textbf{Step 4.} \textit{The sequence \((f_n\tilde{t}_n)\) converges towards \(f_{\rm eff}\tilde{t}_{\rm eff}\) in \(\mathscr{M}([-1,1])\) weak-*.}

Let \((\tilde{t}_n,u_n)\in \cup_{p>1}L^p(-1,1)\times \cup_{p>1}W^{1,p}(-1,1)\) be the unique solution of problem~\(\mathscr{P}_{\!\rm o}\) associated with friction coefficient \(f_n\).  By Proposition~\ref{thm:eqPb}, we have:
\[
\bigl\| u_n' \bigr\|_{L^{\infty}(-1,1)} \leq \bigl\| g' \bigr\|_{L^{\infty}(-1,1)}.
\]
By Ascoli's theorem, the sequence \(u_n\) converges in \(C^0([-1,1])\) towards some limit \(\bar{u}\) (upon extracting a subsequence).  We are going to prove that \(\bar{u}=u_{\rm eff}\).  Of course, we have \(\bar{u}\leq g\).  As:
\[
\int_{-1}^1 \bigl(\tilde{t}_{\rm eff}-Pt_{\rm eff}^0\bigr)\bigl(\bar{u}-g\bigl)\,{\rm d}x = \lim_{n \rightarrow \infty}\int_{-1}^1 \bigl(\tilde{t}_n-Pt_{n}^0\bigr)\bigl(u_n-g\bigl)\,{\rm d}x = 0,
\]
we have \(\bar{u}=g=u_{\rm eff}\) on \(\text{supp} (\tilde{t}_{\rm eff}-Pt_{\rm eff}^0)\).  Consider an arbitrary \(x_0\in \left]-1,1\right[\setminus \text{supp} (\tilde{t}_{\rm eff}-Pt_{\rm eff}^0)\) (in the case where there exists some).  We denote by \((\tilde{t}_n,v_n)\in \cup_{p>1}L^p(-1,1)\times \cup_{p>1}W^{1,p}(-1,1)\) (resp. \((\tilde{t}_{\rm eff},v_{\rm eff})\)) the unique solution of problem~\(\mathscr{P}_{\!\rm a}\) associated with friction coefficient \(f_n\) (resp. \(f_{\rm eff}\)).  Set:
\begin{align*}
	w_n(x) & := v_n(x) - \int_{-1}^x \frac{g'\,e^{-\bar{\tau}_n}}{\scriptstyle \sqrt{1+f_n^2}}\,{\rm d}s - \sup_{x\in\left]-1,1\right[} \biggl\{v_n(x) - \int_{-1}^x \frac{g'\,e^{-\bar{\tau}_n}}{\scriptstyle \sqrt{1+f_n^2}}\,{\rm d}s \biggr\},\\
	& = \int_{-1}^x \frac{e^{-\bar{\tau}_n}}{\scriptstyle \sqrt{1+f_n^2}} \bigl(u_n'-g'\bigr)\, {\rm d}s - \sup_{x\in\left]-1,1\right[} \int_{-1}^x \frac{e^{-\bar{\tau}_n}}{\scriptstyle \sqrt{1+f_n^2}} \bigl(u_n'-g'\bigr)\, {\rm d}s,\\
	& = \int_{-1}^x \frac{e^{-\bar{\tau}_n}}{\scriptstyle \sqrt{1+f_n^2}} \biggl(-\frac{1}{\pi}\text{\rm pv}\frac{1}{x}*\tilde{t}_n + f_n\tilde{t}_n - g'\biggr)\, {\rm d}s -\mbox{} \\
	& \hspace*{5cm} \mbox{} - \sup_{x\in\left]-1,1\right[} \int_{-1}^x \frac{e^{-\bar{\tau}_n}}{\scriptstyle \sqrt{1+f_n^2}} \biggl(-\frac{1}{\pi}\text{\rm pv}\frac{1}{x}*\tilde{t}_n + f_n\tilde{t}_n - g'\biggr)\, {\rm d}s, 
\end{align*}
and \(w_{\rm eff}\) by replacing all the `$n$' by `eff' in the above definition.  By Proposition~\ref{thm:eqPb}, the function \(w_{\rm eff}\) vanishes on \(\text{supp} (\tilde{t}_{\rm eff}-Pt_{\rm eff}^0)\) and is negative on \(\left]-1,1\right[\setminus \text{supp} (\tilde{t}_{\rm eff}-Pt_{\rm eff}^0)\).  Hence, \(w_{\rm eff}(x_0) < 0\).  By the analysis led in step~2 and step~3, the sequence \((w_n)\) converges towards \(w_{\rm eff}\) in \(C^0([-1,1])\).  Therefore:
\[
\exists N\in\mathbb{N},\quad\exists\eta>0,\quad\forall n\geq 0, \quad\forall x\in [x_0-2\eta,x_0+2\eta],\qquad w_n(x) < 0.
\]
Hence, for all \(n\geq N\), \([x_0-\eta,x_0+\eta]\cap \text{supp} (\tilde{t}_{n}-Pt_{n}^0)=\varnothing\), and, of course \([x_0-\eta,x_0+\eta]\cap \text{supp} (\tilde{t}_{\rm eff}-Pt_{\rm eff}^0)=\varnothing\).  As, \(\text{supp} (\tilde{t}_{\rm eff}-Pt_{\rm eff}^0)\) is connected, it lies either on the left, or the right of 
\([x_0-\eta,x_0+\eta]\).  Let us assume that it lies on the right.  Then, the same is true of \(\text{supp} (\tilde{t}_{n}-Pt_{n}^0)\) whenever \(n\) is large enough, by the weak-* convergence in \(\mathscr{M}([-1,1])\) of \((\tilde{t}_{n}-Pt_{n}^0)\).  Then:
\[
\exists U_n\in \mathbb{R},\quad\forall x\in \left]-1,x_0+\eta\right[,\qquad u_n(x) = U_n -\frac{1}{\pi}\int_{x_0+\eta}^{1} \log|x-s|\,\Bigl(\tilde{t}_{n}(s)-Pt_{n}^0(s)\Bigr)\,{\rm d}s,
\]
where the sequence \((U_n)\) must be bounded.  We can therefore take the limit to obtain:
\[
\exists U\in \mathbb{R},\quad \forall x\in \left]-1,x_0+\eta\right[,\qquad \bar{u}(x) = U + u_{\rm eff}(x),
\]
and therefore:
\[
\exists U\in \mathbb{R},\quad \forall x < \inf \text{supp} (\tilde{t}_{\rm eff}-Pt_{\rm eff}^0),\qquad \bar{u}(x) = U + u_{\rm eff}(x).
\]
As \(\bar{u}\) and \(u_{\rm eff}\) are continuous and equal on \(\text{supp} (\tilde{t}_{\rm eff}-Pt_{\rm eff}^0)\), we have \(U=0\).  All in all, we have proved \(\bar{u} = u_{\rm eff}\).

The sequence \((f_n\tilde{t}_n)\) being bounded in \(\mathscr{M}([-1,1])\) converges in \(\mathscr{M}([-1,1])\) weak-*, upon extracting a subsequence, towards some limit \(\bar{\bar{t}}\in\mathscr{M}([-1,1])\).  We are going to prove that \(\bar{\bar{t}}=f_{\rm eff}\tilde{t}_{\rm eff}\).  Going to the limit in the formula:
\[
-\frac{1}{\pi}\text{\rm pv}\frac{1}{x}*\tilde{t}_n + f_n\tilde{t}_n = u_n' ,\qquad \text{in } \mathscr{D}'(\left]-1,1\right[),
\]
we obtain:
\[
-\frac{1}{\pi}\text{\rm pv}\frac{1}{x}*\tilde{t}_{\rm eff} + \bar{\bar{t}} = u_{\rm eff}' = 	-\frac{1}{\pi}\text{\rm pv}\frac{1}{x}*\tilde{t}_{\rm eff} + f_{\rm eff}\tilde{t}_{\rm eff} ,\qquad \text{in } \mathscr{D}'(\left]-1,1\right[).
\]
Hence, we have \(\bar{\bar{t}}=f_{\rm eff}\tilde{t}_{\rm eff}\) as elements of \(\mathscr{D}'(\left]-1,1\right[)\).  To obtain the equality as elements of \(\mathscr{M}([-1,1])\), we need to prove that the measure \(\bar{\bar{t}}\) has no atom at \(-1\) and \(1\).  But this is ensured by:
\[
\|f_1\|_{L^\infty(-1,1)} \bigl(\tilde{t}_n-Pt_n^0\bigr) \leq f_n\bigl(\tilde{t}_n-Pt_n^0\bigr) \leq - \|f_1\|_{L^\infty(-1,1)} \bigl(\tilde{t}_n-Pt_n^0\bigr),
\]
which yields:
\[
\|f_1\|_{L^\infty(-1,1)} \bigl(\tilde{t}_{\rm eff}-Pt_{\rm eff}^0\bigr) \leq \bar{\bar{t}} - Pf_{\rm eff}t_{\rm eff}^0 \leq - \|f_1\|_{L^\infty(-1,1)} \bigl(\tilde{t}_{\rm eff}-Pt_{\rm eff}^0\bigr),
\]
in \(\mathscr{M}([-1,1])\), as \((t_n^0)\) and \((f_{n}t_n^0)\) converge respectively towards \(t_{\rm eff}^0\) and \(f_{\rm eff}t_{\rm eff}^0\) in \(\mathscr{M}([-1,1])\) weak-*, thanks to Theorem~\ref{thm:Homog}.  All in all, we have proved that the sequence \((f_n\tilde{t}_n)\) converges towards \(f_{\rm eff}\tilde{t}_{\rm eff}\) in \(\mathscr{M}([-1,1])\) weak-*.
\qed

\section*{Appendix A: Marcel Riesz's \(L^p\)-theory of the Hilbert transform}

The distributional derivative of the locally integrable function \(\log|x|\) is denoted by \(\text{pv}\,1/x\) (meaning `principal value').  It satisfies:
\[
\forall\varphi\in C_c^\infty(\mathbb{R}),\qquad \Bigl\langle \text{pv}\,\frac{1}{x},\varphi\Bigr\rangle = \lim_{\varepsilon \rightarrow 0+} \int_{\mathbb{R}\setminus\left]-\varepsilon,\varepsilon\right[} \frac{\varphi(x)}{x}\,{\rm d}x,
\]
(\(C_c^\infty\) stands for the space of \(C^\infty\) compactly supported test-functions), where the limit exists thanks to the differentiability of \(\varphi\) at \(0\).  The distribution \(\text{pv}\,1/x\) is a tempered distribution whose Fourier transform is given by:
\[
\mathcal{F}\bigl[\text{pv}\,1/x\bigr] = -i \sqrt{\frac{\pi}{2}}\,\text{sgn},
\]
where \(\text{sgn}\) is the sign function and the following convention is used to define the Fourier transform of any integrable function \(f\):
\[
\mathcal{F}\bigl[f\bigr](\xi) := \frac{1}{\sqrt{2\pi}}\int_{-\infty}^\infty f(x)\,e^{-ix\xi}\,{\rm d}x.
\]
As the Fourier transform is an isometry of \(L^2(\mathbb{R})\), the same applies to the so-called Hilbert transform defined by the convolution product:
\[
\mathcal{H}\bigl[f\bigr] = -\frac{1}{\pi}\,\text{pv}\frac{1}{x}*f,
\]
whose inverse is minus itself:
\[
\forall f\in L^2(\mathbb{R}),\qquad\mathcal{H}\bigl[\mathcal{H}[f]\bigr] = -f.
\]

Actually, these results classically extend to \(L^p(\mathbb{R})\), provided that \(p\in \left]1,+\infty\right[\).  A proof of the following two theorems of Marcel Riesz can be found for example in~\cite{titch}.

\begin{theo}
	[M. Riesz] 
	\label{thm:Riesz1}
	Let \(p\in\left]1,\infty\right[\) and \(f\in L^p(\mathbb{R})\). 
	Then:
	\begin{itemize}
		\item[(i)] For almost all \(x\in\mathbb{R}\), 
		\[
		\mathcal{H}[f](x) := -\frac{1}{\pi}\lim_{\varepsilon\rightarrow 0+} 
		\int_{\mathbb{R}\setminus\left]x-\varepsilon,x+\varepsilon\right[} \frac{f(s)}{x-s}\,\,{\rm d}s ,
		\]
		exists and is finite.
		\item[(ii)] The function \(\mathcal{H}[f]\), defined almost everywhere in \(\mathbb{R}\) by the above formula, is in \(L^p(\mathbb{R})\).
		\item[(iii)] The mapping \(\mathcal{H}:L^p \rightarrow L^p\) is linear continuous.
		\item[(iv)] The following identity holds:
		\[
		\forall f\in L^p(\mathbb{R}),\qquad\mathcal{H}\bigl[\mathcal{H}[f]\bigr] = -f,
		\]
		which shows that \(\mathcal{H}\) is an isomorphism of \(L^p(\mathbb{R})\).
		\item[(v)] The function \(\mathcal{H}[f]\) coincides with the convolution product:
		\[
		\mathcal{H}\bigl[f\bigr] = -\frac{1}{\pi}\,\text{\rm pv}\frac{1}{x}*f,
		\]
		and is accordingly named the Hilbert transform of \(f\in L^p(\mathbb{R})\) (\(p\in\left]1,+\infty\right[\)).
	\end{itemize}
\end{theo}

The Hilbert transform is closely related to the study of some class of 
holomorphic functions defined in the open upper-half plane:
\[
\Pi^+ := \Bigl\{ z\in\mathbb{C} \, ; \, \Im (z) >0 \Bigr\},
\]
as seen in the next theorem. 

\begin{theo}
	[M. Riesz] 
	\label{thm:Riesz2}
	\mbox{}
	\begin{itemize}
		\item[(i)] Let \(p\in\left]1,\infty\right[\) and \(f\in L^p(\mathbb{R})\).  Then, the function 
		\(\Phi : \Pi^+ \rightarrow \mathbb{C}\) defined by:
		\[
		\Phi (z) := \frac{1}{\pi}\int_{-\infty}^{\infty} \frac{f(t)}{t-z}\,\,{\rm d}t  ,
		\]
		is holomorphic in \(\Pi^+\) and satisfies:
		\begin{equation}
			\label{defHardy}
			\exists\,K \in \mathbb{R},\quad \forall\, y>0,\qquad \int_{-\infty}^{\infty} \bigl| 
			\Phi(x+iy) \bigr|^p \, {\rm d}x < K  .
		\end{equation}
		In addition, \(\Phi(x+iy)\) converges, as \(y\rightarrow 0+\), in \(L^p(\mathbb{R};\mathbb{C})\), but 
		also for almost all \(x\in \mathbb{R}\), towards the limit:
		\[
		\Phi(x+i0) = \mathcal{H}[f](x) + i f(x).
		\]
		\item[(ii)] Reciprocally, if \(\Phi\) is some holomorphic function in \(\Pi^+\) satisfying the 
		condition: 
		\[
		\exists\,K \in \mathbb{R},\quad \forall\, y>0,\qquad \int_{-\infty}^{\infty} \bigl| 
		\Phi(x+iy) \bigr|^p \, {\rm d}x < K ,
		\]
		then, \(\Phi(x+iy)\) converges, as \(y\rightarrow 0+\), in \(L^p(\mathbb{R};\mathbb{C})\), but 
		also for almost all \(x\in \mathbb{R}\), towards some limit 
		\(\Phi(x+i0)\), which satisfies in 
		addition:
		\[
		\Re \bigl\{\Phi(x+i0)\bigr\} =  \mathcal{H} \Bigl[\Im \bigl\{\Phi(x+i0)\bigr\}\Bigr] .
		\]
	\end{itemize}
\end{theo}

Using Theorem~\ref{thm:Riesz2}(i) and the following indentities:
\begin{align*}
	\Re \Bigl(\bigl(\mathcal{H}[f_1]+if_1\bigr)\bigl(\mathcal{H}[f_2]+if_2\bigr)\Bigr) &= \mathcal{H}[f_1]\,\mathcal{H}[f_2]-f_1\,f_2,\\
	\Im \Bigl(\bigl(\mathcal{H}[f_1]+if_1\bigr)\bigl(\mathcal{H}[f_2]+if_2\bigr)\Bigr) &= \mathcal{H}[f_1]\,f_2 + f_1\,\mathcal{H}[f_2],
\end{align*}
we get the following corollary.

\begin{coro}
	[Poincaré-Bertrand-Tricomi]
	\label{thm:pbt}
	Let \(f_{1}\in L^{p_{1}}(\mathbb{R})\) and \(f_{2}\in 
	L^{p_{2}}(\mathbb{R})\), with \(p_1,p_2\in \left]1,\infty\right[\) and \(\frac{1}{p_{1}}+\frac{1}{p_{2}}<1\).  Then:
	\[
	\mathcal{H}\Bigl[ \mathcal{H}[f_{1}]\,f_{2} + f_{1}\,\mathcal{H}[f_{2}] \Bigr] = 
	\mathcal{H}[f_{1}]\,\mathcal{H}[f_{2}] - f_{1}\, f_{2} .
	\]
\end{coro}

\begin{coro}
	\label{thm:noteTricomi}
	Let \(\Phi:\Pi^+\rightarrow\mathbb{C}\) be a holomorphic function such that \(\Phi(x+iy)\) converges, as \(y\rightarrow 0+\), towards a limit \(\Phi(x+i0)\), for almost all \(x\in \mathbb{R}\).  Suppose that \(\Phi(x+i0)\in L^p(\mathbb{R};\mathbb{C})\), for some \(p\in \left]1,\infty\right[\) and \(\Phi\) is bounded at infinity.  Then:
	\[
	\Re \bigl\{\Phi(x+i0)\bigr\} =  \mathcal{H} \Bigl[\Im \bigl\{\Phi(x+i0)\bigr\}\Bigr] .
	\]
\end{coro}

\noindent\textbf{Proof.} Let \(f:=\Im\Phi(x+i0)\in L^p(\mathbb{R})\).  Then, the function \(\tilde{\Phi}(z):=\Phi(z)-\frac{1}{\pi}\int_{-\infty}^\infty\frac{f(t)}{t-z}\,{\rm d}t\) is holomorphic in \(\Pi^+\).  By Theorem~\ref{thm:Riesz2}(i), \(\Im\tilde{\Phi}(x+iy)\) converges, as \(y\rightarrow0+\), towards \(0\), for almost all \(x\in \mathbb{R}\).  This fact can be used to extend the harmonic function \((x,y)\mapsto \Im\tilde{\Phi}(x+iy)\) as a harmonic function defined on the whole plane (by reflection).  As this harmonic extension is bounded at infinity, it is bounded and is therefore constant by the Liouville theorem for harmonic functions.  Hence, the holomorphic function \(\tilde{\Phi}\) must be constant on \(\Pi^+\).  This constant must equal \(\Re\tilde{\Phi}(x+i0)\in L^p(\mathbb{R})\) and is therefore \(0\).  Hence, \(\Phi(z)=\frac{1}{\pi}\int_{-\infty}^\infty\frac{f(t)\,{\rm d}t}{t-z}\) and the claim is now obtained from Theorem~\ref{thm:Riesz2}(i). \qed

\section*{Appendix B: the spaces $H^{1/2}$ and $H^{-1/2}$}

Let \(I\) be an arbitrary (bounded or unbounded) open interval in \(\mathbb{R}\).

\begin{defi}
	The space:
	\[
	H^{1/2}(I) := \biggl\{u\in L^2(I) \;\Bigm|\; \int_{I\times I} \Bigl(\frac{u(x)-u(y)}{x-y}\Bigr)^2{\rm d}x\,{\rm d}y < \infty\biggr\},
	\]
	endowed with the norm:
	\[
	\bigl\| u\bigr\|_{H^{1/2}}^2 := \bigl\| u\bigr\|_{L^{2}}^2 + \int_{I\times I} \Bigl(\frac{u(x)-u(y)}{x-y}\Bigr)^2{\rm d}x\,{\rm d}y
	\]
	is a Hilbert space.  The space \(C_c^\infty(I)\) of \(C^\infty\) functions with compact support in \(I\) is dense in \(H^{1/2}\).  Its dual space is therefore a space of distributions in \(I\).  It is denoted by \(H^{-1/2}(I)\).
\end{defi}

It follows directly from the definition of \(H^{1/2}\) that the restriction to \(\left]c,d\right[\subset\left]a,b\right[\) of some \(u\in H^{1/2}(\left]a,b\right[)\) defines an element of \(H^{1/2}(\left]c,d\right[)\).  As a corollary, the extension by zero of some \(t\in H^{-1/2}(\left]c,d\right[)\) defines an element of \(H^{-1/2}(\left]a,b\right[)\).  However, the extension by zero of an arbitrary element of \(H^{1/2}(\left]c,d\right[)\) is not, in general, an element of \(H^{1/2}(\left]a,b\right[)\).  As a corollary, the restriction to \(\left]c,d\right[\) of some \(t\in H^{-1/2}(\left]a,b\right[)\) is not, in general, in \(H^{-1/2}(\left]c,d\right[)\).

In the particular case where \(I=\mathbb{R}\), the spaces \(H^{1/2}(\mathbb{R})\) and \(H^{-1/2}(\mathbb{R})\) can be equivalently defined in terms of the Fourier transform.

\begin{prop}
	\label{thm:defSobolevFourier}
	Denoting by \(\mathscr{S}'(\mathbb{R})\) the space of tempered distributions and by \(\hat{t}=\mathcal{F}[t]\) the Fourier transform of an arbitrary tempered distribution \(t\in \mathscr{S}'(\mathbb{R})\), we have:
	\begin{align*}
		H^{1/2}(\mathbb{R}) & = \Bigl\{ u\in 
		\mathscr{S}'(\mathbb{R})\;\bigm|\; (1+|\xi|^2)^{1/4}\,\hat{u}(\xi) \in L^2(\mathbb{R};\mathbb{C})\Bigr\},\\
		H^{-1/2}(\mathbb{R}) & = \Bigl\{ t\in 
		\mathscr{S}'(\mathbb{R})\;\bigm|\; (1+|\xi|^2)^{-1/4}\,\hat{t}(\xi)\in L^2(\mathbb{R};\mathbb{C})\Bigr\},
	\end{align*}
	and the corresponding natural norms are equivalent to those of \(H^{1/2}(\mathbb{R})\) and \(H^{-1/2}(\mathbb{R})\).  In particular, the Hilbert transform is an isomorphism of both \(H^{1/2}(\mathbb{R})\) and \(H^{-1/2}(\mathbb{R})\).
\end{prop}

\begin{theo}
	[Sobolev embeddings]
	Assume that $a$ and $b$ are finite.  Then:
	\[
	\forall p\in \left[1,\infty\right[,\qquad H^{1/2}(a,b)\subset L^p(a,b),
	\]
	with continuous embedding.  On the dual side:
	\[
	\forall p\in \left]1,\infty\right],\qquad  L^p(a,b)\subset H^{-1/2}(a,b),
	\]
	with continuous embeddings.
\end{theo}

The space \(H^{1/2}(-1,1)\) contains unbounded functions such as \(\log|\log|x/2||\), but not functions with jumps such as the Heaviside function.  As a corollary, a measure \(t\in \mathscr{M}([-1,1])\cap H^{-1/2}(-1,1)\) has no atoms.

\section*{Appendix C: Fourier transform and Poisson integral}

The following convention is chosen to define the Fourier transform of any integrable function \(f\):
\[
\mathcal{F}\bigl[f\bigr](\xi) := \frac{1}{\sqrt{2\pi}}\int_{-\infty}^\infty f(x)\,e^{-ix\xi}\,{\rm d}x.
\]

\begin{prop}
	\label{thm:FourierTransforms}
	We denote by \(\text{\rm pv}1/x\) (`principal value') the distributional derivative of \(\log|x|\) and by \(\text{\rm fp}1/|x|\) (`finite part') that of \(\text{\rm sgn}(x)\,\log|x|\) ({\rm sgn} is the sign function).  Then, the Fourier transforms of the following distributions are known explicitly.
	\begin{align*}
		\mathcal{F}\Bigl[\log\sqrt{x^2+a^2}\Bigr] & = -\sqrt{\frac{\pi}{2}}\,e^{-a|\xi|}\,\text{\rm fp}\frac{1}{|\xi|} - \sqrt{2\pi}\,\Gamma_{\rm Eul}\,\delta, \qquad & 
		\mathcal{F}\Bigl[{\scriptstyle \frac{2}{\pi}}\arctan (x/a)\Bigr] & = - i \sqrt{\frac{2}{\pi}}\,e^{-a|\xi|}\,\text{\rm pv}\frac{1}{\xi},\\
		\mathcal{F}\Bigl[\log|x|\Bigr] & = -\sqrt{\frac{\pi}{2}}\,\text{\rm fp}\frac{1}{|\xi|} - \sqrt{2\pi}\,\Gamma_{\rm Eul}\,\delta, &
		\mathcal{F}\Bigl[\text{\rm sgn}(x)\Bigr] & = - i \sqrt{\frac{2}{\pi}}\,\text{\rm pv}\frac{1}{\xi},
	\end{align*}
	where \(a> 0\) is an arbitrary positive real constant, \(\delta\) the Dirac measure at \(0\) and \(\Gamma_{\rm Eul}:=\lim_{n\rightarrow\infty}-\log n+\sum_{k=1}^n 1/k\) is the Euler-Mascheroni constant.
\end{prop}

As a consequence of the above Fourier transforms, we have the following lemma whose proof is to be found in \cite[Theorem 3]{bibiJiri}.

\begin{lemm}
	\label{thm:lemmHminusHalf}
	For all \(t\in H^{-1/2}(-1,1)\) extended by zero on \(\mathbb{R}\), the convolution product \(\log|x|*t\) is a locally square integrable function, whose restriction to \(\left]-1,1\right[\) is in \(H^{1/2}(-1,1)\).  The symmetric bilinear form:
	\[
	t_1,t_2 \mapsto -\bigl\langle \log|x|*t_1,t_2\bigr\rangle,
	\]
	is a scalar product on the Hilbert space \(H^{-1/2}(-1,1)\) which induces a norm that is equivalent to that of \(H^{-1/2}\).
\end{lemm}

Given an arbitrary \(u\in L^p(\mathbb{R})\), the function:
\begin{equation}
	\label{eq:defPoissonIntegral}
	\Phi_u (x,y) := \frac{1}{\pi} \int_{-\infty}^\infty \frac{y\,u(s)}{(x-s)^2+y^2}\,{\rm d}s = \frac{1}{\pi}\frac{y}{x^2+y^2}\stackrel{x}{*}u = \frac{1}{\pi}\arctan \frac{x}{y} * u',
\end{equation}
is a well-defined \emph{harmonic} function in \(\Pi^+\), called the \emph{Poisson integral} of \(u\).  The following proposition gives some of its properties in the case where \(u\in H^{1/2}(\mathbb{R})\).

\begin{prop}
	\label{thm:prolAnal}
	Let \(u\in H^{1/2}(\mathbb{R})\) be arbitrary.  Then, its Poisson integral \(\Phi_u\) is harmonic on \(\Pi^+:=\mathbb{R}\times\left]0,+\infty\right[\), \(\nabla\Phi_u\in L^2(\Pi^+)\) and \(\Phi_u\in H^1(\mathbb{R}\times\left]0,Y\right[)\), for all \(Y>0\).  Moreover, \(u\) coincides with the trace of \(\Phi_u\) on \(\partial \Pi^+\simeq \mathbb{R}\) and the Hilbert transform of \(u'\) coincides with the derivative \(\partial \Phi_u/\partial y\) on the boundary, in the sense:
	\[
	\forall \varphi\in H^1(\Pi^+),\qquad \int_{\Pi^+} \biggl(\frac{\partial\Phi_u}{\partial x}\frac{\partial\varphi}{\partial x} + \frac{\partial\Phi_u}{\partial y}\frac{\partial\varphi}{\partial y}\biggr){\rm d}x\,{\rm d}y = \biggl\langle \frac{1}{\pi} \, \text{\rm pv}\frac{1}{x} * u'\;,\; \varphi(x,0) \biggr\rangle.
	\]
\end{prop}

\noindent\textbf{Proof.} First, note that, at fixed \(y>0\), the function \(x\mapsto y/(x^2+y^2)\) is in \(L^1(\mathbb{R})\).  Therefore, its convolution product with the function \(u\) is well-defined and the function \(x\mapsto \Phi_u(x,y)\) is in \(L^2(\mathbb{R})\).  As the function \((x,y)\mapsto y/(x^2+y^2)\) is harmonic in \(\Pi^+\), the same is true of the function \(\Phi_u\).

By Proposition~\ref{thm:FourierTransforms}, the Fourier transform \(\hat{\Phi}_u(\xi,y)\) of \(\Phi_u\) with respect to the $x$-variable is given by:
\[
\hat{\Phi}_u(\xi,y) = \sqrt{2\pi}\,\frac{e^{-y|\xi|}}{\sqrt{2\pi}}\,\hat{u}(\xi) = e^{-y|\xi|}\,\hat{u}(\xi),
\]
which shows \(\Phi_u\in C^0([0,Y],L^2(\mathbb{R}))\subset L^2(\mathbb{R}\times\left]0,Y\right[)\), for all \(Y>0\).  By:
\[
\Bigl|-i\xi\,\hat{\Phi}_u(\xi,y)\Bigr|^2 + \biggl| \frac{\partial \hat{\Phi}_u}{\partial y}(\xi,y) \biggr|^2 = 2e^{-2y|\xi|}\,|\xi|^2|\hat{u}(\xi)|^2,
\]
we obtain \(\nabla\Phi_u\in L^2(\Pi^+)\) and:
\begin{equation}
	\label{eq:estExt}
	\|\nabla\Phi_u\|_{L^2(\Pi^+)}^2 =\int_{-\infty}^{+\infty}|\xi|\,|\hat{u}(\xi)|^2\,{\rm d}\xi\leq C \| u\|_{H^{1/2}(\mathbb{R})}^2,
\end{equation}
for some real constant \(C\), independent of \(u\).  Hence, we have proved that \(\Phi_u\in H^1(\mathbb{R}\times \left]0,Y\right[)\) and that its trace on \(\partial \Pi^+\simeq \mathbb{R}\) is \(u\).  Taking a sequence \((u_n)\) in \(C_c^\infty(\mathbb{R})\) that converges strongly towards \(u\) in \(H^{1/2}(\mathbb{R})\), the sequence \((\nabla\Phi_{u_n})\) converges strongly towards \(\nabla\Phi_u\) in \(L^2(\Pi^+)\), thanks to estimate~\eqref{eq:estExt}.  But:
\[
\frac{\partial \hat{\Phi}_{u_n}}{\partial y}(\xi,y) = - e^{-y|\xi|}\,|\xi|\,\hat{u}_n(\xi),
\]
which shows that \(\Phi_{u_n}\in C^\infty(\overline{\Pi^+})\), and also that:
\begin{align*}
	\frac{\partial \hat{\Phi}_{u_n}}{\partial y}(\xi,0) & = - \frac{\sqrt{2\pi}}{\pi}\Bigl(i\xi\,\hat{u}_n(\xi)\Bigr)\biggl(-i\sqrt{\frac{\pi}{2}}\,\text{sgn}(\xi)\biggr),\\
	& = - \frac{\sqrt{2\pi}}{\pi} \; \widehat{u_n'}(\xi) \; \widehat{\text{pv}\,1/x}(\xi),
\end{align*}
which is nothing but:
\[
\frac{\partial \Phi_{u_n}}{\partial y}(x,0) = - \frac{1}{\pi}\,\text{pv}\,\frac{1}{x}*u_n'.
\]
By Green's formula, this gives:
\[
\forall \varphi\in C_c^\infty(\overline{\Pi^+}),\qquad \int_{\Pi^+} \biggl(\frac{\partial\Phi_{u_n}}{\partial x}\frac{\partial\varphi}{\partial x} + \frac{\partial\Phi_{u_n}}{\partial y}\frac{\partial\varphi}{\partial y}\biggr){\rm d}x\,{\rm d}y = \biggl\langle \frac{1}{\pi} \, \text{\rm pv}\frac{1}{x} * u_n'\;,\; \varphi(x,0) \biggr\rangle,
\]
and, by taking the limit \(n \rightarrow +\infty\) using the fact that the Hilbert transform maps continuously \(H^{-1/2}(\mathbb{R})\) onto itself:
\[
\forall \varphi\in C_c^\infty(\overline{\Pi^+}),\qquad \int_{\Pi^+} \biggl(\frac{\partial\Phi_u}{\partial x}\frac{\partial\varphi}{\partial x} + \frac{\partial\Phi_u}{\partial y}\frac{\partial\varphi}{\partial y}\biggr){\rm d}x\,{\rm d}y = \biggl\langle \frac{1}{\pi} \, \text{\rm pv}\frac{1}{x} * u'\;,\; \varphi(x,0) \biggr\rangle.
\]
The expected identity now follows from the density of \(C_c^\infty(\overline{\Pi^+})\) in \(H^1(\Pi^+)\).\qed

\section*{Appendix D: pseudomonotone variational inequalities}

The content of this appendix particularizes some more general results from \cite{BrezisIeq}, for easy reference and the convenience of the reader.

Let \(E\) be a separable Banach space and \(E^*\) be its (topological) dual space, the duality product being denoted as usual by \(\langle \cdot,\cdot\rangle\).  In the sequel, \(A:E^*\rightarrow E\) denotes a bounded linear mapping.

\begin{defi}
	\label{thm:defMonot}
	The linear mapping \(A:E^*\rightarrow E\) is said \emph{monotone} if:
	\[
	\forall x\in E^*,\qquad \bigl\langle Ax  , x \bigr\rangle \geq 0,
	\]
	and is said \emph{strictly monotone} if:
	\[
	\forall x\in E^*\setminus\{0\},\qquad \bigl\langle Ax , x \bigr\rangle > 0.
	\]
\end{defi}

\begin{defi}
	\label{thm:defPseudoMonot}
	The bounded linear mapping \(A:E^*\rightarrow E\) is said \emph{pseudomonotone} or \emph{pseudomonotone in the sense of Brézis}, if all sequence \((x_n)\) in \(E^*\), weakly-* convergent to \(x\in E^*\) and such that:
	\[
	\limsup_{n\rightarrow+\infty} \bigl\langle A x_n,x_n-x\bigr\rangle \leq 0, 
	\]
	has the property:
	\[
	\forall y\in E^*,\qquad \bigl\langle A x,x-y\bigr\rangle \leq \liminf_{n\rightarrow+\infty} \bigl\langle A x_n,x_n-y\bigr\rangle
	\]
\end{defi}

As a particular case of Proposition~23 of \cite{BrezisIeq}, we have:

\begin{prop}
	If the bounded linear mapping \(A:E^*\rightarrow E\) is monotone, then it is pseudomonotone.
\end{prop}

The following theorem is a particular case of a more general result (corollary~29 of \cite{BrezisIeq}) of Brézis.

\begin{theo}
	\label{thm:theoBrez}
	Let \(K\) be a weak-* compact convex subset of \(E^*\) such that \(0\in K\), \(A:E^*\rightarrow E\) be a bounded linear pseudomonotone mapping, and \(\varphi :E^*\rightarrow \mathbb{R}\cup\{+\infty\}\) a convex, weakly-* lower semicontinuous function such that \(\varphi(0)=0\).  We assume:
	\[
	\forall x\in E^*\setminus K,\qquad \langle Ax,x\rangle + \varphi(x) > 0.
	\]
	Then, there exists \(u\in K\) such that:
	\[
	\forall v\in E^*,\qquad \bigl\langle Au, v-u\bigr\rangle + \varphi(v) - \varphi(u) \geq 0.
	\]
	If, in addition, \(A\) is strictly monotone, then \(u\) is unique. 
\end{theo}

Making use of the fact that the bounded closed convex subsets of \(E^*\) are weak-* compact, we have the obvious following corollary (which is not given in \cite{BrezisIeq} where a similar statement for reflexive Banach spaces is given instead).

\begin{coro}
	\label{thm:coroBrez}
	Let \(K\) be a closed convex subset of \(E^*\) that contains \(0\), \(A:E^*\rightarrow E\) be a bounded linear pseudomonotone mapping, and \(\varphi :E^*\rightarrow \mathbb{R}\cup\{+\infty\}\) a convex, weakly-* lower semicontinuous function such that \(\varphi(0)=0\).  We assume:
	\[
	\lim_{\substack{\|x\|_{E^*}\rightarrow +\infty\\
	x\in K}} \frac{\langle Ax,x\rangle+\varphi(x)}{\|x\|_{E^*}} = +\infty.
	\]
	Then, there exists \(u\in K\) such that:
	\[
	\forall v\in K,\qquad \bigl\langle Au, v-u\bigr\rangle + \varphi(v) - \varphi(u) \geq 0.
	\]
	If, in addition, \(A\) is strictly monotone, then \(u\) is unique.
\end{coro}

\noindent\textbf{Proof.}
Let \(K\) and \(\varphi\) be according the hypotheses of this corollary.  There exists \(R>0\) such that:
\[
\|x\|_{E^*} > R\;\text{ and }\;x\in K  \qquad \Rightarrow\qquad \langle Ax,x\rangle +\varphi(x) > 0.
\]
We set:
\[
K' = K \cap B^*(0,R),\qquad \varphi' = \varphi + I_K,
\]
where \(B^*(0,R)\) is the closed ball of radius \(R\) in \(E^*\) and \(I_K\) is the function taking the value \(0\) all over \(K\) and \(+\infty\) outside.  Then, \(K'\) and \(\varphi'\) satisfy the hypotheses of Theorem~\ref{thm:theoBrez}, so that there exists \(u\in K'\subset K\) such that:
\[
\forall v\in E^*,\qquad \bigl\langle Au, v-u\bigr\rangle + I_K(v) + \varphi(v) - \varphi(u) \geq 0,
\]
which yields the claim. \qed

\section*{Acknowledgements}

Flaviana Iurlano thanks the project ``Tremplins nouveaux entrants \& nouvelles entrantes de la Faculté des Sciences et Ingénierie --- Sorbonne Université''.

\end{document}